\documentclass[12pt,twoside,a4paper]{article}
\usepackage[utf8]{inputenc}
\usepackage{blindtext}
\usepackage{soul}
\usepackage{latexsym} 
\usepackage{amsmath,amssymb,amsfonts,amsthm} 
\usepackage{graphics,psfrag} 
\usepackage{fancyhdr,lastpage} 
\usepackage{color,verbatim} 
\usepackage[table]{xcolor} 
\usepackage{metalogo}
\usepackage{etoolbox}
\usepackage{cite}
\usepackage{mathrsfs}
\usepackage[affil-it]{authblk}

\usepackage{amsthm}     
\usepackage{thmtools}   
\usepackage{hyperref} 
\usepackage{cleveref}   
\usepackage{orcidlink}

\usepackage{xfrac}
\usepackage{bbm}
\patchcmd{\thebibliography}{\section*{\refname}}{}{}{}

\declaretheorem[name=Theorem,numberwithin=section]{theorem}
\newtheorem{proposition}{Proposition}[section]
\newtheorem{corollary}{Corollary}[theorem]
\newtheorem{lemma}[theorem]{Lemma}

\newtheorem{definition}[theorem]{Definition}
\newtheorem{example}[theorem]{Example}
\newtheorem{remark}[theorem]{Remark}

\hypersetup{
    colorlinks=true,        
    linkcolor=green,         
    citecolor=magenta,          
    filecolor=magenta,      
}

\usepackage{calc}

\usepackage[a4paper,left=2cm, right=2cm, top=3cm, bottom=3cm]{geometry}
\let\oldbibliography\thebibliography
\renewcommand{\thebibliography}[1]{%
  \oldbibliography{#1}%
  \setlength{\itemsep}{-1.5mm}%
}

\def\R{\mathbb{R}}

\def\N{\mathbb{N}}
\def\P{\mathbb{P}}

\def\E{\mathbb{E}}

\newcommand{\be}{\begin{equation}}
\newcommand{\ee}{\end{equation}}
\newcommand{\bea}{\begin{eqnarray}}
\newcommand{\eea}{\end{eqnarray}}
\newcommand{\beann}{\begin{eqnarray*}}
\newcommand{\eeann}{\end{eqnarray*}}
\newcommand{\benn}{\begin{equation*}}
\newcommand{\eenn}{\end{equation*}}
\newcommand{\DG}{Dawson--G\"artner }

\newcommand{\cK}{{\mathcal K}}  
\newcommand{\cM}{{\mathcal M}}  
\newcommand{\cP}{{\mathcal P}}  
\newcommand{\cS}{{\mathcal S}}  
\newcommand{\cW}{{\mathcal W}}  
\newcommand{\Square}{{S^n_i\times S^n_j}}
\newcommand{\SquareDiag}{{\Pi^n_{diag}}}
\newcommand{\SquareDiagI}{{S^n_i\times S^n_i}}
\newcommand{\SquareUnion}{{(S^n_i\times S^n_j)\cup (S^n_j\times S^n_i)}}

\newcommand{\ER}{Erd\H{o}s-Rényi }

\newcommand{\ie}{i.e.\ }

\newcommand{\Space}{\mathbf{Z}}
\newcommand{\Topo}{\mathcal{O}_\Space}
\newcommand{\Borel}{\mathcal{B}(\Space)}

\newcommand{\signed}{\pm}

\newcommand{\Proba}{\mathcal{P}(\Space)}
\newcommand{\Probam}{\mathcal{M}_1(\Space)}
\newcommand{\SubProba}{\mathcal{M}_{\leq1}(\Space)}
\newcommand{\Meas}{\mathcal{M}_+(\Space)}
\newcommand{\SignedMeas}{\mathcal{M}_{\pm}(\Space)}

\newcommand{\CbFunct}{C_b(\Space)}

\newcommand{\Graphon}{\mathcal{W}_1}

\newcommand{\Kernel}{\mathcal{W}_{\signed}}

\newcommand{\UGraphon}{\widetilde{\mathcal{W}}_1}

\newcommand{\UGraphond}{\widetilde{\mathcal{W}}_{1,d}}

\newcommand{\TotalMass}[1]{\norm{#1}_\infty}
\newcommand{\TM}[1]{\norm{#1}_\infty}

\newcommand{\InvRelabel}{S_{[0,1]}}
\newcommand{\Relabel}{\bar{S}_{[0,1]}}

\newcommand{\NcutR}[1]{\Vert#1\Vert_{\square,\R}}

\newcommand{\NcutRSymbol}{\Vert\cdot\Vert_{\square,\R}}

\newcommand{\simd}{\sim_{d}}

\newcommand{\dd}{\delta_{\square}}

\newcommand{\Ugo}{\mathcal{X}^*}

\renewcommand{\P}{\mathbb{P}}

\newcommand{\drv}{\mathrm{d}}
\newcommand{\rd}{\mathrm{d}}

\newcommand{\norm}[1]{\Vert#1\Vert}

\fancypagestyle{plain}{%
  \fancyhf{}%
  \fancyfoot[C]{\footnotesize Page \thepage\ of \pageref{LastPage}}%
}

\newcommand{\EE}{\mathbb{E}}
\newcommand{\lep}{\left(}
\newcommand{\rip}{\right)}

\newcommand{\keywords}[1]{\noindent\textbf{Keywords:} #1}
\newcommand{\subjclass}[1]{\noindent\textbf{AMS2020 subject classification:} #1}

\usepackage{fancyhdr}
\title{Large deviations for probability graphons}

\usepackage{authblk}

\author[1]{Pierfrancesco Dionigi\,\orcidlink{0000-0003-2180-8669}\thanks{Pierfrancesco Dionigi is supported by the ERC Synergy under Grant No. 810115 - DYNASNET.\\ \texttt{email: dionigi@renyi.hu}}}
\author[2,3]{Giulio Zucal\,\orcidlink{0009-0000-8261-1291}\thanks{Giulio Zucal \\ \texttt{email: zucal@mpi-cbg.de}}}
\affil[1]{HUN-REN Alfréd Rényi Institute of Mathematics, Budapest, Hungary}
\affil[2]{Max Planck Institute of Molecular Cell Biology and Genetics, Dresden, Germany}
\affil[3]{Center for Systems Biology Dresden, Germany}

\date{\today}

\begin{document}
\maketitle

\begin{abstract}
We establish a large deviation principle (LDP) for \emph{probability graphons}, which are symmetric functions from the unit square into the space of probability measures. This notion extends classical graphons and provides a flexible framework for studying the limit behavior of large dense weighted graphs. In particular, our result generalizes the seminal work of Chatterjee and Varadhan (2011), who derived an LDP for Erdős–Rényi random graphs via graphon theory. We move beyond their binary (Bernoulli) setting to encompass arbitrary edge-weight distributions. Specifically, we analyze the distribution on probability graphons induced by random weighted graphs in which edges are sampled independently from a common reference probability measure supported on a compact Polish space. We prove that this distribution satisfies an LDP with a good rate function, expressed as an extension of the Kullback–Leibler divergence between probability graphons and the reference measure. This theorem can also be viewed as a Sanov-type result in the graphon setting. Our work provides a rigorous foundation for analyzing rare events in weighted networks and supports statistical inference in structured random graph models under distributional edge uncertainty.
\end{abstract}

\keywords{Graph limits, Large networks, Probability graphons, Edge-decorated graphs, Weighted graphs, Dense weighted graph sequences, Random matrices, Multiplex networks, Large deviations, Relative entropy, Sanov’s theorem}  

    \vspace{0.2cm}
\subjclass{05C80 (Random Graphs), 60B20 (Random matrices), 60B10 (Convergence of measures), 60F10 (Large deviations), 60C05 (Combinatorial probability), 28A33  (Spaces of measures, convergence of measures)}
\newpage
\section{Introduction}
The past 40 years saw the emergence, rise, and success of complex systems. Retrospectively it is not a surprise that this happened. The last half century witnessed an unprecedented pace of development globally, giving birth to a world profoundly different from the beginning of past century: a connected and unbelievably small globe, fastly evolving and complex to analyze. Networks and their applications became a pivotal tool in analyzing and decrypting some  of new and difficult challenges that emerged in the last years. From social sciences (where the application of networks' topological properties dates back at least to the 60's, see \cite{milgram1967small}) and economics, to biology, neuroscience, epidemiology, energy and communication infrastructures (one for all, internet) their presence is ubiquitous. Even the most recent and successful approaches to complexity via machine learning techniques are deeply intertwined with networks and graph theory (see Graph neural network etc). Together with (and because of) their systematic use, networks put through some of the most interesting challenges in probability and discrete mathematics of the last century. We will mention three of these issues, which are the main motivation of our work, and we will explain how our work naturally answers to part of these interesting mathematical questions.

Firstly, driven by the average size of the networks present in real world applications, mathematicians needed to change perspective thinking about very large graphs and their global behavior, in contrast with the beginning of graph theory when the focus was more on finite size and local questions. Big numbers in math call for limit theories, especially in probability, when the underlying object is random and built by a multitude of smaller parts. The hope is that the limit object is somewhat easier to analyze thanks to cancellations or averaging effects (think about the law of large numbers or central limit theorems). Graph limit theories are a natural answer to the above problem and represented one of the most interesting mathematical developments in this area in the recent past. In this area, the literature naturally divides according to the density of the underlying graph sequence, distinguishing between the dense, sparse, and intermediate regimes. The first two are by now the most well-developed: the theory of dense graph limits is firmly established \cite{BORGS20081801,Lovsz2007SzemerdisLF, LOVASZ2006933,borgs2011convergentAnnals}, and the theory of sparse graph limits has also reached a mature stage \cite{BenjaminiLimit,local-global1,Hatami2014LimitsOL}. Both have been extensively studied and applied since the early 2000s; see also the monograph \cite{LovaszGraphLimits}. By contrast, the search for effective limit theories in the intermediate-density regime remains an active and evolving area of research \cite{MarkovSpaces,KUNSZENTIKOVACS20191,frenkel2018convergence,backhausz2018action,veitch2015classrandomgraphsarising,janson2016graphons,caron2017sparse,borgs2018sparse,borgs2020identifiability,borgs2019sampling,10.1214/18-AOS1778,JANSON2022103549}, where no consensus has yet emerged on a unifying framework.

While the mathematical literature has largely concentrated on binary graphs, weighted networks have been extensively explored in the network science community. Therefore, a fact that has been somewhat overlooked in the mathematical literature is that most real networks carry richer information between two nodes than a mere binary indicator of connection. In practice, many networks are naturally weighted: edge weights may represent quantitative or physical attributes (such as distance, flow volume, or electrical potential difference) or qualitative ones (for instance, edges distinguished by types or colors, as in multiplex networks). This perspective highlights that real-world networks are far from being simple $0/1$ arrays and underscores the relevance of connecting random graph theory with random matrix theory. As a consequence, most of the theoretical tools regarding graph limit theory have been developed systematically only for binary graphs until recently. In fact, even if limit objects for edge-decorated graphs have been interestingly considered from the rise of graph limit theory, see the unpublished work \cite{lovász2010limits}, a complete convergence theory for these objects, parallel to the binary graphs case, has been lacking in the literature. Exceptions are the notable works \cite{falgasravry2016multicolour,rath2011multigraph} limited to the case in which the space of decorations is discrete and some applications of these works \cite{rath2012Configmultigraph,rath2012time}. The connection with exchangable random arrays and Aldous-Hoover theorem \cite{kallenberg1992symmetries, hoover1979relations, aldous1981representations, aldous2010exchangeability,austin2008exchangeable,diaconis2007graph} was also realized (and is also connected to limit theories for more general combinatorial objects as hypergraphs \cite{hypergrELEK20121731,HypergraphsSzegedy2,HypergraphonsZhao,zucal2023action,zucal2024probabilitygraphonspvariablesequivalent}). Recently, the development of a complete limit theory for dense sequences of edge-decorated graphs finally received a lot of attention \cite{abraham2023probabilitygraphons,zucal2024probabilitygraphonsrightconvergence,zucal2024probabilitygraphonspvariablesequivalent,athreya2023pathconvergencemarkovchains,KUNSZENTIKOVACS2022109284}. In some of these papers the authors called these limit objects probability graphons \cite{abraham2023probabilitygraphons,zucal2024probabilitygraphonsrightconvergence,zucal2024probabilitygraphonspvariablesequivalent} to underline the connection with probability measures. These objects have already found applications in diverse areas \cite{dufour2024inferencedecoratedgraphsapplication,andrade2025monochromaticsubgraphsrandomlycolored,ganguly2025meanfieldanalysislatentvariable}. Probability graphons are functions from $[0,1]^2$ to $\mathcal{P}(\Space),$ the space of probability measures over $\Space,$ where $\Space$ is the set of values that the edges can take. They can be thought of as the natural limit objects of weighted graphs. It is indeed possible to think of a random weighted graph to be a random realization of a discrete object where instead of having an array of $0/1$ random variables we have an array of probability distributions. When the number of vertices of the graph is going to infinity this objects can be seen converging to probability graphons. The reader already familiar with the details of graph limit theory for binary graphs might immediately realize the difficulties that an object that takes values in an infinite dimensional space might create. 

Most of the scientific investigation resources in network science were directed to the comprehension of how to build effective and faithful network models of real world ones. That was by no surprise the top-ranked task: good models allow scientists to forecast structure and patterns in existing networks, and predict their possible evolution or behavior. But there is a trade-off between how close to reality a network is and its analytical simplicity. The more we require structure appearing in real life networks (e.g. number of triangles, cliques, or homophily between groups), the more the model becomes complicated and difficult to work with. Probabilistic models are usually prioritizing independence of the edges' appearance, that can't give the correct results when it comes to patterns or higher structures withing the network. An effective approach to this problem is asking how unlikely it is for our simple models to look like the snapshot of what we are able to witness in real world networks. Large deviations are such type of tool that permit us to describe from a simple model how unlikely it is to observe a precise realization with some given characteristic (e.g. a given number of triangles in excess), and how, on average, will look like such network conditioning on this event. The importance of a large deviation principle for large dense binary networks was already highlighted by the work of \cite{ChatterjeeVaradhan2011} (see also the more didactic book \cite{chatterjee2017large}) and the line of research it gave raise to \cite{grebik2023large,Borgs_Chayes_Gaudio_Petti_Sen_2025,lubetzky2015replica,ChatterjeeDiaconisExponential2013,chatterjeeRandomGraphsGiven2011b,bhamidi2008mixing,Den_Hollander_Markering_2023,radinPhaseTransitionsExponential2013b,10.1214/18-EJP135}.

Our work stands at the meeting point of the three challenges mentioned above. We develop an LDP for probability graphons as exposed in the next subsection. The reader should be able to understand our results as exposed below, but is referred to the following sections for the precise definitions and technical statements about probability graphons and our results. We believe that our work can lead to the development of analogous results for weighted graphs that have been already studied for the binary case, closing therefore the gap between the two fields. For instance the development of theory of weighted exponential graphs similarly to what has been done in \cite{WILSON201737,bhamidiWeightedExponentialRandom2018} or to the study of large deviations of evolving weighted graphs, in a similar fashion to \cite{braunsteinsSamplepathLargeDeviation2023,bhamidi2025large}.
\subsection{Results and main theorem}
\paragraph{Weighted graphs.} We now introduce one of the central objects of this paper: weighted graphs. Although their definition is straightforward, we spell it out explicitly in order to fix notation.  
A weighted graph of size $n$ is a triple $g_{n,M} = (V,E,M)$, where $V$ is the vertex set with $|V|=n$, $E = (e_{i,j})_{1 \leq i,j \leq n}$ denotes the edge set, and $M = (m_{i,j})_{1 \leq i,j \leq n}$ assigns a weight to each edge. In many applications, the weights take values in a type space $\Space$. It is common in the literature to consider the complete edge set, leaving the presence or absence of an edge to be encoded by its weight in $M$. Conventions also differ as to whether self-loops $(e_{i,i},m_{i,i})_{1 \leq i \leq n}$ are included. Since our theory is unaffected by this choice, we adopt the customary convention of excluding self-loops.  
It is natural to consider a random graph model where the edge weights are sampled independently from a distribution $\nu \in \mathcal{P}(\Space)$. This induces a natural distribution $\mu_{n,\nu}$ on the space $\mathcal{G}_n$ of weighted graphs of size $n$. A random graph drawn from $\mu_{n,\nu}$ will be denoted by $g_{n,\nu} \equiv g_{n,M} = \bigl(V,E,(m_{i,j} \sim \nu)_{1 \leq i,j \leq n}\bigr)$.
\begin{remark}\label{rmk:Zero}
When dealing with $g_{n,\nu}$, in the case in which we want to consider the complete edge set, we are leaving the presence or absence of an edge to be encoded by its weight in $(m_{i,j}\sim \nu)_{1 \leq i,j \leq n}$. Therefore, one should make sure that $0$ is an element of $\Space$, or alternatively considering $\Space\cup\{0\}$ when it makes sense. For Polish spaces $\Space$ where $\Space\cup\{0\}$ does not make sense, it is customary to introduce a cemetery state $\partial$, see \cite[Remark 6.2]{abraham2023probabilitygraphons}. In the following, without loss of generality, we will assume that the state $0$ is included in $\Space$.
\end{remark}

\paragraph{What is a probability graphon?} Before stating our main result without going into the technical details of the following Sections, it will be helpfull for the reader to have an intuition of what is a probability graphon. Recall that \emph{real graphon} is a function $w:[0,1]^2\to [0,1]$. The power of real graphon is that they are the dense limit of binary graphs: for each $(x,y)$, $w(x,y)$ represents the density of connection between the vertices labeled $x$ and $y$. Now, what if we have a weighted graph, with a possibly uncountable type space? What is the point $(x,y)$ representing? It turns out that the best way to represent the limit of such objects is to attach to each point $(x,y)$ the limit distribution of the weights associated to the edge between $x$ and $y$.

Let $\Space$ be a compact Polish space and $\mathcal{P}(\Space)$ the space of probability measures on $\Space$. Let $\nu\in\mathcal{P}(\Space)$ be a probability measure and $\widetilde{W}$ an \emph{unlabeled probability graphon}, i.e. measurable function from $[0,1]\times[0,1]$ to $\mathcal{P}(\Space)$ (see Definition \ref{DefProbGraphon} and Section \ref{subsection_unlabeled_cut_distance} to understand why we need the \emph{unlabeled} version; in general we indicate with $\,\widetilde{}\,$ always elements that are \emph{unlabeled} to distinguish them from their labeled counterpart). In \cite{abraham2023probabilitygraphons} it was showed that the natural space that contains probability graphons $\widetilde{\mathcal{W}}_1(\Space)$, has nice features and can be equipped, similarly with the space of real graphons, with a cut metric such that the space is compact (see details in Section \ref{Sec2ProbGraphons}). It will be convenient to view a weighted graph as an edge-decorated graph whose edge values are probability measures, that is, as a function 
\[
g_{n,M} : V^2 \to \mathcal{P}(\Space).
\]
This representation is obtained by associating to each weight $m_{i,j} \in \Space$ the Dirac measure $\delta_{m_{i,j}} \in \mathcal{P}(\Space)$ centered in $m_{i,j}$. In this way, every weighted graph $g_{n,M}$ can be naturally identified with a probability graphon $\widetilde{W}_{g_{n,M}}$, after partitioning $[0,1]^2$ into $n^2$ equal squares $\Square$ (see Section~\ref{subsec:stepfunctions} for the rigorous definition of these squares) of size $\tfrac{1}{n^2}$ and associating to every square $\Square$ the value $\delta_{m_{i,j}}$
so as to map the edge set onto it.
\paragraph{Our main results.} Consider the family of distributions $\{\mu_{n,\nu}\}_{n \in \mathbb{N}}$, each inducing the random weighted graph model $g_{n,\nu}$ on $\mathcal{G}_n$. As discussed above, every weighted graph can be identified with a probability graphon. Denote by $\widetilde{W}_{g_{n,\nu}}$ the unlabeled random probability graphon model associated with sampling a weighted graph of size $n$ from $\mu_{n,\nu}$. 
From a classical perspective, one expects that repeated sampling of the graphs $g_{n,\nu}$ yields, for each box in $[0,1]^2$ corresponding to a pair $(i,j)$, an empirical measure formed by the $\delta_{m_{i,j}}$’s which converges to the underlying distribution $\nu$. This is directly analogous to the law of large numbers for i.i.d.\ samples. Large deviations from this behavior are described by Sanov’s theorem, which provides an LDP for the empirical distribution of i.i.d.\ samples with rate function given by the relative entropy with respect to the sampling law (see, e.g., \cite[Theorem~6.2.10]{dembo2009large}).  
Our main result can be interpreted as a Sanov-type generalization to the space of probability graphons $\widetilde{\mathcal{W}}_1(\Space)$. In this setting, the role of the empirical measure is played by the empirical distribution of edge weights across the graph, encoded in $\widetilde{W}_{g_{n,\nu}}$. Unlike the classical case, where the empirical measure is defined through a global averaging factor $1/n$, here the natural averaging procedure is local: it is expressed through integrals of $\widetilde{W}_{g_{n,\nu}}$ over subsets $U \subset [0,1]^2$, thereby capturing the distributional structure of weights across the network.
To state formally our results, let $\mu_{n,\nu}$ be the distribution of the weighted graphs of size $n$ with weights sampled from $\nu$: this automatically defines a measure  $\mu_{n,\nu}$ in the space $\widetilde{\mathcal{W}}_1(\Space)$. Indeed here and in the following it will be crucial interpreting weighted graphs as probability graphons, and, viceversa, probability graphons with a block structure to realizations or empirical distributions of weighted graphs. The next step is to find what functional gives the right variational problem that identifies the rate function, i.e. the function that describes the exponential cost we have to pay to deviate from the \emph{expected} behaviour. For a probability measure $\mu\in\mathcal{P}(\Space)$ absolutely continuous with respect to $\nu$, $\mu\ll\nu$, we recall the classic relative entropy (or Kullback-Leibler divergence):
\[
\mathcal{H}(\,\mu\mid \nu\,) = 
\int \log \left(\frac{d\mu}{d\nu}\right) \, d\mu.
\]
We can define the \emph{probability graphon relative entropy} of $\widetilde{W}$ with respect to $\nu$ the classic Kullback-Leibler divergence applied to the probability graphon $\widetilde{W}:$
\[
\widetilde{\mathcal{H}}(\widetilde{W} \mid \nu) = \int_{[0,1]^2} \mathcal{H}(\widetilde{W}(x,y) \mid \nu) \, dx \, dy.
\]
We will interpret $\widetilde{\mathcal{H}}(\widetilde{W} \mid \nu)\equiv I_\nu\lep\widetilde{W}\rip$ as a rate function, similarly for what happens in Sanov's theorem. 
\begin{remark}
    Along all of this work we will assume without loss of generality that $\nu$ has as support all of $\Space$.
\end{remark} 
Our main result is the following

\begin{restatable}{theorem}{MainTheorem}[LDP for probability graphons.]\label{thm_mainresult}
    Let $\Space$ be a compact Polish space. Consider the space $\UGraphon(\Space)$ of probility graphons equipped with the unlabeled cut metric $\delta_\square$ defined in Section \ref{subsection_unlabeled_cut_distance}. Then for any closed set $\widetilde{\mathcal{C}} \subseteq \widetilde{\mathcal{W}}_1(\Space)$ and any open set $\widetilde{\mathcal{O}} \subseteq \widetilde{\mathcal{W}}_1(\Space)$ we have
\[
-\inf_{\widetilde{W} \in \widetilde{\mathcal{O}}} I_\nu(\widetilde{W})\leq\liminf_{n \to \infty} \frac{2}{n^2} \log \widetilde{\mu}_{n,\nu}(\widetilde{\mathcal{O}}\,)\leq\limsup_{n \to \infty} \frac{2}{n^2} \log \widetilde{\mu}_{n,\nu}(\widetilde{\mathcal{C}}\,) \leq -\inf_{\widetilde{W} \in \widetilde{\mathcal{C}}} I_\nu(\widetilde{W}) \,.
\]
\end{restatable}
The above result can be interpreted as follows. Suppose we have a random weighted graph model we sample from. Each of the edges has a weight distributed according to $\nu.$ As $n$ grows, $\mu_{n,\nu}$ can be interpreted as the empirical distribution of the sample of weighted graphs we have obtained. Suppose we have a very rare event that we want to study, the event $\widetilde{\mathcal{B}}$ (which in the result above is described for technical reasons through the use of open and closed sets in the space of probability graphons $\UGraphon(\Space)$). The logarithm of the probability of $\widetilde{\mathcal{B}}$ is asymptotically controlled by $-I_\nu(\widetilde{W})$, up to matching upper and lower bounds, where $\widetilde{W}$ is chosen in $\widetilde{\mathcal{B}}$ in order to minimize $I_\nu(\widetilde{W})$ (i.e. maximize the probability). This is a classic mantra in large deviations, where rare events happen in the most probable between all the unlikely ways. It is easy to see how the above result can also be useful in hypothesis testing, when we have a guess for $\nu$, and we measure $\mu_{n,\nu}$ through sampling and we can ask how our prior distribution describes well the rare events of the empirical one $\mu_{n,\nu}$.

One of the most striking and useful features of large deviations is that once we have an LDP is particularly straightforward to study the measure conditioned on rare events. A situation where this can turn out very convenient, in the storytelling of above, is when we see from our sampled data a substantial deviation of weights of triangles from what we would expect given $\nu$. We could therefore condition on this rare event, tilting the distribution and asking how does the \emph{typical} weighted graph looks like under this conditioning. It is not hard to realize that the typical weighted graph is represented by the probability graphon (or graphons) minimizing $I_\nu$ on the rare event $\widetilde{\mathcal{B}}$. The main result of Section \ref{section:CondDist} is about how the natural distance $\delta_\Box$ (the cut distance see Definition \ref{def:ddcut}) of the probability graphon space decreases between the minimizer (or the minimizers) of $I_\nu$ in $\widetilde{\mathcal{B}}$ and instances of a random weighted graph of size $n$ and weights distributed as $\nu$ when interpreted as elements of $\UGraphon(\Space)$. 

\begin{restatable}{theorem}{SecondTheorem}\label{thm:CondDist}
 Let $\Space$ be a compact Polish space. Take any $\nu \in \mathcal{P}(\Space)$ and $n \geq 1$. Let $g_{n,\nu}$ be a random weighted graph of size $n$ distributed according to $\mu_{n,\nu}$, and let $\widetilde{W}_{g_{n,\nu}}$ its associated probability graphon. Let $\widetilde{\mathcal{B}}$ be a closed subset of $\UGraphon(\Space)$ satisfying equation \eqref{eq:closedsubset}. Let $\mathscr{M}_{\widetilde{\mathcal{B}}}$ be the subset of $\widetilde{\mathcal{B}}$ where $I_\nu$ is minimized. Then $\mathscr{M}_{\widetilde{\mathcal{B}}}$ is non-empty and compact, and for each $n$, and each $\varepsilon > 0$,
\[
\mathbb{P}\lep\left.\delta_{\square}(\widetilde{W}_{g_{n,\nu}}\,, \mathscr{M}_{\widetilde{\mathcal{B}}}\,) \geq \varepsilon \,\right|\, \widetilde{W}_{g_{n,\nu}} \in \widetilde{\mathcal{B}}\,\rip \leq e^{-C(\varepsilon, \widetilde{\mathcal{B}}\,)\, n^2}
\]
\textit{where $C(\varepsilon, \widetilde{\mathcal{B}}\,)$ is a positive constant depending only on $\varepsilon$ and $\widetilde{\mathcal{B}}$ and $\delta_{\square}(\widetilde{W}^w\,, \mathscr{M}_{\widetilde{\mathcal{B}}}\,)$ is defined as in \eqref{eq:cutmetricsets}. In particular, if $\mathscr{M}_{\widetilde{\mathcal{B}}}$ contains only one element $\widetilde{W}^*$, then the conditional distribution of $\widetilde{W}_{g_{n,\nu}}$ given $\widetilde{W}_{g_{n,\nu}}\in \widetilde{\mathcal{B}}$ converges to the point mass at $\widetilde{W}^*$ as $n \to \infty$.}
\end{restatable}
The theorem substantially says that once a rare event happens, the typical random weighted graph interpreted as a probability graphon sampled through  $\mu_{n,\nu}$ will asymptotically converge in a neighbor of one of the minimizers of $I_\nu$ in $\widetilde{\mathcal{B}}$.

\begin{remark}
    Many of the results presented in this work can be generalized to vertex-weighted graphs but we make the choice of presenting the results in the setting of edge-weighted only graphs, discussing the extension to vertex-weighted graphs to some remarks, see Remarks \ref{Rmk:Vertex-Weights1}, \ref{rmk:VertexWeights2}, \ref{rmk:VertexWeightedGraphs3} and \ref{rmk:VertesWeights4}. This choice is made in order not to make the presentation not too technical and the notation not too heavy. 
\end{remark}

\begin{remark}
\label{Rmk:WhyCompact}
Our results are stated for $\Space$ compact. Nevertheless, many distributions of high interest in applications have unbounded support. We expect our results to hold for general Polish spaces (possibly with some additional mild conditions) as Sanov's theorem holds for general (non-compact) Polish spaces. However, to extend our proofs to the non-compact case would require additional technical work that we will develop in future work, also to avoid becoming even more technical here.  
\end{remark}

\subsection{Sketch of the proof and structure of the paper.}
The technical difficulty of proving the above results relies in different points. Despite similar conclusions hold for real graphons through the seminal work of Chatterjee and Varadan \cite{ChatterjeeVaradhan2011}, the topological difficulties of the space of probability graphons pose a bigger challange. Probability graphons are functions that take values in an infinite dimensional space whose metric structure is complicated. The key point of our proof, in the same fashion of the real graphon's one, relies in the compactness induced by the cut metric defined on the probability graphon space (see Remarks \ref{rmk:metrics} and \ref{rmk:metrics_2}). We will talk about the technicalities of graphon space in Section \ref{Sec2ProbGraphons}. 

Besides the technical difficulties coming with the space, the second step was identifying with sufficient generality the rate function for our problem. Indeed while the formula of \cite{ChatterjeeVaradhan2011} has a clear interpretation for the real graphon case, we had to find a formula that made sense for general distributions. We found this to be the relative entropy between the target distribution $\nu$ and the model distribution $\mu_{n,\nu}$. In section \ref{Sec:relent}, we define relative entropy, highlighting its role in large deviations, and then explain how to generalize in the probability graphon space $\UGraphon(\Space)$. 

The proof of our results relies upon the argument used by \cite{ChatterjeeVaradhan2011}, carefully generalized and adapted to our setting. Indeed, it is immediately clear that one can hope to use the original argument under the case where $\Space$ is a finite space. Our tactic is a classical \emph{divide and conquer}: we first tackle the finite case generalizing the existing results and then we use Dawson--G\"artner theorem to achieve our result for compact spaces. In section \ref{Sec:LDP_finite} we prove Theorem \ref{thm_mainresult} under the assumption that $\Space$ is finite. The proof goes through 3 different steps. First an upper bound in the weak topology is proved, then the proper upper bound is derived from this. In both steps we heavily use the continuity property of the space, approximating probability graphons with block approximants, and controlling the error of the approximation. Once this is done, we can use standard large deviations results adapted to our case. After, a lower bound is obtained through a standard tilting of measure argument.  

The next step is extending the proof to not finite $\Space$. This is done in Section \ref{sec:dawsongartner} through Dawson--G\"artner theorem. 

Finally in Section \ref{section:CondDist} we prove Theorem \ref{thm:CondDist}.

\vspace{15pt}
\textbf{Acknowledgements:}  The authors would like to thank Luisa Andreis for helpful discussions about the Dawson--G\"artner theorem and Shankar Bhamidi for pointing out useful references about weighted random graphs. 

\section{Probability graphons}\label{Sec2ProbGraphons}
In this section we introduce probability graphons. Our point of view will be general enough for the reader to appreciate the technicalities of the field, but we will limit ourselves to the results needed for our paper. The reader interested to the details of the theory and its application to, e.g., graph sampling is referred to \cite{abraham2023probabilitygraphons}. The reader that is instead familiar to the topic can skip this section or consult it for fixing the notation.

The first thing that needs to be noted is that the space of probability measures with support $\Space$, $\mathcal{P}(\Space)$ is a complex space. The topological and metric structure of $\mathcal{P}(\Space)$ heavily depend on the finiteness, compactness and metric structure chosen on $\Space$. We will therefore start with an \emph{intermezzo} on measure and probability spaces in Section \ref{sec:measurespaces}. We then continue with \ref{sec:realvaluedgraphon} to revise parent theory of real graphons, i.e. the graph limits for binary graphs. From Section \ref{section_def_graphons} we start the journey with probability graphons with their definition, we then continue with the natural distance in this space in Section \ref{section_def_dcut} and the unlabeled space, fundamental for our results, in Section \ref{subsection_unlabeled_cut_distance}.  We then conclude in Section \ref{sec:reglemma} with the statement of (one version of) regularity lemma, a fundamental results for our main theorem.

\subsection{A measure space intermezzo}\label{sec:measurespaces}
In this section, we introduce basic measure-theoretic concepts and fix the notation following \cite{abraham2023probabilitygraphons}. 

A function $\varphi : \Omega_1 \to \Omega_2$ between two probability spaces $(\Omega_i, \mathcal{A}_i, \P_i)$, $i = 1, 2$, is called measure-preserving if it is measurable and if for every set $A \in \mathcal{F}_2$ the condition $\P_2(A) = \P_1(\varphi^{-1}(A))$ holds. This is equivalent to saying that for any measurable non-negative function $f : \Omega_2 \to \mathbb{R}$, the following equation is satisfied:
\begin{equation}\label{eq:InvBijInt}
    \int_{\Omega_1} f(\varphi(x))\, d\P_1(x) = \int_{\Omega_2} f(x)\, d\P_2(x).
\end{equation}
In this work, we will always consider the unit interval $[0,1]$ and the unit square $[0,1]^2=[0,1]\times [0,1]$ equipped with the Euclidean norm, the Borel $\sigma-$algebra and the Lebesgue measure. We will denote the Lebesgue measure on $[0,1]$ with $\lambda$ and the Lebesgue measure on $[0,1]^2$ with $\text{Leb}.$ These are probability spaces.
We denote by $\InvRelabel$ the set of all bijective measure-preserving maps from $[0,1]$ to itself, while $\Relabel$ denotes the set of measure-preserving maps from $[0,1]$ to itself.
\medskip

Let $(\Space,d_{\Space})$ be a non-empty Polish metric space. Let $\Topo$ be the topology generated by the metric $d_{\Space}$ and $\Borel$ be the Borel $\sigma$-algebra on $\Space$ generated by the topology $\Topo$. Furthermore, we let $\CbFunct$ be the  space of continuous bounded real-valued functions on $(\Space,\Topo)$ equipped with the supremum norm $\|\cdot\|_{\infty}.$ We denote by $\SignedMeas$ the space of finite signed measures on $(\Space,\Borel)$, $\Meas$ as the subspace of measures, $\SubProba$ as the subspace of measures with total mass at most $1$, and $\Proba$ (or, equivalently by notation $\Probam$) as the subspace of probability measures. The following chain of inclusions holds:
\[
  \Proba \subset \SubProba \subset \Meas \subset \SignedMeas .
\]

For two measures $\mu,\nu\in  \Meas$ we write $\mu \ll\nu$ to say that $\mu$ is absolutely continuous w.r.t $\nu.$ If $\mu\ll\nu,$ we denote with $\frac{d\mu}{d\nu}$ the Radon-Nikodyn derivative of $\mu$ w.r.t\ $\nu.$ For a signed measure $\mu\in\SignedMeas$, we recall the unique Hahn-Jordan decomposition $\mu  = \mu^+ -  \mu^-$ where $\mu^+, \mu^-  \in\Meas$ are mutually singular measures (that  is $\mu^+(A)=0$ and $\mu^-(A^c)=0$ for some measurable set $A$).
The total variation measure $\vert\mu\vert$ of $\mu\in\SignedMeas$ is defined as $\vert\mu\vert  =   \mu^+  +  \mu^-  \in  \Meas$. In particular, for a measure $\mu\in\Meas$, we have $\vert  \mu \vert  = \mu$. For a signed measure $\mu\in \SignedMeas$ and a real-valued measurable function $f$ defined on $\Space$, we denote by $\mu[f]=\mu(f)=\langle \mu, f \rangle=
\int f \, \rd \mu=\int_{\Space} f(x)\, \rd\mu$ the integral of $f$ with respect to $\mu$ when well defined. For a signed measure  $\mu\in\SignedMeas$, the total mass of $\mu$ is $\|{\mu}\|_{TV}=\mu^+(\Space) + \mu^-(\Space),$ which also equals the supremum of $\mu(f)$ taken over all measurable functions $f$ with values in $[-1,1]$.

We consider the space $\SignedMeas$ equipped with the weak convergence topology, i.e. the smallest topology such that for each $f\in\CbFunct$ the map $\mu \mapsto \mu(f)$ is continuous. In particular, a sequence  of signed measures $(\mu_n)_{n\in\N}$  weakly converges to some $\mu\in\SignedMeas$ if and only if, for each function $f\in   \CbFunct$, we have $\lim_{n\to +\infty} \mu_n(f) = \mu(f)$. Recall that $\Meas$ and $\Proba$ equipped with the weak convergence topology are Polish spaces. 

The maps $\mu   \mapsto  \mu^+$ and $\mu   \mapsto  \mu^-$ (and thus also $\mu \mapsto \vert\mu\vert$) are measurable (see \cite[Theorem~2.8]{MeasurableSetsMeasures} and Remark 2.4 in \cite{abraham2023probabilitygraphons}). As a consequence, the map $\mu \mapsto \|\mu\|_{TV}$ is also measurable.
Observe that $\Proba$ and $\Meas$ are closed, and therefore measurable subsets of $\SignedMeas.$

The interested reader can find more details about the weak convergence of signed measures in Section 2 in \cite{abraham2023probabilitygraphons} and in the standard references \cite{Bogachev,BogachevMT1,BogachevMT2,Ethier,MeasurableSetsMeasures,Varadarajan}.
A possible metric that has the important property of metrizing weak convergence on the space of measures is the following:
 
\begin{definition}[Lévy-Prokhorov metric]\label{LevyProk}
 The \emph{Lévy-Prokhorov Metric} $d_{\mathcal{LP}}$ on the space of measures $\Meas$ is for $\eta_1,\eta_2\in \Meas$
$$\begin{aligned}
d_{\mathcal{LP}}\left(\eta_{1}, \eta_{2}\right)=&\inf \left\{\varepsilon>0: \eta_{1}(U) \leq \eta_{2}\left(U^{\varepsilon}\right)+\varepsilon \text{ and } \right.\\
&\left.\eta_{2}(U) \leq \eta_{1}\left(U^{\varepsilon}\right)+\varepsilon  \text{ for all } U \in \Topo\right\},
\end{aligned}$$

where we recall that $\Topo$ is the Borel $\sigma$-algebra on $\Space$ generated by $d_{\Space}$, $U^{\varepsilon}$ is the set of points that have a distance $d_{\Space}$ smaller than $\varepsilon$ from $U$.
\end{definition}

\begin{remark}\label{rmk:metrics}
Besides the Lévy--Prokhorov distance, other classical metrics can be used to metrize the weak topology on the space of probability measures. 
Typical examples are the Kantorovich--Rubinstein norm (also known as the bounded Lipschitz distance), the Fortet--Mourier norm, and norms built from convergence determining sequences of bounded continuous functions. See, for example, \cite{abraham2023probabilitygraphons}.
\end{remark} 
\subsection{Real-valued graphons}\label{sec:realvaluedgraphon}

Before dwelling upon probability graphons, we briefly recall here the notion of (real-valued) graphons from dense graph limit theory \cite{LovaszGraphLimits,LOVASZ2006933,BORGS20081801,borgs2011convergentAnnals}. Being probability graphons a generalization of real-valued graphons, we briefly go through the theory so that unfamiliar readers are introduced to both theories and can appreciate the step forward that probability graphons make. The reader familiar with the theory can directly skip this section.
\begin{definition}[Real-valued graphon]
   Let $\Omega_0$ be a probability space. A \emph{real-valued graphon} is a measurable function
    $$w:\Omega_0\times \Omega_0\rightarrow [0,1].$$
\end{definition}

\begin{remark}
Sometimes, one requires real-valued graphons to be symmetric in the two variables or $\Omega_0=[0,1]$. We do not impose this condition a priori here. 
\end{remark}

Let's suppose now that $\Omega_0=[0,1].$

We denote by $\NcutRSymbol$ the real-valued cut norm defined for (linear combinations of) real-valued graphons as:
\begin{equation}
	\label{eq_def_NcutR}
\NcutR{w} = \sup_{S,T \subset [0,1]} \left\vert \int_{S\times T} w(x,y)\ \drv x\drv y \right\vert
\end{equation}
where $w$ is a (linear combination) of real-valued graphons, see  \cite[Chapter 8]{LovaszGraphLimits}. This induces a distance that we call $d_{\Box,\R}$.

A key point of the theory of graphons is that the space can be made compact defining equivalence classes of graphons.  

\begin{definition}[Relabeling action on graphons]
For $w:[0,1]^2\to[0,1]$ measurable and $\varphi\in\InvRelabel$ a measure preserving bijection, define the \emph{relabeling} (or pullback)
\[
w^{\varphi}(x,y)\;:=\; w\big(\varphi(x),\varphi(y)\big),\qquad (x,y)\in[0,1]^2.
\]
This defines a right action of $\InvRelabel$ on the space $\mathcal{W}_0$ of (not necessarily symmetric)
$[0,1]$-valued graphons on $[0,1]^2$, with $(w^{\psi})^{\varphi}=w^{\psi\circ\varphi}$ and $w^{\mathrm{id}}=w$.
\end{definition}
From the above we can define the pseudo-metric $\delta_{\Box,\R}$:
$$\delta_{\square,\R}(u,w)=\inf_{\varphi\in \InvRelabel}\NcutR{u-w^{\varphi}},$$
which turns out to be a central object of the theory with nice properties.

\begin{remark}[Cut norm invariance under relabeling]
For any $u,w\in\mathcal{W}_0$ and $\varphi,\psi\in\InvRelabel$,
\[
\big\|u^{\psi}-w^{\varphi}\big\|_{\square,\R}=\big\|u-w^{\varphi\circ\psi^{-1}}\big\|_{\square,\R},
\]
hence in particular $\|w^{\varphi}\|_{\square,\R}=\|w\|_{\square,\R}$ and
$\delta_{\square,\R}(u^{\psi},w^{\varphi})=\delta_{\square,\R}(u,w)$.
\end{remark}

\begin{definition}[Quotient by relabelings and weak isomorphism]\label{def:quotient}
Define an equivalence relation on $\mathcal{W}_0$ by
\[
u\sim w \quad \Longleftrightarrow \quad \delta_{\square,\R}(u,w)=0.
\]
(The relation $\sim$ is often called \emph{weak isomorphism}.) The \emph{quotient space of graphons}
is the set of equivalence classes
\[
\widetilde{\mathcal{W}}_0\;:=\;\mathcal{W}_0/\!\sim \;=\;\{\widetilde{w}: w\in\mathcal{W}_0\}.
\]
We endow $\widetilde{\mathcal{W}}_0$ with the metric induced by the cut distance:
for representatives $u,w$,
\[
\delta_{\square,\R}\big(\widetilde{u},\widetilde{w}\big)\;:=\;\tilde{\delta}_{\square,\R}(u,w).
\]
This is well defined because $\delta_{\square,\R}$ is invariant under relabelings and vanishes
exactly on equivalence classes.
\end{definition}

\begin{lemma}[Pseudometric vs metric]
On $\mathcal{W}_0$, $\delta_{\square,\R}$ is a pseudometric. On the quotient
$\widetilde{\mathcal{W}}_0$ of Definition~\ref{def:quotient}, $\delta_{\square,\R}$ is a genuine
metric.
\end{lemma}

\begin{remark}[Symmetric graphons]\label{rem:symmetric_subspace}
If desired, one may restrict to the subspace $\mathcal{W}_0^{\mathrm{sym}}$ of symmetric graphons
($w(x,y)=w(y,x)$ a.e.). All definitions above and the compactness result below remain valid on the
corresponding quotient $\widetilde{\mathcal{W}}_0^{\mathrm{sym}}$.
\end{remark}

\begin{theorem}[Compactness in the cut metric]\label{thm:compactness}
The metric space $\big(\widetilde{\mathcal{W}}_0,\delta_{\square,\R}\big)$ is compact. Equivalently,
every sequence $(w_n)_{n\ge 1}$ of graphons admits a subsequence $(w_{n_k})_{k\ge 1}$ and
relabelings $\varphi_k\in\InvRelabel$ such that, for some graphon $w$,
\[
\lim_{k\to\infty}\big\|\,w_{n_k}^{\varphi_k}-w\,\big\|_{\square,\R}=0.
\]
In particular, $\widetilde{\mathcal{W}}_0$ is complete and totally bounded.
\end{theorem}

The proof of this combines a weak regularity lemma (approximating graphons by bounded-complexity step functions in the
cut norm) with a compactness argument on the finite-dimensional simplex of step functions modulo
relabelings, followed by a diagonal extraction.

\begin{remark}[Normalization and ranges]\label{rem:graphonrange}
The statements above hold verbatim for graphons with values in a compact interval $[-K,K]$ with the
obvious modification of the cut norm. The theory can be as well developed for general probability spaces $\Omega_1$, see \cite{jansonGraphonsCutNorm2013}. We fixed $[0,1]$-valued graphons here for convenience.
\end{remark}

See \cite{LovaszGraphLimits} for more details about the real-valued graphons, the real-valued cut norm, the real-valued cut metric and the compact space contruction.

\subsection{Definition of probability graphons}\label{section_def_graphons}

We can now introduce the theory of probability graphons as developed in \cite{abraham2023probabilitygraphons}. Other references are \cite{zucal2024probabilitygraphonsrightconvergence,KUNSZENTIKOVACS2022109284,lovász2010limits,rath2011multigraph}.
Following the notation of the previous subsections, for a measurable set $\Space\subset \R$, we denote with $\mathcal{P}(\Space)$ the space of probability measures on $\Space$.

\begin{definition}[Probability graphon]\label{DefProbGraphon}
Let $\Space\subset \R$ be a Borel measurable subset of $\R$ and $\Omega_1$ be a probability space. A \emph{probability graphon} on $\Space$ (or $\mathcal{P}(\Space)$-valued kernel) is a map $W$ from $\Omega_1\times \Omega_1$ to $\mathcal{P}(\Space)$, such that:
\begin{itemize}
    \item $ W$ is a probability measure in $\mathrm{d} z$: for every $(x_1, x_2) \in \Omega_1^2, \ W(x_1, x_2 ; \cdot)$ belongs to $\mathcal{P}(\Space)$.
    \item $ W$ is measurable in $(x_1, x_2)$: for every measurable set $E \subset \Space$, the function $(x_1, x_2) \mapsto$ $ W(x_1, x_2 ; E)$ defined on $\Omega_1\times \Omega_1$ is measurable.
\end{itemize} 
\end{definition}
As before is of particular interest (and generality) the case with $\Omega_1=[0,1]$.  

\begin{remark}[Probability graphons $ W : \Omega_1\times \Omega_1 \to \cP(\Space)$]As for Remark \ref{rem:graphonrange} the choice of $\Omega_1=[0,1]$ endowed with the Lebesgue measure $\lambda$ is a matter of convenience: (see Definition 3.1 in \cite{abraham2023probabilitygraphons} or Definition 3.3 in \cite{zucal2024probabilitygraphonsrightconvergence}). However, the theory can be developed for the case of a general probability space $\Omega_1$, see Remark 3.4 in \cite{abraham2023probabilitygraphons}. 
\end{remark}

We will focus on this case and we indicate with $\Graphon(\Space)$ the space of probability graphons, i.e.\ probability graphons from $[0,1]\times [0,1]$ to $\mathcal{P}(\Space)$, where we identify probability graphons that are  equal almost everywhere on $[0,1]^2$ (with respect to the Lebesgue measure). 

We will sometimes need to consider subsets of the space of probability graphons $\Graphon(\Space)$. For a generic subset $\cM\subset \cP(\Space)$, we  denote  by $\cW  _\cM$  the subset  of probability graphons  $W\in  \Graphon(\Space)$  which are  $\cM$-valued: $W(x,y; \cdot)\in \cM$ for every $(x,y)\in [0, 1]^2$.

 One can easily identify real graphons with the space $\mathcal{W}_{\mathcal{P}(\{0,1\})}$:

\begin{example}[On real-valued kernels]\label{rem:real-valued-kernels}
Every real-valued  graphon  $w$ can be represented as a probability graphon $W$ in the following way. Let's consider $\Space =\{ 0, 1\}$ with the discrete metric and the probability graphon $W$ defined as $W(x,y;\drv z)  = w(x,y)  \delta_1(\drv z) +  (1-w(x,y)) \delta_0(\drv z)$ for every $x,y\in[0,1]$, where we recall that $\delta_s$ is the Dirac mass located at $s\in \R$. In particular, we have $$w(x,y)=W(x,y; \{1\})=\int_{\Space}z W(x,y,\drv z)$$ for $x,y \in [0, 1]$. 
\end{example}
\begin{example}[Finitely edge colored graphs]\label{rem:finitecoloredgraphs}
 For every $x,y\in[0,1]$, let $\Space=\{0,1,\ldots,r\}$ with the discrete metric we can define the probability graphon: 
$$
W(x,y;\mathrm{d}z)  =w_{r}(x,y)\delta_r(\mathrm{d} z)+\ldots +w_1(x,y)  \delta_1(\mathrm{d} z) +  w_0(x,y) \delta_0(\mathrm{d} z)$$ where $w_i(x,y)\geq 0$ and $\sum^r_{i=0}w_i(x,y)=1$ for every $x,y\in[0,1].$
As above, integrating away the variables one by one we can obtain an overlap of $r$ distinct real-graphons.
\end{example}
\begin{definition}\label{DefSymProbGr}
We say that a measure-valued kernel or graphon $W$ is symmetric 
if for almost every $x,y \in [0,1]$, $W(x,y;\cdot)=W(y,x;\cdot)$.
\end{definition}

\begin{remark}[Symmetric kernels]
We will consider in general non-symmetric probability graphons to handle also directed graphs whose
adjacency matrices are thus \emph{a priori} non-symmetric.
\end{remark}
For a probability graphon $W$ we can then define its action on any bounded function $f\in \CbFunct$. We denote by $W[f]$ the real-valued graphon defined by 
$$
W[f](x,y)=W(x,y;f)=\int_{\Space}f(z)W(x,y,\mathrm{d}z)
$$
for each $x,y\in [0,1].$

We conclude this subsection regarding extensions for vertex-weighted graphs.

\begin{remark}[Extension to vertex-weights]
\label{Rmk:Vertex-Weights1}The framework presented here for probability-graphons could easily be extended to add weights ($\Space^{\prime}$-valued decorations) on the vertices.
In this case, graph limits for graphs with $\Space-$valued decorations on the edges and $\Space^{\prime}$-valued decorations on the vertices are represented by a probability-graphon $W^{\text{e}}:[0,1]^2 \to \mathcal{P}(\Space)$ for edge-weights as before, and a one-variable kernel $W^{\text{v}} : [0,1] \to \mathcal{P}(\Space^{\prime})$ for vertex-weights. See also \cite[Remark 1.7]{abraham2023probabilitygraphons}.
\end{remark}

\subsection{The cut distance}\label{section_def_dcut}
Analogously to the cut norm for real-valued graphons, we can introduce the \emph{cut  distance}, a distance function on $\Graphon(\Space)$.
For a probability graphon $W\in\Graphon(\Space)$ and a measurable subset $A\subset [0,1]^2$, we define $W(A;\cdot),$ the measure on $\Space\subset\R$
\[
  W(A;\cdot) = \int_{A} W(x,y;\cdot)\ \drv x \drv y.
\]
We can then define  the following semidistance on the space of probability graphons, which will be a central object for our results.

\begin{definition}[The cut semi-distance $d_{\square}$]
Let $d_{\mathcal{LP}}$ be the Lévy-Pokhorov metric on $\mathcal{P}(\mathbf{Z}).$ The associated cut semi-distance $d_{\square}$ is the function defined on $\Graphon(\Space)^2$ by:
\begin{equation}\label{cutsemi-disteq}
    d_{\square}(U, W)=\sup _{S, T \subset[0,1]} d_{\mathcal{LP}}(U(S \times T ; \cdot), W(S \times T ; \cdot)),
\end{equation}
where the supremum is taken over all measurable subsets $S$ and $T$ of $[0,1]$.
\end{definition}
This distance has been proposed in \cite{abraham2023probabilitygraphons} where the authors have also  shown that it would be topologically equivalent to consider other distances inducing the weak convergence on probability measures instead of the Lévy-Pokhorov metric to define the cut semi-distance $d_{\square}.$ A similar metric was also considered in \cite{athreya2023pathconvergencemarkovchains}. 
In this work, we focus only on the well-known Lévy-Pokhorov metric for simplicity.

\begin{remark}\label{rmk:metrics_2}
All the choices in Remark \ref{rmk:metrics} lead to well-defined cut distances on probability-graphons, and it was shown in \cite[Proposition~1.1]{abraham2023probabilitygraphons} that the resulting topologies are in fact equivalent. 
In particular, one may work interchangeably with the Lévy--Prokhorov distance or these alternative metrics, depending on which formulation is most convenient for a given application.
\end{remark}

\begin{example}
\label{rmk:finitecoloredMetrics}
For the case $\Space=\{0,1,\ldots,r\},$ for 
 $$
W(x,y;\mathrm{d}z)  =w_{r}(x,y)\delta_r(\mathrm{d} z)+\ldots +w_1(x,y)  \delta_1(\mathrm{d} z) +  w_0(x,y) \delta_0(\mathrm{d} z)$$
and 
$$
U(x,y;\mathrm{d}z)  =u_{r}(x,y)\delta_r(\mathrm{d} z)+\ldots +u_1(x,y)  \delta_1(\mathrm{d} z) +  u_0(x,y) \delta_0(\mathrm{d} z)$$
as in Example \ref{rem:finitecoloredgraphs}, $d_{\square}$ is topologically equivalent to the distance $d_{\square_k},$ considered also in \cite{falgasravry2016multicolour}, defined as
\[
d_{\square_k}(W,U)=\sup_{S,\,T \subseteq [0,1]}
\sum_{i=0}^{r} 
\left|\int_{S \times T} \bigl( w_i(x,y) - u_i(x,y) \bigr) \, dx \, dy \right|.
\]
This is the cut norm for real-graphons for $r=1.$
\end{example}

\subsection{Weak isomorphism and the unlabeled  cut distance}
	\label{subsection_unlabeled_cut_distance}
As for the real-graphon case, the next step is to identify probability graphons which differ by some relabeling. This identification will lead again to a quotient space that will be our main playground for proving our results.
Weak isomorphism is the continuum version of the relationship
between adjacency matrices of isomorphic graphs.
Recall $\Relabel$ denotes the set of measure-preserving maps from $[0, 1]$ to $[0, 1]$ equipped with the Lebesgue measure, and $\InvRelabel$ denotes the set of bijective measure-preserving maps from $[0, 1]$ to $[0, 1]$.

We denote by $W^\varphi$ the relabeling of a probability graphon $W$ by a measure-preserving map $\varphi\in \Relabel$, i.e.\
the probability graphon defined for every $x,y\in [0,1]$ and every measurable set $A\subset \Space$ by:
\[
  W^\varphi(x,y;A) = W(\varphi(x),\varphi(y);A)
  \quad \text{for $x,y\in [0,1]$ and $ A\subset \Space$ measurable}.
\]

\begin{definition}[Weak isomorphism]\label{def_weak_isomorphism}
We say that two probability graphons $U$ and $W$ are \emph{weakly isomorphic}
(and we note $U\sim W$) if there exists two measure-preserving maps $\varphi, \psi\in\Relabel$
such that $U^\varphi(x,y; \cdot) = W^\psi(x,y;\cdot)$ for almost every $x,y\in [0,1].$

We denote by $\UGraphon(\Space)=\Graphon (\Space)\,/ \sim$ the space of unlabeled probability graphons, i.e.\
the space of probability graphons where we identify probability graphons that are weakly isomorphic. 
\end{definition}

A natural question, is to ask weather the relabeling modifies in any way the cut (semi-)distance between two probability graphons. This leads to the following definition:
\begin{definition}[The unlabeled cut distance $\delta_{\square}$]\label{def:ddcut}
 The \emph{unlabeled cut distance} is the premetric $\dd$ on $\Graphon(\Space)$ such that for two probability graphons $U$ and $W$:
\begin{equation}\label{def_ddcut}
\dd(U,W) = \inf_{\varphi\in\InvRelabel} d_{\square}(U,W^\varphi)	
= \inf_{\varphi\in\InvRelabel} d_{\square}\left(U^{\varphi},W\right)		 .\end{equation}
\end{definition}

The unlabeled cut distance  $\dd$ is symmetric and satisfies the triangular inequality.  Therefore, $\dd$ defines a distance (that will still be denoted by $\dd $) on the quotient space $\UGraphond(\Space) = \UGraphon (\Space)/ \simd$ of the space of probability graphons with respect to the equivalence relation $\simd$ defined by $U\simd W$  if and only if  $\dd(U,W)=0$.

Is it therefore possible to connect the notion of weak isomorphism and $\delta_\Box$? The following theorem answers to this quesiton.

\begin{theorem}[Lemma 3.17 in \cite{abraham2023probabilitygraphons}]
   \label{theo:Wm=W}
   Two probability graphons are weakly isomorphic, \ie $U \sim W$, if and only if
   $U \simd W$, \emph{\ie}  $\dd(U,W) = 0$.

   Furthermore, the map $\dd$ is a distance on $\UGraphon(\Space)=\UGraphond(\Space).$
\end{theorem}

The unlabeled cut metric can be defined in several different ways as the following lemma states. The following lemma is a special case of Proposition 3.18 in \cite{abraham2023probabilitygraphons}.

\begin{lemma}[Proposition 3.18 in \cite{abraham2023probabilitygraphons}]\label{thm_min_dist}
 For the unlabeled cut distance $\dd$ on $\Graphon$ we have the following equality:
\begin{equation}
  \label{eq_premetric}
  \begin{aligned}
\dd(U,W) 
& = \underset{\varphi\in \InvRelabel}{\inf} d_{\square}(U,W^\varphi) 
  = \underset{\varphi\in \Relabel}{\inf} d_{\square}(U,W^\varphi)\\
& = \underset{\psi\in \InvRelabel}{\inf} d_{\square}(U^\psi,W) 
= \underset{\psi\in \Relabel}{\inf} d_{\square}(U^\psi,W) 	\\
& = \underset{\varphi, \psi\in \InvRelabel}{\inf} d_{\square}(U^\psi,W^\varphi) 
= \underset{\varphi,\psi\in \Relabel}{\min} d_{\square}(U^\psi,W^\varphi).
\end{aligned}
\end{equation}
\end{lemma}

\begin{remark}[Metric for vertex-weighted graphs]
\label{rmk:VertexWeights2}
A similar metric to $\delta_{\square}$ can be defined analogously for pairs $W^{\text{v}}$ and $W^{\text{e}}$ to study also vertex-weighted graph limits, recall Remark \ref{Rmk:Vertex-Weights1} and the notation considered there. Consider pairs $(W^{\text{v}},W^{\text{e}})$ and $(U^{\text{v}},U^{\text{e}})$ of the type of Remark \ref{Rmk:Vertex-Weights1}. We can define a metric in the following way \[
\delta_{\square_{\text{v},\text{e}}}\left((W^{\text{v}},W^{\text{e}}),(U^{\text{v}},U^{\text{e}})\right)=\inf_{\varphi}\left(\sup_{T\subset [0,1]}d_{\mathcal{LP}}\left(\int_T W^{\text{v}}(x)dx,\int_T U^{\text{v}}(\varphi(x))dx\right)+d_{\square}(W^{\text{e}},(U^{\text{e}})^{\varphi})\right)
\]
Observe that in the previous expression the infimum is taken over the same measure-preserving map $\varphi : [0,1]\to[0,1]$ for both kernels $W^{\text{v}}$ and probability graphons $W^{\text{e}}$ when relabeling. This was also hinted to in \cite[Remark 1.7]{abraham2023probabilitygraphons}. All the properties of probability graphons should also hold for these more general objects with minimal changes in the proofs.
\end{remark}

A fundamental result of real-graphon theory is the compactness of the quotient space. The next theorem is an analogous compactness result for families of probability graphons. This is a slight reformulation of Proposition 5.2 in \cite{abraham2023probabilitygraphons}.

\begin{theorem}[Proposition 5.2 in \cite{abraham2023probabilitygraphons}]\label{ThmRelCompact}
Let $\cK\subset \UGraphon(\Space).$ The set $\cK$ is relatively compact for $\delta_{\square}$ if and only if it is \emph{tight}, i.e.\ the set of measures $\{ M_{\widetilde W} : \widetilde W\in\mathcal{K} \} \subset \mathcal{P}(\Space)$ is tight, where $M_{\widetilde W}$ for $\widetilde W\in \UGraphon(\Space)$ is the measure
\begin{equation}
  \label{eq:def-MW}
M_{\widetilde W}(\drv z) 
=  \widetilde W  ([0,1]^2; \drv z)
= \int_{[0,1]^2} \widetilde W (x,y;\drv z) \ \drv x \drv y.
\end{equation}
\end{theorem}
From the theorem follows the next result, pivotal for out theory:

\begin{corollary}[{\cite[Theorem 5.1]{abraham2023probabilitygraphons}}\,]\label{cor:compact}
If $\Space$ is compact then $\UGraphon(\Space)$ is compact. 
\end{corollary}

\begin{remark}
For our results we will focus on the case when $\Space$ is compact. Indeed, as already mentioned, the result of Corollary \ref{cor:compact} is important for our results. As highlighted in Remark \ref{Rmk:WhyCompact} we expect a generalization to the case to non compact $\Space$ to be possible.
\end{remark}

\subsection{Step functions and \texorpdfstring{$n$}{n}-th level approximations.} \label{subsec:stepfunctions}
 An important special case of probability graphons are step-functions probability graphons which are often used for approximation. We dedicate a subsection to this topic because it is going to be a recurring tool that we will use in our proofs. In particular, while for real-graphons the convergence with respect to the cut norm is often powerful enough to determine a convergence in other functional senses (see \cite{chatterjee2017large}) here the underlying structure of the space requires an extra care.

\begin{definition}[Step-functions probability graphons]
  \label{def:stepfunction}
A probability graphon $W\in \Graphon(\Space)$ is a \emph{step-function} if there exists a finite partition of $[0,1]$ into measurable (possibly empty) sets, say $\Pi=\{P_1,\cdots,P_m\}$, such that $W$ is constant on the sets $P_i \times  P_j$, for $1\leq i,j\leq  k$.
\end{definition}

 We will be particularly interested in the following step-function probability graphons obtained from a probability graphon and a given partition of the unit interval. 
 
\begin{definition}[The stepping operator]
	\label{def_stepping_operator}
Let $W\in \Graphon(\Space)$ be a probability graphon and  $\Pi=\{S_1,\cdots,S_k\}$ be a finite partition of $[0,1]$ and recall that $\lambda$ denotes the Lebesgue measure over $[0,1].$
We define the probability graphon $W_\Pi$ adapted to the
partition  $\Pi$ by averaging $W$ over the partition subsets:
\[
  W_\Pi(x,y;\cdot) = \frac{1}{\lambda(S_i)\lambda(S_j)}\,
  W(S_i \times S_j;\cdot)
  \qquad \text{for $ x\in S_i, y\in S_j$,}
\]
when $\lambda(S_i)\neq 0$ and $\lambda(S_j)\neq 0$, and $W_\Pi(x,y;\cdot) = 0$ the null measure otherwise.
\end{definition}

\begin{remark}
The value of $ W_\Pi(x,y;\cdot)$ for  $ x\in
S_i, y\in S_j$ when $\lambda(S_i)\lambda(S_j)=0$ is not relevant as probability graphons are defined up to an almost everywhere equivalence. 
\end{remark}

The following is an approximation result for probability graphons using stepping operators.

\begin{lemma}[Approximation using the stepping operator, {\cite[Lemma 4.5]{abraham2023probabilitygraphons}}]
\label{lem:approx-W-p}
Let $W\in \Kernel$ be a signed measure-valued kernel
 (which is bounded by definition). 
Let  $(\Pi_n  )_{n\in\N}$  be  a  refining sequence  of   finite
partitions  of  $[0,1]$  that  generates the Borel $\sigma$-field on $[0,1]$.   
Then,  the  sequence $(W_{\Pi_n}  )_{n\in\N}$  
is uniformly  bounded by  $\TM{W}$, and
weakly converges to $W$ almost everywhere (on $[0,1]^2$).

In particular, $(W_{\Pi_n}  )_{n\in\N}$ converges to $W$ in labelled cut metric $d_{\square}.$
\end{lemma}
\proof
The first part of the lemma is Lemma 4.5 in \cite{abraham2023probabilitygraphons}. That $(W_{\Pi_n}  )_{n\in\N}$ converges to $W$ in labelled cut metric $d_{\square}$ follows directly using \cite[Lemma 3.13]{abraham2023probabilitygraphons}. 
\endproof

Despite the theory could be developed for general partitions of $[0,1]$, it will be particularly convenient for us to consider the partition $\widehat{\Pi}_n$ of sets of the form \[
S^n_i=\left( \frac{i-1}{n}, \frac{i}{n} \right) ,\]
for $i\in[n]$. This partition generates a partition of the square $[0,1]^2$ into rectangles of the type

\[S^n_i\times S^n_j=\left( \frac{i-1}{n}, \frac{i}{n} \right)\times \left( \frac{j-1}{n}, \frac{j}{n} \right)
\]
together with the remaining set of measure zero. Indeed these partitions are a natural mapping for the vertex set and the edge set of a graph respectively. We will denote by $\widehat{\Pi}_n$ both the partition on $[0,1]$ and the induced one on $[0,1]^2$.

\begin{definition}[$n$--th level approximant]\label{def:Approximants}
The $n$--th level approximant of $W$ is the probability graphon $\widehat{W}_n= W_{\widehat{\Pi}_n},$ where $\widehat{\Pi}_n$ is the partition of $[0,1]^2$ introduced above and $W_{\widehat{\Pi}_n}$ is the stepping operator with respect to that partition.\end{definition}

For $n \geq 1$, let $\cW^n_{1}(\Space)$ be the set of all probability graphons that are constant on the open squares of the form
\[
\left( \frac{i-1}{n}, \frac{i}{n} \right) \times \left( \frac{j-1}{n}, \frac{j}{n} \right)
\]
for every $1 \leq i,j \leq n$. Moreover, we denote with $S^n_{[0,1]}$ the set of all measure-preserving bijections of $[0,1]$
that maps in an injective way any interval of the form $\left(\frac{i-1}{n}, \frac{i}{n}\right)$ to an interval of
the same type, and do not move points of the form $\frac{i}{n}$. Observe that if $\varphi \in S^n_{[0,1]}$ and
$\widetilde{W} \in \cW^n_{1}(\Space)$, then $\widetilde{W}^{\varphi} \in \cW^n_{1}(\Space)$. Moreover, we note that the set $S^n_{[0,1]}$ is in natural bijection with $S_n,$ the
group of permutations of $[n],$ and, therefore, the cardinality of $S^n_{[0,1]}$ is  $|S^n_{[0,1]}| = n!$. Moreover we will denote by $U_n$ denote the set of all measurable functions from $[0,1]$ mapping into $\mathbb{R}$ which are constant in intervals of the form $\left(\frac{i-1}{n}, \frac{i}{n}\right)$. Note that if $u \in U_n$ and $\varphi \in S^n_{[0,1]}$, then $u^{\varphi} \in U_n$, where $(u^\varphi)(x) = u(\varphi(x))$. 

\subsection{Weak regularity lemma }\label{sec:reglemma}
A fundamental notion at the very foundation of dense graphs theory is Szemerédi regularity lemma (see \cite{szemeredi1975regular} and \cite{frieze1999quick} for the weak version). In the dense case for binary graphs it is one of the main ingredients to prove the existence of an LDP in real-graphon space.

The following theorem is a weak Szemerédi regularity lemma for probability graphons. There exists many different versions of the Szemerédi regularity lemma for edge colored graphs \cite[Theorem 1.18]{komlos1995szemeredi}, see also \cite{axenovich2011Szemeredi_version, robert_spectral_szemeredi}, and probability graphons \cite[Corollary 4.14]{abraham2023probabilitygraphons}, see also \cite[Lemma 7.5]{falgasravry2016multicolour} for finitely decorated probability graphons. We give here a slightly modified version. The main difference with Regularity lemma for probability graphons \cite[Corollary 4.14]{abraham2023probabilitygraphons} is point 3 of the theorem.

\begin{theorem}[Weak Regularity Lemma for Probability Graphons]\label{thm:reglemma}
Let $\Space$ be compact. Given any $\varepsilon \in (0,1)$, there exists a set $\cS(\varepsilon) \subset \Graphon(\Space)$ with the following properties:
\begin{enumerate}
    \item The set $\cS(\varepsilon)$ is finite and
    \[
     |\cS(\varepsilon)| \leq M(\varepsilon),
    \]
where $M(\varepsilon)>0$ depends only on $\varepsilon$ and $|\cS(\varepsilon)|$ denotes the cardinality of the set $  |\cS(\varepsilon)|.$ 
    \item For any $\widetilde{W} \in \Graphon$, there exists $\varphi \in \InvRelabel$ and $\widetilde{U}\in \cS(\varepsilon)$ such that
    \[
    d_{\square}(\widetilde{W}^{\varphi}, \widetilde{U}) < \varepsilon.
    \]
In particular, the metric space $(\UGraphon, \delta_{\square})$ is compact and $\cS(\varepsilon)$ is an $\varepsilon-$net for this space.
    \item If $\widetilde{W} \in \cW^n_{1}$, then the function $\varphi$ in part 2 can be chosen to be in the set $S^n_{[0,1]}$ defined in the previous subsection.
\end{enumerate}
\end{theorem}
\proof
We know that the space of probability graphons $\UGraphon(\Space)$ is compact when $\Space$ is compact (Theorem \ref{ThmRelCompact}) and that the space of $\Space-$edge decorated graphs is dense in the space of probability graphons  \cite[Lemma 6.12]{abraham2023probabilitygraphons}. Therefore, the $\Space-$edge decorated graphs with at most $N_{\varepsilon}$ vertices are an $\varepsilon-$net for the space of probability graphons. 

We obtain then the first two points of the Theorem.

We want now to prove point $3.$ Let's assume first that $\Space$ is finite. As we know that the $\Space-$edge decorated graphs with at most $N_{\varepsilon}$ (with $N_{\varepsilon}$ big enough to be fixed later) vertices are an $\varepsilon-$net for the space of probability graphons. Therefore, we have that for any probability graphon $\widetilde{W} \in \cW^n_{1}$ there exists a $\Space-$edge decorated graph $g_{n,M}$ on $n$ vertices such that $\widetilde{W} =\widetilde{W}_{g_{n,M}}. $ Observe that if $n\leq N_{\varepsilon}$ then $\widetilde{W} =\widetilde{W}_{g_{n,M}}$ is an element of the $\varepsilon-$net. Instead, if $n> N_{\varepsilon}$ for $N_{\varepsilon}>0$ big enough depending only on $\varepsilon$, then by the Szemerédi regularity lemma for edge colored graphs \cite[Theorem 1.18]{komlos1995szemeredi}, see also \cite{axenovich2011Szemeredi_version, robert_spectral_szemeredi}, there exists a $\Space-$edge decorated graph $H$ with $m\leq N_{\varepsilon}$ vertices such that $d_{\square}(\widetilde{W}^{\varphi_n}_{g_{n,M}},\widetilde{W}_{H})\leq \varepsilon$ for some $\varphi_n\in S^n_{[0,1]}.$

Let's consider now the general case where $\Space$ is a compact Polish space. Then for any $\gamma>0$ we can choose a finite set $\{x_1,\ldots, x_{N_{\gamma}}\}$ of cardinality $N_{\gamma}$ of elements of $\Space$ such that the set of the probability measures of the form \begin{equation}\label{eq:EmpMeasureProof}
   \sum^{N_{\gamma}}_{i=1}a_i\delta_{x_i}
\end{equation}
where $a_i\geq 0$ and $\sum^{N_{\gamma}}_{i=1}a_i=1$  is a $\gamma-$net in the space of probability measures in $\cP(\Space)$ equipped with the weak topology (or equivalently the Levy-Prokhorov metric). Therefore, for every $\mu\in \cP(\Space)$ there exists a measure $\nu$ as in \eqref{eq:EmpMeasureProof} such that $d_{\mathcal{LP}}(\mu,\nu)<\gamma.$
Now let's consider a probability graphon $\widetilde{W}$ on $\cW^n_1.$ On every open square of the form
\[
\Square=\left( \frac{i-1}{n}, \frac{i}{n} \right) \times \left( \frac{j-1}{n}, \frac{j}{n} \right)
\]
for every $1 \leq i,j \leq n$ the probability graphon is constant, i.e.\ takes as value a probability measure $\mu_{i,j}$ for every $(x,y)\in\Square.$ Then we can approximate this measure with a measure $\nu_{i,j}$ of the form \eqref{eq:EmpMeasureProof} as $d_{\mathcal{LP}}(\mu_{i,j},\nu_{i,j})<\gamma.$ Therefore, we obtain (using quasi-convexity of the Levy-Prokhorov metric, see for example \cite[Lemma 3.21]{abraham2023probabilitygraphons})
\begin{equation*}
    d_{\square}(\widetilde{W},\widetilde{U}_{\gamma})<\gamma
\end{equation*}
for $U_{\gamma}$ taking values only in measures of the form  \eqref{eq:EmpMeasureProof}. The probability graphon $\widetilde{U}_{\gamma}$ is, therefore, taking values in $\cP(\{x_1,\ldots,x_{M_{\gamma}}\}),$ i.e.\ on the space of probability measures on the finite space $\{x_1,\ldots,x_{M_{\gamma}}\}.$ Therefore, we can apply now the result for finite spaces to obtain a $\{x_1,\ldots,x_{M_{\gamma}}\}-$edge colored graph with at most $N_{\varepsilon}$ vertices (or more precisely its probability graphon representation) $\widetilde{U}_{\gamma,\varepsilon}$ such that  $d_{\square}(\widetilde{U}^{\varphi}_{\gamma},\widetilde{U}_{\gamma,\varepsilon})<\varepsilon$ such that $\varphi\in S^n_{[0,1]}$ from result for finite graphs. From the triangular inequality we obtain that 
\begin{equation*}
\begin{aligned}
& d_{\square}(\widetilde{W}^{\varphi},\widetilde{U}_{\gamma,\varepsilon})  
\\
&\leq
 d_{\square}(\widetilde{W}^{\varphi},\widetilde{U}^{\varphi}_{\gamma})+d_{\square}(\widetilde{U}^{\varphi}_{\gamma},\widetilde{U}_{\gamma,\varepsilon})\\
&
=d_{\square}(\widetilde{W},\widetilde{U}_{\gamma})+d_{\square}(\widetilde{U}^{\varphi}_{\gamma},\widetilde{U}_{\gamma,\varepsilon})
 \\
 &<\gamma+\varepsilon.
 \end{aligned}
\end{equation*}
Therefore, we can choose $\varphi\in S^n_{[0,1]}$ for any $\widetilde{W}\in \cW^n_1$ and this proves also the third point in the general case $\Space$ compact. 
\endproof

\section{Relative entropy}\label{Sec:relent}
One of the key steps for proving an LDP is to find an explicit (or at least calculable) formula for the rate function. While in the case of real-graphons arising from \ER graphs the rate function is an easy guess, for general distributions on the edge the right candidate is found to be the relative entropy or Kullback-Leibler distance. This property of the theory makes our result a direct generalization of the classical Sanov's theorem for empirical distributions (see for example \cite[Theorem 6.2.10]{dembo2009large}) on one side and the large deviation principle for \ER graphs by Chatterjee and Varadhan \cite{ChatterjeeVaradhan2011} on the other side. In the following we will try to present our result for the most general $\Space$, but some of them will require finite $\Space$. For some of these results, a generalization to continuous $\Space$ will follow from Section \ref{sec:dawsongartner}, see Remark \ref{rem:continousSpace}.
  
\subsection{Relative entropy for measures}

We begin with the definition of relative entropy for classical measures and some properties of it.

\begin{definition}\label{Def:relEnt}
Let $\nu$  be a probability measure over $\space$, for any signed measure $\omega$ over $\Space$, the relative entropy of $\omega$ with respect to $\nu,$ denoted as $\mathcal{H}(\omega \mid \nu),$ is   \[
\mathcal{H}(\omega \mid \nu) = 
\begin{cases} 
\int \log \left(\frac{d\omega}{d\nu}\right) \, d\omega & \text{if } \omega \ll \nu \text{ and } \omega \text{ is a probability measure} \\
+ \infty & \text{else},
\end{cases}
\]
where $\frac{d\omega}{d\nu}$ is the Radon--Nikodym derivative of $\omega$ with respect to $\nu$.
\end{definition}

In the case of finite dimensional probability spaces the relative entropy takes a particularly simple form.

\begin{example}[Finite case]
Let $\Space=\{0,v_1,\dots, v_n\}$ and let's $\omega$ and $\nu$ the probability measures over $\Space$ uniquely defined by $\omega(v_i)=p_i$ and $\nu(v_i)=q_i$ for every $i\in [n].$ The relative entropy between $\omega$ and $\nu$ then is
\[
\mathcal{H}(\omega \mid \nu) = \sum^n_{i=1} p_i \log \left(\frac{p_i}{q_i}\right)+p_0\log\left(\frac{p_0}{q_0}\right).
\]
\end{example}

Relative entropy enjoys some nice basic properties.
\begin{lemma}[{\cite[Lemma~6.2.12]{dembo2009large}}\,]\label{Lemm:PropRelEntMeas}
The relative entropy $\mathcal{H}(\cdot \mid \nu),$ considered as a function from the space of finite signed measures, equipped with the weak topology, to $\mathbb{R}\cup \{\infty\}$ is a measurable function. Moreover, the relative entropy  function $\mathcal{H}(\cdot \mid \nu)$ is convex and $\mathcal{H}(\omega \mid \nu)\geq 0$ for any finite signed measure $\omega$ and $\mathcal{H}(\omega \mid \nu)= 0$ if and only if $\omega=\nu.$
\end{lemma}
The relative entropy has also the following representation, see for example Dembo-Zeitouni \cite[Lemma~6.2.13]{dembo2009large}.
\begin{lemma}\label{Lemm:RelativeEntropLegendre}
Let $\omega$ be a finite signed measure and $\nu$ be a probability measure on $\Space$. The following equality holds for the relative entropy
\[
\mathcal{H}(\omega \mid \nu) =  \sup_{a\in C_b(\Space)} \left[ \int_\Space a(z) \omega(dz) - \log \int_\Space e^{a(z)} \, \nu(dz) \right].
\]
\end{lemma}
\begin{remark}
In other words, Lemma \ref{Lemm:RelativeEntropLegendre} states that the relative entropy of $\omega$ with respect to $\nu$ is the Legendre transform of the function $\log \int_\Space e^{a(z)} \, \nu(d z).$
\end{remark}

\subsection{Relative entropy for probability graphons.}\label{subsection:RelatEntroProbGra}
We can naturally extend the notion of relative entropy with respect to a probability measure from probability measures to probability graphons in the following way. 

\begin{definition}\label{DefRelatProbGraphon}
Let's consider a probability measure $\nu\in \mathcal{P}(\Space)$ and a probability graphon $W$ from $[0,1]$ into $\mathcal{P}(\Space),$ the relative entropy of $W,$
\[
\widetilde{\mathcal{H}}(W \mid \nu) = \int_{[0,1]^2} \mathcal{H}(W(x,y) \mid \nu) \, dx \, dy.
\]
\end{definition}
\begin{remark}
Observe that the function $\mathcal{H}(W \mid \nu)$ is a measurable function from $[0,1]^2$ into $\mathbb{R}$ as composition of $W$ and $\mathcal{H}$ which are measurable functions. Moreover, as $\mathcal{H}(W(x,y) \mid \nu)\geq 0$ for all $(x,y)\in [0,1]^2,$ the function $\mathcal{H}(W \mid \nu)$ is also integrable. Therefore, the relative entropy for probability graphons $\widetilde{\mathcal{H}}(W \mid \nu)$ is well defined.
\end{remark}
\begin{remark}
The relative entropy from Definition \ref{DefRelatProbGraphon} can be interpreted as the average relative entropy on the space of edges (represented by $[0,1]^2$). In particular, for a constant probability graphon Definition \ref{DefRelatProbGraphon} coincides with the classical relative entropy for probability measures (Definition \ref{Def:relEnt}).
Moreover, as said in the introduction to this section, Definition \ref{DefRelatProbGraphon} is a generalization of what the authors did in \cite{ChatterjeeVaradhan2011}. Their rate function is equivalent to our $\widetilde{\mathcal{H}}$ when the space is $\Space=\{0,1\}$.  In Definition 6.1 in \cite{falgasravry2016multicolour} the authors explore a connected notion, the one of entropy for finite $\Space$.  Similarly in \cite[Definition 5 and Definition 9]{skeja2024quantifyingmultivariategraphdependencies} the authors explore a notion of entropy for multiplex networks. With our work we complete the picture introducing the relative entropy for probability graphons and, furthermore, we generalize these mentioned notions to the continuous type case (see also Remark \ref{rem:continousSpace}).
\end{remark} 

\begin{remark}
\label{rmk:VertexWeightedGraphs3}
One can naturally define a relative entropy also for a probability measure $\nu^{\prime}\in\mathcal{P}(\Space^{\prime})$ and a kernel $W^{\text{v}}:[0,1]\rightarrow \mathcal{P}(\Space^{\prime})$ in the following way  \[
\widetilde{\mathcal{H}}_{\text{v}}(W^{\text{v}}\mid \nu^{\prime})=\int_{[0,1]}
\mathcal{H}(W^{\text{v}}(x)\mid \nu^{\prime})dx.
\]
A special case of this object has been considered in a different context in \cite[Definition 1]{kuehn2024global}. For two probability measures $\nu^{\prime}\in\mathcal{P}(\Space^{\prime})$ and $\nu\in\mathcal{P}(\Space)$ and pairs $(W^{\text{v}},W^{\text{e}})$ as in Remarks \ref{Rmk:Vertex-Weights1} and \ref{rmk:VertexWeights2} the relative entropy can be naturally defined as \begin{equation}\label{eq:RelEntropyVertexWeight}
   \widetilde{\mathcal{H}}_{\text{v},\text{e}}((W^{\text{v}},W^{\text{v}})\mid (\nu^{\prime},\nu))=\widetilde{\mathcal{H}}_v(W^{\text{v}}\mid \nu^{\prime})+\widetilde{\mathcal{H}}(W^{\text{e}}\mid \nu). 
\end{equation}
\end{remark}

The following proposition is a consequence of Theorem \ref{theo:Wm=W} or Lemma \ref{thm_min_dist}. Recall the unlabeled cut metric $\delta_{\square}$ from Definition \ref{def:ddcut} and recall $\Relabel$ denotes the space of measure preserving transformations from $[0,1]$ into $[0,1].$

\begin{proposition}\label{PropWellPosRelEntropyQuot}
Let $\Space$ be a Polish space. If $\delta_{\square}(W, V) = 0$ for two probability graphons $W, V,$ then $\widetilde{\mathcal{H}}(W \mid \nu)=\widetilde{\mathcal{H}}(V \mid \nu)$. 
\end{proposition}
\proof
Let's assume $\delta_{\square}(W, V) = 0.$ Then there exists a $\varphi \in \Relabel$ such that  \[
d_{\square}(W, V^{\varphi}) = 0
\]
from Theorem \ref{theo:Wm=W} or Lemma \ref{thm_min_dist}.

We observe that 
\begin{equation}
    \widetilde{\mathcal{H}}(V \mid \nu) =\int_{[0,1]^2} \mathcal{H}(V(x,y) \mid \nu)\ dx\ dy=\int_{[0,1]^2} \mathcal{H}(V(\varphi(x),\varphi(y))) \mid \nu)\ dx\ dy=\widetilde{\mathcal{H}}(V^{\varphi} \mid \nu) 
\end{equation}
for any $\varphi$ in $\Relabel$ for the properties of integrals and measure-preserving bijections, see equation \eqref{eq:InvBijInt}.
\endproof

As a consequence of Proposition \ref{PropWellPosRelEntropyQuot} we get that the relative entropy is well defined on the space of probability graphons equipped with the unlabeled cut metric $\delta_{\square}$ as well. 

\begin{definition}\label{DefRelativeEntropyProbGraphonsRelabe}
Let $\Space$ a Polish space and $\nu \in \cP(\Space)$ be a probability measure.
The relative entropy on $\UGraphon(\Space)$ is the function
\begin{equation*}
    \widetilde{\mathcal{H}}(\cdot \mid \nu):\UGraphon(\Space) \to \mathbb{R} 
\end{equation*}
where for any $\widetilde{W}\in \UGraphon(\Space),$ then $\widetilde{\mathcal{H}}(\widetilde{W} \mid \nu)=\widetilde{\mathcal{H}}(W\mid \nu),$ where $W$ in the right hand side of this equality is any representative of the class $\widetilde{W}\in \UGraphon(\Space).$ 
\end{definition}

The following corollary is a direct consequence of Lemma \ref{Lemm:RelativeEntropLegendre}.

\begin{corollary}\label{CorProbGraphonRelENTROPLegendre}
Let's consider a probability measure $\nu\in \mathcal{P}(\Space)$ and a probability graphon $W$ from $[0,1]$ into $\mathcal{P}(\Space),$ the relative entropy of $W$, satisfies \[
\widetilde{\mathcal{H}}(W \mid \nu) = \int_{[0,1]^2} \sup_{a \in C_b(\Space)} \left[ \int_\Space a(z) W(x,y,dz) - \log \int_\Space e^{a(z)} \, d\nu(z) \right] \, dx \, dy.
\]
\end{corollary}
Let's define $A$ to be a function from $[0,1]^2$ to $C_b(\Space)$. We have the quantity $J_{A,\nu}: \mathcal{X}\to \mathbb{R}$ defined as 
\[
J_{A,\nu}(W)=\int_{[0,1]^2}  \left[ \int_\Space A(x,y;z) W(x,y,dz) - \log \int_\Space e^{A(x,y;z)} \, d\nu(z) \right]\, dx \, dy.
\]
In general we need $A$ to be regular enough to work with, so we require 
\[
A\in\mathcal{X}^*=(L^2([0,1]^2;\mathcal{M}_{\pm}(\Space)))^*\cong L^2([0,1]^2;\mathbb{R}^s) \cong L^2([0,1]^2;C_b(\Space)).
\] 
We remark that the identification of $\mathcal{X}^*$ with these spaces is only true in the case of $\Space$ finite.
\begin{remark}\label{RemkRelEntropIneqJ}
From Corollary \ref{CorProbGraphonRelENTROPLegendre} it directly follows 
\begin{equation}\label{eq:IneqUgoRelEnt}
    \widetilde{\mathcal{H}}(W \mid \nu) \geq \sup_{A \in  \Ugo}J_{A,\nu}(W).
\end{equation}
\end{remark}
We prove that in the finite dimensional case, the inequality \eqref{eq:IneqUgoRelEnt} is actually an equality. The following is a generalization of \cite[Lemma 5.2]{chatterjee2017large}.

\begin{lemma}\label{LemmEquivLemm5.2Chatterjee}
  Let $\Space=\{v_1,\ldots,v_n\}$ a finite space. Let $\nu\in \cP(\Space)$ and $W$ be a probability graphon from $[0,1]$ into $\cP(\Space).$ The following equality holds \[
\widetilde{\mathcal{H}}(W \mid \nu) = \sup_{A \in  \Ugo}J_{A,\nu}(W).
\]
\end{lemma}
\proof
We already observed the inequality \eqref{eq:IneqUgoRelEnt}.
We now want to prove the other inequality \[
\widetilde{\mathcal{H}}(W \mid \nu) \leq \sup_{A \in { \Ugo}}J_{A,\nu}(W).
\]
Let's consider now  for every $x,y\in [0,1]^2$ the function $A^*$ from $[0,1]^2$ to  $\CbFunct$ defined as \begin{equation}\label{eq:AstarLemmaRelEntrGraphon}
    A^*(x,y;z)=\log \left(\frac{dW(x,y)}{d\nu}(z)\right)-\log \left(\frac{dW(x,y)}{d\nu}(0)\right).
\end{equation}

The function $A^*$ is well defined and measurable (as it is clear from Example \ref{rem:finitecoloredgraphs}), because we are considering the case of $\Space=\{0,v_1,\ldots , v_n\}\cong  \{0\}\cup [n].$ We can also write down the right-hand side in a more explicit way in this case. For $i\in  [n],$ we have

\[
\begin{aligned}
A^*_i(x,y)=A^*(x,y;v_i)&= \log \left( \frac{w_i(x,y)}{p_i} \right) - \log \left( \frac{w_0(x,y)}{p_0} \right)\\
&=\log \left( \frac{w_i(x,y)}{p_i} \right) - \log \left( \frac{1 - \sum^n_{j=1} w_j(x,y)}{1 - \sum^n_{j=1}p_j} \right),
\end{aligned}
\]
where $p_i=\mathbb{P}(\{v_i\})$ for $i\in [n]$ and $p_0=\mathbb{P}(\{0\})$ and for every $i\in [n]$ the real-valued graphon $w_i$ is defined as $w_i(x,y)=W(x,y;v_i)$ and $w_0(x,y)=W(x,y;0).$ In the following we will denote $0$ as $v_0$ when convenient.

Substituting $A^*(x,y)$ in $J_{A^*,\nu}(W)$ we obtain
\[
\begin{aligned}
J_{A^*,\nu}(W)&=
\int_{[0,1]^2} \left[ \int_{\Space} A^*(x,y;z) W(x,y,dz) - \log \int_\Space e^{A^*(x,y;z)} \, d\nu(z) \right]\, dx \, dy\\
&=\int_{[0,1]^2}  \left[\left[ \int_{\Space} \log \left(\frac{dW(x,y)}{d\nu}(z)\right)-\log \left(\frac{dW(x,y)}{d\nu}(0)\right) \right]W(x,y,dz)\right.\\  &\qquad \qquad \qquad\left. - \log \int_\Space  \frac{\left(\frac{dW(x,y)}{d\nu}(z)\right)}{\left(\frac{dW(x,y)}{d\nu}(0)\right)} \, d\nu(z) \right]\, dx \, dy\\
&=\int_{[0,1]^2}  \left[ \int_{\Space} \log \left(\frac{dW(x,y)}{d\nu}(z)\right)W(x,y,dz)\right]\, dx \, dy\\&=
\int_{[0,1]^2}\mathcal{H}(W \mid \nu) \ dx \, dy
\\
&=\widetilde{\mathcal{H}}(W \mid \nu).
\end{aligned}
\]
This would prove the statement of the Lemma if the function $A^*$ would be in $\Ugo$, however, we do not have any guarantee about this. Therefore, we need now to construct a sequence of functions $A^*_n$ such that $J_{A^*_n,\nu}(W)$ converges to $J_{A^*,\nu}(W)=\widetilde{\mathcal{H}}(W \mid \nu).$ 

For all $i\in \{0\}\cup [n]$ we define the measurable subsets of $[0,1]^2$
\begin{align*}
A_{\varepsilon}^i &= \{(x,y) \in [0,1]^2 \mid w_i(x,y) \geq \varepsilon\}, \\
B_{\varepsilon}^i &= \{(x,y) \in [0,1]^2 \mid 0 < w_i(x,y) < \varepsilon \}, \\
E^i &= \{(x,y) \in [0,1]^2 \mid w_i(x,y) = 0 \},\\
P&= \{(x,y) \in [0,1]^2 \mid w_0(x,y) = 1 \}= \{(x,y) \in [0,1]^2 \mid \sum^n_{i=1}w_i(x,y) = 0 \}
\end{align*}
Let's also define the sets 
\begin{align*}
A_{\varepsilon} &= \cap^n_{i=0}A_{\varepsilon}^i=\{(x,y) \in [0,1]^2 \mid w_i(x,y) \geq \varepsilon \text{ for all }i\in \{0\}\cup [n]\}, \\
B_{\varepsilon} &=\cup^n_{i=0}B_{\varepsilon}^i=\{(x,y) \in [0,1]^2 \mid \text{ there exists an } i\in \{0\}\cup [n] \text{ such that }0 < w_i(x,y) < \varepsilon\},  \\
E &= \cup^n_{i=1}E^i=\{(x,y) \in [0,1]^2 \mid \text{ there exists an } i\in [n] \text{ such that }w_i(x,y) =0\}.
\end{align*}

Observe that $P=\cap^n_{i=1}E^i\subset E$ and $P\subset (E^0)^c.$

For each $\varepsilon \in (0,1)$ and $M \geq 1$, define the function from $[0,1]^2$ into $\CbFunct$ defined as
\[
A_{M,\varepsilon}(x,y;\cdot) =
\begin{cases}
A^*(x,y;\cdot) & \text{if } (x,y) \in A_{\varepsilon}, \\
0 & \text{if } (x,y) \in B_{\varepsilon}, \\
A_M^*(x,y;\cdot) & \text{if } (x,y) \in E\cup E^0,
\end{cases}
\]
where the function $A_M^*(x,y;\cdot)$ is defined for  for every $i\in [n]$ as 
\[
A^*_{M}(x,y;v_i) =
\begin{cases}
A^*(x,y;v_i) & \text{if } (x,y) \notin E^0\cup E^i,  \\
\log \left( \frac{w_i(x,y)}{p_i} \right)+M & \text{if } (x,y) \in E^0\setminus E^i,\\
0 & \text{if } (x,y) \in E^0\cap E^i,\\
-M & \text{if } (x,y) \in E^i\setminus E^0
\end{cases}
\]and 

Clearly, $|A_{M,\varepsilon}(x,y;v_i)|\leq 2M$ is bounded for each $(x,y)\in [0,1]^2$ and $i\in \{0\}\cup [n].$  Therefore, the function $A_{M,\varepsilon}$ is in $\Ugo$. Moreover, the function $A_{M,\varepsilon}$ is symmetric since the real-graphons $w_i$ are symmetric for every $i\in \{0\}\cup [n].$ We can now consider the function 
\begin{equation}\label{EqKMProofE0}
\begin{aligned}
    K_{M,\varepsilon}(x,y) &= \mathcal{H}(W(x,y) \mid \nu) - \left[ \int_{\Space} A_{M,\varepsilon}(x,y;z) W(x,y;dz) - \log \int_\Space e^{A_{M,\varepsilon}(x,y;z)} \, d\nu(z) \right]\\
    &=\left(\sum^n_{i=1} w_i(x,y) \log \left(\frac{w_i(x,y)}{p_i}\right)+w_0(x,y)\log\left(\frac{w_0(x,y)}{p_0}\right)\right)\\
    &\qquad \qquad -\left(\sum^n_{i=1}A_{M,\varepsilon}(x,y;v_i)w_i(x,y)-\log(\sum^n_{i=1}e^{A_{M,\varepsilon}(x,y;v_i)}p_i+p_0)\right).
\end{aligned}
\end{equation}

Observe that $K_{M,\varepsilon}(x,y)\geq 0$ from Remark \ref{RemkRelEntropIneqJ}.
The proof is now concluded showing that the integral of $K_{M,\varepsilon}$ tends to zero when $M$ tends to $+\infty$ and $\varepsilon$ tends to $0.$ 
Observe in fact that $K_{M,\varepsilon}(x,y)=0$ when $(x,y)\in A_{\varepsilon}.$ It is easy to verify that $\widetilde{\mathcal{H}}(\cdot \mid \nu)$ is a bounded function on the space of probability graphons over $\Space$ when the space $\Space$ is finite and $\nu$ has all of $\Space$ as support. 
Therefore for each $i\in \{0\}\cup [n]$ we have $0 \leq K_{M,\varepsilon}(x,y;v_i) \leq C(\nu)$ on $B_{\varepsilon}\cup E$, where $C(\nu)$ is a constant depending only on $\nu,$ i.e.\ depending only on the values $p_i>0$ for $i\in \{0\}\cup [n].$ 

Next, for $(x,y) \in  E^0$, we can partition the set $ [n]$ in subsets $N_{(x,y)}$ and $S_{(x,y)}=N^c_{(x,y)}$ such that $(x,y) \in E^i$ if $i\in N_{(x,y)}$ and $(x,y) \in E^0\setminus E^i$ if $i\in S_{(x,y)}.$ Therefore, we can explicitly rewrite equation \eqref{EqKMProofE0} in the following way:

\begin{align*}
    K&_{M,\varepsilon}(x,y)=\\
    =&\left(\sum^n_{i=0} w_i(x,y) \log \left(\frac{w_i(x,y)}{p_i}\right)\right) -\left(\sum^n_{i=1}A_{M,\varepsilon}(x,y;v_i)w_i(x,y)-\log\lep\sum^n_{i=1}e^{A_{M,\varepsilon}(x,y;v_i)}p_i+p_0\rip\right)\\
 =&\left(0+\sum_{i\in S_{(x,y)}}w_i(x,y) \log \left(\frac{w_i(x,y)}{p_i}\right)+\sum_{i\in N_{(x,y)}}0\right)-\left(\sum_{i\in S_{(x,y)}}w_i(x,y) \left(\log \left(\frac{w_i(x,y)}{p_i}\right)+M\right)\right.\\
    &\left. -\log\lep\sum_{i\in S_{(x,y)}}e^{\left(\log \left(\frac{w_i(x,y)}{p_i}\right)+M\right)}p_i+\sum_{i\in N_{(x,y)}}e^0p_i+p_0\rip\right)\\
 =&\left(\sum_{i\in S_{(x,y)}}w_i(x,y) \log \left(\frac{w_i(x,y)}{p_i}\right)\right) -\left(\sum_{i\in S_{(x,y)}}w_i(x,y) \left(\log \left(\frac{w_i(x,y)}{p_i}\right)\right)+M\sum_{i\in S_{(x,y)}}w_i(x,y)\right.\\
    &-\left.\log(e^M\sum_{i\in S_{(x,y)}}\frac{w_i(x,y)}{p_i}p_i+\sum_{i\in N_{(x,y)}}p_i+p_0)\right)\\
    =&\left(-M\sum_{i\in S_{(x,y)}}w_i(x,y) +\log\lep e^{M}\sum_{i\in S_{(x,y)}}w_i(x,y)+\sum_{i\in N_{(x,y)}}p_i+p_0\rip\right)\\
    =&-M +\log\lep e^{M}+\sum_{i\in N_{(x,y)}}p_i+p_0\rip\\
    =&\log\left(1+e^{-M}(\sum_{i\in N_{(x,y)}}p_i+p_0)\right)\longrightarrow 0.
\end{align*}

Similarly, for $(x,y) \in  E\setminus E_0,$ we can partition the set $ [n]$ in subsets $N_{(x,y)}$ and $S_{(x,y)}=N^c_{(x,y)}$ such that $(x,y) \in E^i$ if $i\in N_{(x,y)}$ and $(x,y) \in E\setminus E^i$ if $i\in S_{(x,y)}.$ Therefore, we can explicitly rewrite equation \eqref{EqKMProofE0} in the following way:

\begin{align*}
    K&_{M,\varepsilon}(x,y)=\\
    =&\left(\sum^n_{i=0} w_i(x,y) \log \left(\frac{w_i(x,y)}{p_i}\right)\right) -\left(\sum^n_{i=1}A_{M,\varepsilon}(x,y;v_i)w_i(x,y)-\log\lep\sum^n_{i=1}e^{A_{M,\varepsilon}(x,y;v_i)}p_i+p_0\rip\right)
    \end{align*}
\begin{align*}
 =&\left(w_0(x,y) \log \left(\frac{w_0(x,y)}{p_0}\right)+\sum_{i\in S_{(x,y)}}w_i(x,y) \log \left(\frac{w_i(x,y)}{p_i}\right)+\sum_{i\in N_{(x,y)}}0\right)\\
 & -\left(\sum_{i\in S_{(x,y)}}w_i(x,y)   \left( \log \left(\frac{w_i(x,y)}{p_i}\right)-\log\left(\frac{w_0(x,y)}{p_0}\right)\right)+\sum_{i\in N_{(x,y)}}w_i(x,y) \left(-M\right)\right.\\
    &-\left. \log\lep\sum_{i\in N_{(x,y)}}e^{\left(-M\right)}p_i+\sum_{i\in S_{(x,y)}}e^{ \left( \log \left(\frac{w_i(x,y)}{p_i}\right)-\log\left(\frac{w_0(x,y)}{p_0}\right)\right)}+p_0\rip\right)\\
 =&\left(w_0(x,y)+\sum_{i\in S_{(x,y)}}w_i(x,y)\right) \log \left(\frac{w_0(x,y)}{p_0}\right)\\
 & -\left( -\log\lep e^{-M}\sum_{i\in N_{(x,y)}}p_i+\sum_{i\in S_{(x,y)}}\frac{p_0}{w_0(x,y)}w_i(x,y)+p_0\rip\right)\\ 
    =&\log \left(\frac{w_0(x,y)}{p_0}\right) + \log\left(e^{-M}\sum_{i\in N_{(x,y)}}p_i+\frac{p_0}{w_0(x,y)}\right)\\ 
    =&\log\left(e^{-M}\sum_{i\in N_{(x,y)}}p_i+\frac{p_0}{w_0(x,y)}\right)-\log \left(\frac{p_0}{w_0(x,y)}\right)\longrightarrow 0.
\end{align*}

With these estimates we get the bound  $0\leq  K_{M,\varepsilon}(x,y)\leq C(\nu, M)$ for all $(x,y)\in E\cup E^0,$ where $C(\nu, M)$ is a constant that tends to zero as $M \to \infty$. Combining all the previous observations, we obtain
\begin{equation*}
\begin{aligned}
     &\int_{[0,1]^2} \left(\mathcal{H}(W(x,y) \mid \nu) - \left[ \int_{\Space} A_{M,\varepsilon}(x,y;z) W(x,y;dz) - \log \int_\Space e^{A_{M,\varepsilon}(x,y;z)} \, d\nu(z) \right] \right)\, dx \, dy\\
     &\qquad \qquad=\widetilde{\mathcal{H}}(W \mid \nu)-J_{A_{M,\varepsilon},\nu}(W)=\\
     &\qquad \qquad\int_{[0,1]^2} K_{M}(x, y) \, dx \, dy \\
     &\qquad \qquad \leq C(\nu) \, \text{Leb}(B_{\varepsilon}) + C(\nu, M),
\end{aligned}
\end{equation*}

where $\text{Leb}(B_{\varepsilon})$ is the Lebesgue measure of $B_{\varepsilon}$. Therefore, we get
\[
\widetilde{\mathcal{H}}(W \mid \nu)\leq \sup_{A \in \Ugo} J_{\nu,A}( W) + C(\nu) \, \text{Leb}(B_{\varepsilon}) + C(\nu, M).
\]

Observe that the sets $B_{\varepsilon}$ are bounded sets that monotonically decrease as $\varepsilon \to 0$, and the intersection of all these sets is empty. Therefore, $\text{Leb}(B_{\varepsilon})$ converges to zero as $\varepsilon\to 0$. Recalling that $C(\nu, M)$ converges to zero as $M \to \infty,$ this completes the proof. 
\endproof

From the representation obtained in Lemma \ref{LemmEquivLemm5.2Chatterjee} we obtain the following corollary.

\begin{corollary}[Induced by the cut metric $d$ on $W$]\label{Cor:LowSemicontRelEnt}
Let $\Space$ be a finite space. The function $\widetilde{\mathcal{H}}(\cdot \mid \nu)$ is lower semi-continuous with respect to the topology induced by the cut metric $d_{\square}$ on the space $\Graphon(\Space).$
\end{corollary}

\proof
Directly by the definition of the weak topology we have that the function $J_{\nu,A}(\cdot)$ is continuous under the weak topology for each $A\in \Ugo.$ We also recall that the topology induced by the labelled cut metric $d_{\square}$ is stronger than the weak topology (Lemma \ref{lemm:AppendCutStrongerWeak}). Therefore, these functions are continuous also for the labelled cut metric $d_{\square}.$ Therefore, by Lemma \ref{LemmEquivLemm5.2Chatterjee}, the relative entropy $\widetilde{\mathcal{H}}(\cdot \mid \nu)$ is lower semi-continuous under the cut topology on $\Graphon(\Space),$ recalling that the supremum of any family of lower semi-continuous functions is lower semi-continuous. 
\endproof

From the lower semicontinuity of the relative entropy of a probability graphon and Corollary \ref{Cor:LowSemicontRelEnt} one obtains also the lower semicontinuity of the relative entropy on $\UGraphon(\Space).$ 

\begin{corollary}\label{cor:LowSemicontUGraphon}
Let $\Space$ be a finite space. The relative entropy $\widetilde{\mathcal{H}}(\cdot \mid \nu)$ as a function on $\UGraphon(\Space),$ as in Definition \ref{DefRelativeEntropyProbGraphonsRelabe}, is lower semi-continuous.
\end{corollary}

Corollary \ref{cor:LowSemicontUGraphon} holds more generally for $\Space$ compact but we will see this only in Section \ref{sec:dawsongartner}.

We conclude this section with some additional properties of relative entropy that will be useful in the following sections.

\begin{lemma}\label{lemm:EntropIneqDisc}
Let $\Space$ be a Polish space and $W\in \Graphon(\Space)$ be a probability graphon. Let $(\widehat{W}_n)_n$ be the sequence of $n-$th level approximants of $W.$  The following inequality holds:
\begin{equation*}
    \widetilde{\mathcal{H}}(W \mid \nu)\geq \widetilde{\mathcal{H}}(\widehat{W}_n \mid \nu).
\end{equation*}
\end{lemma}
\proof
The result follows from the fact that the relative entropy for probability measures is convex (Lemma \ref{Lemm:PropRelEntMeas}) and Jensen inequality as $\widehat{W}_n$ is the conditional expectation with respect a finite partition. 
\endproof

Notice that the relative entropy is not a continuous function in general, as one can easily see from the rate function for real-graphons from \cite{ChatterjeeVaradhan2011}. Nevertheless, it is possible to see that it is well behaved with respect approximants as stated in the following corollary.

\begin{corollary}\label{cor:ConvEntropyAproxim}
Let $\Space$ be a finite space and $W\in \Graphon(\Space)$ be a probability graphon. Let $(\widehat{W}_n)_n$ be the sequence of $n-$th level approximants of $W.$  The following equality holds:
\begin{equation*}
 \lim_{n\rightarrow\infty}\widetilde{\mathcal{H}}(\widehat{W}_n \mid \nu)=   \widetilde{\mathcal{H}}(W \mid \nu).
\end{equation*}
\end{corollary}
\proof
Recall that $W_n$ converges to $W$ in $d_{\square}$ and $\widetilde{\mathcal{H}}(\cdot \mid \nu)$ is lower-semicontinuous (Corollary \ref{Cor:LowSemicontRelEnt}). By Lemma \ref{lemm:EntropIneqDisc}, we have  $\widetilde{\mathcal{H}}(\widehat{W}_n \mid \nu)\leq \widetilde{\mathcal{H}}(W_n \mid \nu)$ for every $n$ and, therefore, $\limsup_{n\rightarrow \infty}I_\nu(W_n)\leq I_\nu(W).$ This directly gives the result.
\endproof

In Corollary \ref{cor:ConvEntropyAproxim}, the only reason because we restricted to $\Space$ finite is to use lower-semicontinuity of relative entropy for probability graphons. Howerver, as already mentioned this property of relative entropy holds more generally (Corollary \ref{cor:SemicontinuityFinal}) for $\Space$ compact. Therefore, also Corollary \ref{cor:ConvEntropyAproxim} also holds more generally for $\Space$ compact.

\begin{remark}\label{rem:continousSpace}
As already described in the introduction we performed a \emph{divide et impera} strategy: we first prove our statements for the finite state space $\Space$ and then we lift our results to the case of continuous $\Space$ using Dawson--G\"artner Theorem. For what concerns the relative entropy above most of the results have been proven for a general $\Space$ unless otherwise stated. Among the results where the restriction of $\Space$ finite is required there are Corollary \ref{Cor:LowSemicontRelEnt}, Corollary \ref{cor:LowSemicontUGraphon} and Corollary \ref{cor:ConvEntropyAproxim}. It follows from Section \ref{sec:dawsongartner} that all the above results hold as well in the case of $\Space$ compact, see Corollary \ref{cor:SemicontinuityFinal}. However, in the case of $\Space$ compact the identification of $\mathcal{X}$ and $\mathcal{X}^*$ with nicer spaces is difficult and therefore we take another approach to prove lower--semicontinuity (and its consequences) obtaining these result directly from Dawson--G\"artner Theorem. We also conjecture that a version for continuous $\Space$ of Lemma \ref{LemmEquivLemm5.2Chatterjee} holds, but we don't pursue this direction here.
\end{remark}

\section{LDP for finite dimensional case}
\label{Sec:LDP_finite}
We will now prove a result that we will need in order to prove the main result Theorem \ref{thm_mainresult}. This is an analogous version of Theorem \ref{thm_mainresult} for finite spaces. We base our proof on the outline of the analogous results for real-graphons in \cite{ChatterjeeVaradhan2011,chatterjee2017large}. Along all of this section we will consider $\Space$ to be a finite set of cardinality $|\Space|=s$ for $s\in \N$. Observe that given that $|\Space|=s$ is finite, the spaces $\mathcal{M}_{\pm}(\Space)$ and $\CbFunct$ can be identified with the space of functions from $[s]=\{1,\ldots s\}$ into $\mathbb{R}$ that we denote with $\mathbb{R}^{[s]}.$ Therefore, we have $\mathcal{M}_{\pm}(\Space)\cong \CbFunct \cong \mathbb{R}^{[s]}\cong \mathbb{R}^{s}.$    

Let's consider the space 
\[
\mathcal{X}=L^2([0,1]^2;\mathcal{M}_{\pm}(\Space))\cong L^2([0,1]^2;\mathbb{R}^s)
\] 
with its dual 
\[
\mathcal{X}^*=(L^2([0,1]^2;\mathcal{M}_{\pm}(\Space)))^*\cong L^2([0,1]^2;\mathbb{R}^s)\cong   L^2([0,1]^2;C(\Space)).
\]
We remark that the above identification holds only because we are considering $\Space$ to be finite. In the following we will call $W\in \mathcal{X}$ a finite signed measure kernel and $A\in\mathcal{X}^*$ an element of the dual, that is a $\CbFunct-$valued kernel. By the representation theorem we will denote by $\lambda_A$ the functional acting on $\mathcal{X}$ through the element of the dual $A$. Fix a distribution $\nu$ on $\Space$; without loss of generality we will assume that $\nu$ has as support all the space $\Space$. Following the notation of the introduction, consider a weighted graph $g_{n,\nu}$ distributed as $\mu_{n,\nu}$, i.e. the model where each edge weight is independently sampled from $\nu$. As already mentioned in the introduction, any sampled graph $g_{n,\nu}$ can be naturally associated with an element of $\Graphon(\Space)$. We denote this object by $W_{g_{n,\nu}}$. Concretely, this is done by considering the partition of $[0,1]^2$ into squares
\[
\Square = \left( \tfrac{i-1}{n}, \tfrac{i}{n} \right) \times \left( \tfrac{j-1}{n}, \tfrac{j}{n} \right),
\]
and assigning to each $\Square$ the value given by the Dirac measure at $m_{ij}$, namely $\delta_{m_{ij}} \in \mathcal{P}(\Space)$. In this way, $W_{g_{n,\nu}}$ is obtained as the element of $\mathcal{X}$ constructed from these assignments. Observe that, in fact, $\Graphon(\Space)$ is contained in $\mathcal{X}$. Therefore sampling graphs from $\mu_{n,\nu}$ automatically defines a rule to sample from $\Graphon(\Space) \subset \mathcal{X}$ that we will call by the same symbol. It makes therefore completely sense considering $\mu_{n,\nu}$ as a measure \emph{lifted} on $\Graphon(\Space)$ (or equivalently on $\mathcal{X}$ putting mass zero outside the original support). In the same way the graphs $g_{n,\nu}$ define a measure $\widetilde{\mu}_{n,\nu}$ on the quotient space $\UGraphon(\Space)$. 
Let's consider now the rate function $I_{\nu}(\widetilde{W})=\widetilde{\mathcal{H}}(\widetilde{W} \mid \nu),$ recalling the definition of relative entropy over $\UGraphon(\Space),$ i.e.\ for (classes) of probability graphons (Definition \ref{DefRelativeEntropyProbGraphonsRelabe}). We will prove the following theorem, which is the finite counterpart of the main Theorem \ref{thm_mainresult} and will be used as a base to prove the main result in Section \ref{sec:dawsongartner}:
\begin{theorem}\label{thm:LDP}
Let $\Space$ be a finite space. Consider the space $\UGraphon(\Space)$ equipped with the unlabeled cut metric $\delta_\square$ defined in Section \ref{subsection_unlabeled_cut_distance}. Then for any closed set $\widetilde{\mathcal{C}} \subseteq \widetilde{\mathcal{W}}_1$ and any open set $\widetilde{\mathcal{O}} \subseteq \widetilde{\mathcal{W}}_1$ we have
\[
-\inf_{\widetilde{W} \in \widetilde{\mathcal{O}}} I_\nu(\widetilde{W})\leq\liminf_{n \to \infty} \frac{2}{n^2} \log \widetilde{\mu}_{n,\nu}(\widetilde{\mathcal{O}}\,)\leq\limsup_{n \to \infty} \frac{2}{n^2} \log \widetilde{\mu}_{n,\nu}(\widetilde{\mathcal{C}}\,) \leq -\inf_{\widetilde{W} \in \widetilde{\mathcal{C}}} I(\widetilde{W}) \,.
\]
\end{theorem}
\subsection{Upper bound large deviation in the weak topology}

In the case of the space $\Space$ finite, we will first derive a general upper bound for the large deviation principle using the weak topology. This in turn will be used to prove the upper bound for the LDP in the probability graphon topology, i.e. the topology generated by the \emph{unlabeled cut metric} $\delta_{\square}$. 

The rate function for the large deviations upper-bound in the weak topology will be the relative entropy on $\Graphon(\Space),$ i.e.\ for (labelled) probability graphons, recall Definition \ref{DefRelatProbGraphon}. Therefore, for a probability graphon $W$ we will denote \begin{equation}\label{Eq:RateFuncRelEnt}
    I_{\nu}(W)=\mathcal{H}(W|\nu)
\end{equation}
and we will define $I_{\nu}(W)=\infty$ for any other $W.$ More generally, for an element 
\[
W\in L^2([0,1]^2;\mathcal{M}_{\pm}(\Space))\cong\mathcal{X}= L^2([0,1]^2;\mathbb{R}^s),
\]
we define $I_{\nu}(W)=\mathcal{H}(W|\nu)$ for $W\in \Graphon(\Space)\subset L^2([0,1]^2;\mathcal{M}_{\pm}(\Space))$ and $+\infty$ in every other case.

\begin{theorem}[Large deviations upper-bound in the weak topology]\label{thm:ldpupperboundweak}
Let $\Space$ be a finite space. Given the definitions above, for every weakly closed set $\mathcal{C}\subseteq\mathcal{X}$ we have
\begin{equation}
    \limsup_{n\to\infty} \frac{2}{n^2}\log \mu_{n,\nu}(C)\leq - \inf_{W\in \,\mathcal{C}} I_{\nu}(W),
\end{equation}
where $I_\nu$ is the function defined in equation \eqref{Eq:RateFuncRelEnt} .
\end{theorem}
\begin{proof}
Recall that when $\Space$ is finite then $\mathcal{X}\cong L^2([0,1]^2;\mathbb{R}^s)\cong\mathcal{X}^*$.
Therefore, let $A\in \mathcal{X}^*$ and define   
\[
\lambda_A( W)=\int_{[0,1]^2}\int_\Space A(x,y;z) W(x,y;\mathrm{d}z)\mathrm{d}x\mathrm{d}y
\]
for all $ W \in \mathcal{X}$. Thus, $\lambda_A$ is a continuous linear functional from $L^2([0,1]^2;\mathcal{M}_{\pm}(\Space))$ to $\mathbb{R},$ i.e.\ $\lambda_A\in \mathcal{X}^*.$ 
Define for any continuous linear functional $\lambda\in \mathcal{X}^*$ 
\[
\Lambda_n(\lambda)= \log \int_{\mathcal{X}} e^{\lambda( W)}d\mu_{n,\nu}( W),
\]
where $\mu_{n,\nu}$ is the measure lifted on $\mathcal{X}$ obtained from the measure $\mu_{n,\nu}$ defined on the weighted graph space $\mathcal{G}_n$, as described at the beginning of this section and in the introduction.  Define
\[
\bar{\Lambda}(\lambda)=\limsup_{n\to\infty}{\frac{2 \Lambda_n(n^2\lambda/2)}{n^2}},
\]
we will show that from the scaled quantity $\bar{\Lambda}(\lambda)$ we can obtain an upper bound in the weak topology through some careful steps.
First, given the nature of $\mu_{n,\nu}$ on $\mathcal{X}$, we split the integration domain in blocks $[0,1]^2=\cup_{ij,n} \Square,$ where $\Square=[(i-1)/n,i/n]\times[(j-1)/n,j/n]$ and 

\[
\SquareDiag = \bigcup_{i=1}^{n} \SquareDiagI.
\]

For each $n$ our measure $\mu_{n,\nu}$ on $\mathcal{X}$ will give mass to its support, i.e. those elements in $\mathcal{X}$ such that the value of block $\Square$ is the $\delta_{m_{ij}}$ with probability $1$. It is straightforward from this construction that every block is $\delta_{m_{ij}}$-valued.
For $W\in \mathcal{X}$ constant on the blocks $\Square$ with value $\delta_{m_{i,j}}$ for every $i,j\in [n],$ we have that
\begin{equation}\label{eq:ConstantSquareEq}
\begin{aligned}
\int_{\Square}\lep\int_{\Space}A(x,y;z)W(x,y;dz)\rip dx dy&=\int_{\Square}\lep\int_{\Space} A(x,y;z)\delta_{m_{i,j}}(dz)\rip dx dy\\
&=\int_{\Square}A(x,y;m_{i,j}) dx dy.
\end{aligned}
\end{equation}
Therefore, using equation \eqref{eq:ConstantSquareEq} for $W\in \mathcal{X}$ constant on the blocks $\Square$ with value $\delta_{m_{i,j}}$ for every $i,j\in [n],$ we have 
\begin{align*}
  \lambda_A( W)&=\langle A(x,y;z),W(x,y;dz)\rangle\\
&=\sum_{\substack{i,j=1\\
i\neq j}}\int_{\Square}A(x,y;m_{i,j}) dx dy +\sum_{1\leq i\leq n}\int_{\SquareDiagI} A(x,y;m_{i,i}) dx dy
\\&=\sum_{1\leq i\leq j\leq n}\int_{\SquareUnion}A(x,y;m_{i,j}) dx dy .
\end{align*}
Therefore, we can write
\[
{\frac{2 \Lambda_n(n^2\lambda_A/2)}{n^2}}=\frac{2}{n^2}\log \int_{\mathcal{X}} e^{\frac{n^2}{2}\lambda_A( W)}d\mu_{n,\nu}( W)=\frac{2}{n^2}\log \E e^{\frac{n^2}{2}\sum_{i\leq j}\int_{\SquareUnion}A(x,y;m_{i,j}) dx dy}.
\]
Now, define the $n$-th level approximant $\widehat{A}_n(x,y,\cdot)$ of $A(x,y,\cdot)$ as 
\[
\widehat{A}_n(x,y,\cdot)=\frac{1}{\lambda(S_i)\lambda(S_j)}\int_{\Square}A(x,y,\cdot)dxdy=n^2\int_{\Square}A(x,y,\cdot)dxdy
\]
for $(x,y)\in \Square.$ One should think about these as $n-$th level approximants for $\CbFunct-$kernels, that are different from $n-$th level approximants for probability graphons (Definition \ref{def:Approximants}).

To conclude the proof we will need the following lemma.

\begin{lemma}\label{lemm:L2approx}
Let $\Space$ be a Polish space. In the topology of $L^2([0,1]^2;\CbFunct),$ $\hat{A}_n\to A$ as $n$ goes to infinity.
\end{lemma}
\proof
The proof follows directly from Lemma 1 of \cite[Page 67-68]{diestel1977vector} as $\CbFunct$ is a Banach space. 
\endproof
Now we have 
\[
\begin{split}
    &\frac{2}{n^2}\log \E e^{\frac{n^2}{2}\sum_{i\leq j}\int_{\SquareUnion}A(x,y;m_{i,j}) dx dy}\\=&\frac{2}{n^2}\log \E e^{\frac{n^2}{2}\sum_{i\leq j}\int_{\SquareUnion}\widehat{A}_n(x,y;m_{i,j}) dx dy}\\
    =&\frac{2}{n^2}\log\E \prod_{i\leq j}e^{\frac{n^2}{2}\int_{\SquareUnion}\widehat{A}_n(x,y;m_{i,j}) dx dy}
\end{split}
\]
Using independence we have
\[
\begin{split}
    &\frac{2}{n^2}\log \prod_{i\leq j}\EE  e^{\frac{n^2}{2}\int_{\SquareUnion}\widehat{A}_n(x,y;m_{i,j}) dx dy}=\\
    &\frac{2}{n^2}\sum_{i\leq j}\log \EE  e^{\frac{n^2}{2}\int_{\SquareUnion}\widehat{A}_n(x,y;m_{i,j}) dx dy}=\\
    &\sum_{i\leq j}\int_{\SquareUnion}\log \EE e^{\widehat{A}_n(x,y;m_{i,j})}dx dy=\\
    &\int_{[0,1]^2}\log \EE e^{\widehat{A}_n(x,y;m_{i,j})}=\\
    &\int_{[0,1]^2\backslash \SquareDiag}\log \EE e^{\widehat{A}_n(x,y;m_{i,j})}+ \int_{ \SquareDiag}\log \EE e^{\widehat{A}_n(x,y;m_{i,j})}
\end{split}
\]
Using the Lemma \ref{lem:BoundLogMomen}, we have 
\[|\log\EE e^{A(x,y;m_{i,j})}-\log\EE e^{\hat{A}(x,y;m_{i,j})}|\leq|A(x,y;m_{i,j})-\widehat{A}_n(x,y;m_{i,j})|\]
and, therefore, using Lemma \ref{lemm:L2approx}, we obtain that 

\[
\begin{aligned}
&\int_{[0,1]^2}|\log\EE e^{A(x,y;m_{i,j})}-\log\EE e^{\widehat{A}_n(x,y;m_{i,j})}|\mathrm{d}x\mathrm{d}y\leq \\ & \int_{[0,1]^2}|A(x,y;m_{i,j})-\widehat{A}_n(x,y;m_{i,j})|\mathrm{d}x\mathrm{d}y\leq\\ & \|A(\cdot,\cdot;m_{i,j})-\widehat{A}_n(\cdot,\cdot;m_{i,j})\|_2 \leq \\ & \|A-\widehat{A}_n\|_2\rightarrow 0    
\end{aligned}
\]
where we used the Cauchy-Schwartz inequality in the second inequality and in the third line $\|\cdot\|_2$ represents the $L^2-$norm on $L^2([0,1]^2;\R)$ and on the fourth line $\|\cdot\|_2$ represents the $L^2-$norm on $L^2([0,1]^2;\CbFunct)$ (recall $\Space$ is finite). Moreover, from Lemma \ref{lem:BoundLogMomen} again we have 

\[
|\log\EE e^{\widehat{A}_n(x,y;m_{i,j})}|\leq |\widehat{A}_n(x,y;m_{i,j})|
\]
and, therefore,
\[
\begin{aligned}
\int_{\SquareDiag}|\log\EE e^{\widehat{A}_n(x,y;m_{i,j})}|\mathrm{d}x\mathrm{d}y & \leq \int_{\SquareDiag}|\widehat{A}_n(x,y;m_{i,j})|\mathrm{d}x\mathrm{d}y\\ &=\left(Leb(\SquareDiag)\|\widehat{A}_n(\cdot,\cdot;m_{i,j})\|^2_2 \right)^{\frac{1}{2}}\\&
=\frac{1}{n}\|\widehat{A}_n\|_2 \rightarrow 0
\end{aligned}\]
where we used Cauchy-Schwartz inequality and $Leb(\SquareDiag)=\tfrac{1}{n}$ and the fact that $\widehat{A}_n$ converges to $A$ in $L^2$ as $n$ goes to $\infty,$ see Lemma \ref{lemm:L2approx}.

Therefore, the contribution of the diagonal can be neglected and we thus have proved that 
\[
\bar{\Lambda}(\lambda_A)=\int_{[0,1]^2}\log\EE e^{A(x,y;m_{i,j})}dxdy.
\] Now define, for any $w\in\mathcal{X}$
\[
\bar{\Lambda}^*(W)=\sup_{\lambda\in \mathcal{X}^*} (\lambda(W)-\bar{\Lambda}(\lambda)).
\]
We have then by Lemma \ref{LemmEquivLemm5.2Chatterjee}
\[
\bar{\Lambda}^*(W)\geq \sup_A (\lambda_A(W)-\bar{\Lambda}(\lambda_A))=\sup_A J_A(W)=I_{\nu}(W),
\]
where $I_{\nu}(W)$ is the function $I_{\nu}(W)=\widetilde{\mathcal{H}}(W \mid \nu),$  the relative entropy for probability graphons from Definition \ref{DefRelatProbGraphon}.
We observe now that as $\Space$ is finite, the space $\mathcal{X}^*$ (and $\mathcal{X}$ too) is an Hilbert space and thus reflexive. Thus, the weak topology is equivalent to the weak$-*$ topology. Therefore, by the Banach-Alaoglu theorem, bounded sets of $\mathcal{X}^*$ are precompact with respect to the weak$-*$ topology. 
The proof is now completed using the compactness of the weak$-*$ topology and the classical large deviation result Theorem \ref{ThmGeneralLargeDevChatterjee4.1}.
\end{proof}

\subsection{General Upper bound}
Now we will work with classes of probability graphons $\widetilde{W}\in \UGraphon(\Space).$
Let's consider now the rate function $I_{\nu}(\widetilde{W})=\widetilde{\mathcal{H}}(\widetilde{W} \mid \nu),$ recalling the definition of relative entropy for (classes) of probability graphons (Definition \ref{DefRelativeEntropyProbGraphonsRelabe}).

\begin{theorem}[Upper bound]\label{thm:upperbound}
Let $\Space$ be a finite space. Consider the space $\UGraphon(\Space)$ equipped with the cut metric $\delta_\square$ defined in Section \ref{subsection_unlabeled_cut_distance}. Then for any closed set $F \subseteq \UGraphon(\Space)$,
\[
\limsup_{n \to \infty} \frac{2}{n^2} \log \widetilde{\mu}_{n,\nu}(F) \leq -\inf_{\widetilde{W} \in F} I_{\nu}(\widetilde{W}).
\]
\end{theorem}

For $\widetilde{W} \in \UGraphon(\Space)$ and $\eta > 0$, define
\[
S_\square(\widetilde{W}, \eta) = \{ \widetilde{U} \in \UGraphon(\Space) : \delta_\square(\widetilde{W}, \widetilde{U}) \leq \eta \}.
\]
Thanks to Lemma \ref{lemma:upperbound_balls} (\cite[Lemma 4.1]{chatterjee2017large}), to prove our statement is sufficient to prove that for every $\widetilde{W} \in \UGraphon(\Space)$,
\begin{equation}\label{eq:inequalitytoprove}
\lim_{\eta \to 0} \limsup_{n \to \infty} \frac{2}{n^2} \log \widetilde{\mu}_{n,\nu}(S_\square(\widetilde{W}, \eta)) \leq - I_{\nu}(\widetilde{W}). 
\end{equation}
where $I_{\nu}(\cdot)=\widetilde{\mathcal{H}}(\cdot \mid \nu)$ is the function defined in Definition \ref{DefRelativeEntropyProbGraphonsRelabe}.
Let $B(\widetilde{W}, \eta) \subset \UGraphon(\Space)$ be defined as
\[
B(\widetilde{W}, \eta) = \{ U \in \Graphon(\Space) : \widetilde{U} \in S_\square(\widetilde{W}, \eta) \}.
\]
Using the definitions of $\widetilde{\mu}_{n,\nu}$ and $\mu_{n,\nu}$, we have
\[
\mu_{n,\nu}(B(\widetilde{W}, \eta)) = \widetilde{\mu}_{n,\nu}(S_\square(\widetilde{W}, \eta)).
\]

Therefore, for our theorem it suffices to prove that for every $\widetilde{W} \in \UGraphon(\Space)$,
\begin{equation}
\lim_{\eta \to 0} \limsup_{n \to \infty} \frac{1}{n^2} \log \mu_{n,\nu}(B(\widetilde{W}, \eta)) \leq -I(\widetilde{W}).
\end{equation}

Now recall the sets $\mathcal{S}(\varepsilon)$ from Theorem \ref{thm:reglemma}.

\begin{lemma}[bound on the function $\varphi$ of Theorem \ref{thm:reglemma}]\label{lem:regularityupperbound}
Take any $\varepsilon > 0$ and let $\mathcal{S}(\varepsilon)$ be as in Theorem \ref{thm:reglemma}. Then for any $\widetilde{W} \in \UGraphon(\Space)$ and $\eta > 0$,
\[
\mu_{n,\nu}(B(\widetilde{W}, \eta)) \leq n! \, \mu_{n,\nu}(B(\widetilde{W}, \eta) \cap B(\mathcal{S}(\varepsilon), \varepsilon)),
\]
where $B(\mathcal{S}(\varepsilon), \varepsilon) = \{ U \in \Graphon(\Space) : \min_{Q \in \mathcal{S}(\varepsilon)} d_\square(U, W) \leq \varepsilon \}$.
\end{lemma}
\begin{proof}
Let $g_{n,\nu}$ be a random weighted graph, sampled from the distribution $\mu_{n,\nu}$ on $\mathcal{G}_n$ (recall the definitions and the notation at the beginning of this section and the introduction).  By statement $(3)$ of Theorem \ref{thm:reglemma} , there exists $\varphi \in S^n_{[0,1]}$ and $Q \in \mathcal{S}(\varepsilon)$ such that $d_\square(W_{g_{n,\nu}}^\varphi, Q) < \varepsilon$, where $W_{g_{n,\nu}}^\varphi$ is the probability graphon of $g_{n,\nu}$ under $\varphi.$ Let $(\Omega,\P)$ be the probability space where the random graph $g_{n,\nu}$ is defined. Thus,
\[
\begin{split}
    \mu_{n,\nu}(B(\widetilde{W}, \eta)) &= \mathbb{P}(W_{g_{n,\nu}}:W_{g_{n,\nu}} \in B(\widetilde{W}, \eta))\\
&\leq \sum_{\varphi \in S^n_{[0,1]}} \mathbb{P}(W_{g_{n,\nu}}^\varphi \in B(\widetilde{W}, \eta) \cap B(\mathcal{S}(\varepsilon), \varepsilon))\\
&\leq n!\mathbb{P}(W_{g_{n,\nu}} \in B(\widetilde{W}, \eta) \cap B(\mathcal{S}(\varepsilon), \varepsilon))\\
&=  n!\mu_{n,\nu}(B(\widetilde{W}, \eta) \cap B(\mathcal{S}(\varepsilon), \varepsilon)).
\end{split}
\]
where in the second to last inequality we used that $W_{g_{n,\nu}}^\varphi$ has the same probability law as $W^{g_{n,\nu}}$ for any $\varphi \in S^n_{[0,1]}$ and $| S^n_{[0,1]}| = n!$.
\end{proof}

For any $U \in \Graphon(\Space)$ and $\varepsilon > 0$, define
\[
B(U, \varepsilon) := \{ W \in \Graphon(\Space) : d_\square(U, W) \leq \varepsilon \}.
\]

\begin{lemma}[Lemma 5.4 chatterjee]\label{lem:weaklyclosedness}
Let $\Space$ be finite. For any $U \in \Graphon(\Space)$ and $\varepsilon > 0$, $B(U, \varepsilon)$ is weakly closed.
\end{lemma}
\begin{proof}
Suppose that $(U_n)_{n}$ is a sequence in $\Graphon(\Space)$ such that $U_n \in B(W, \varepsilon^{\prime})$ for each $n$ and $U_n \to U$ weakly.  By definition of weak convergence, from $U_n \to U$ weakly, we have that for each $A\in L^2([0,1];\R^s)\cong L^2([0,1];\CbFunct).$
\begin{equation}\label{eq:weakConvUgo}\begin{aligned}
     &\left| \int_{[0,1]^2} \langle A(x,y;\cdot),(U(x,y;\cdot) - W(x,y;\cdot))\rangle_{\R^s} \, dx dy \right|
\\& \hspace{20pt}= \lim_{n \to \infty} \left| \int_{[0,1]^2} \langle A(x,y;\cdot),(U_n(x,y;\cdot) - W(x,y;\cdot))\rangle_{\R^s} \, dx dy \right|,
\end{aligned}
   \end{equation}
where for $v,w\in \R^s$ the scalar product $\langle v,w \rangle_{\R^s}=\sum^n_{i=1}v_iw_i.$  For  any two Borel measurable functions $a, b : [0,1] \to [-1,1]$ and $i\in [n]$ we define the function $A_{a,b,i}\in L^2([0,1];\R^s)$ as
$$A_{a,b,i}(x,y;j)=\begin{cases}
    a(x)b(y) \text{ if } j=i,\\
    0 \text{ else}.
\end{cases}$$
Substituting $A_{a,b,i}$  for $i\in [n]$ in \eqref{eq:weakConvUgo} we obtain
\[
\begin{split}
    &\left| \int_{[0,1]^2} a(x) b(y)(U(x,y;i) - W(x,y;i)) \, dx dy \right|\\
 &\qquad=\lim_{n \to \infty} \left| \int_{[0,1]^2} a(x) b(y)(U_n(x,y;i) - W(x,y;i)) \, dx dy \right| \leq \varepsilon^{\prime}.
\end{split}
\]
Taking supremum over all Borel measurable functions $a, b : [0,1] \to [-1,1]$ gives $d_\square(U, W) \leq \varepsilon$, for an $\varepsilon$ small as much as we want choosing $\varepsilon^{\prime}$ small enough (see also Example \ref{rmk:finitecoloredMetrics}).
\end{proof}

\begin{lemma}\label{lem:upperbndinequalities}
There exists a function $\delta(\widetilde{W}, \varepsilon)$, depending only on $\widetilde{W}$ and $\varepsilon$, with $\delta(\widetilde{W}, \varepsilon) \to 0$ as $\varepsilon \to 0$, such that for each $\widetilde{W} \in \UGraphon(\Space)$, $\eta > 0$ and $\varepsilon > 0$,
\[
\lim_{\eta \to 0} \limsup_{n \to \infty} \frac{1}{n^2} \log \mu_{n,\nu}(B(\widetilde{W}, \eta) \cap B(\mathcal{S}(\varepsilon), \varepsilon)) \leq - I_{\nu}(\widetilde{W}) + \delta(\widetilde{h}, \varepsilon).
\]
\end{lemma}

\begin{proof}
Since $\mathcal{S}(\varepsilon)$ is a finite set, it suffices to show that for fixed $U \in \mathcal{S}(\varepsilon)$,
\[
\lim_{\eta \to 0} \limsup_{n \to \infty} \frac{1}{n^2} \log \mu_{n,\nu}(B(\widetilde{W}, \eta) \cap B(U, \varepsilon)) \leq - I_{\nu}(\widetilde{W}) + \delta(\widetilde{W}, \varepsilon).
\]
If $B(\widetilde{W}, \eta) \cap B(U, \varepsilon)$ is empty for sufficiently small $\eta$, then there is nothing to prove. So let us assume that this is not the case. Then
\begin{equation}\label{eq:UinBall}
U \in B(\widetilde{W}, \varepsilon). 
\end{equation}

By lower semi-continuity of $I_{\nu}$, $I_{\nu}(V) \geq I_{\nu}(\widetilde{W}) - \delta(\widetilde{W}, \varepsilon)$ for $V\in B(\widetilde{W}, 2\varepsilon)$, where $\delta(\widetilde{W}, \varepsilon) \to 0$ as $\varepsilon \to 0$. But by \eqref{eq:UinBall} we have $B(U, \varepsilon) \subseteq B(\widetilde{W}, 2\varepsilon)$ and by Lemma \ref{lem:weaklyclosedness}, $B(U, \varepsilon)$ is weakly closed. Therefore by Theorem \ref{thm:ldpupperboundweak},
\[
\begin{aligned}
\lim_{\eta \to 0} \limsup_{n \to \infty} \frac{1}{n^2} \log \mu_{n,\nu}(B(\widetilde{W}, \eta) \cap B(U, \varepsilon)) 
&\leq \limsup_{n \to \infty} \frac{1}{n^2} \log \mu_{n,\nu}(B(U, \varepsilon)) \\
&\leq - \inf_{V \in B(U, \varepsilon)} I_{\nu}(V) \\
&\leq - \inf_{V \in B(\widetilde{W}, 2\varepsilon)} I_{\nu}(V) \leq - I_{\nu}(\widetilde{W}) + \delta(\widetilde{W}, \varepsilon).
\end{aligned}
\]
\end{proof}

\begin{proof}[Proof of the Upper Bound in Theorem \ref{thm:LDP}]
The combination of Lemma \ref{lem:regularityupperbound} and Lemma \ref{lem:upperbndinequalities} gives the inequality \eqref{eq:inequalitytoprove} after taking $\varepsilon \to 0$, which completes the proof of the upper bound of Theorem \ref{thm:upperbound}.
\end{proof}

\subsection{General Lower bound}

\begin{theorem}[Lower bound]\label{thm:lowerbound}
Let $\Space$ be finite. Consider the space $\widetilde{\mathcal{W}}_{1}(\Space)$ equipped with the unlabeled cut metric $\delta_\square$ defined in Section \ref{subsection_unlabeled_cut_distance}. Then for any open set $\widetilde{\mathcal{O}} \subseteq \widetilde{\mathcal{W}}_1$,
\[
\liminf_{n \to \infty} \frac{2}{n^2} \log \widetilde{\mu}_{n,\nu}(\widetilde{\mathcal{O}}) \geq -\inf_{\widetilde{W} \in \widetilde{\mathcal{O}}} I_\nu(\widetilde{W}).
\]
\end{theorem}
We will proceed with some reductions and prove the above statement at the end of this sections. The reductions performed will allow us to work with $n-$th level approximants, recall Definition \ref{def:Approximants}. We will then connect the rate function $I_\nu(W)$ to the one of $I_\nu(W_n)$ and conclude our theorem with a tilting argument and a control on how probable is that $W_n$ is far from $W$ in the cut distance.
Given that we want to prove a lower bound we are allowed to just consider for all $\widetilde{W} \in \UGraphon(\Space)$ and $\eta \in (0,1)$ the unlabeled ball $\mathcal{S}_{\Box}$ centered at $\widetilde{W}$ with radius $\eta$:
\[
\liminf_{n \to \infty} \frac{1}{n^2} \log \widetilde{\mu}_{n,\nu}(\mathcal{S}_{\Box}(\widetilde{W},\eta)) \geq -I_\nu(\widetilde{W}).
\]
Indeed we have that for any open $\widetilde{\mathcal{O}}$ and $\widetilde{W} \in \widetilde{\mathcal{O}}$, there exists $\eta \in (0,1)$ such that $\mathcal{S}_{\Box}(\widetilde{W},\eta) \subseteq \widetilde{\mathcal{O}}$. We can proceed further. Indeed we can consider sets in the labeled topology. Since $B(W,\eta) \subseteq B(\widetilde{W},\eta)$ and $\mu_{n,\nu}(B(\widetilde{W},\eta)) = \tilde{\mu}_{n,\nu}(\mathcal{S}_{\Box}(\widetilde{W},\eta))$, it suffices to prove that for any $W \in \mathcal{W}$ and $\eta \in (0,1)$,
\[
\liminf_{n \to \infty} \frac{1}{n^2} \log\mu_{n,\nu}(B(W,\eta)) \geq -I_\nu(W).
\]
Not only we can take advantage of the fact that we are proving a lower bound by simplifing the set $\widetilde{\mathcal{O}}$, given an arbitrary $\widetilde{W}$ (and $W$) we can approximate it in multiple steps as we already did in the previous sections. Take any $W \in \Graphon(\Space)$ and $\eta > 0$. Let $\widehat{W}_n$ be the level $n-$th level approximant of $W.$ Therefore, from Lemma \ref{lem:approx-W-p} (Lemma 4.5 in \cite{abraham2023probabilitygraphons}), we have $\widehat{W}_n \to W$ in the labelled cut metric $d_{\square}.$ Consequently, $B(\widehat{W}_n,\varepsilon) \subseteq B(W,\eta)$ for all large $n$. Thus, it suffices to prove that
\[
\liminf_{n \to \infty} \frac{1}{n^2} \log \mu_{n,\nu}(B(\widehat{W}_n,\varepsilon)) \geq -I_\nu(W).
\]

As in the proof of Theorem \ref{thm:ldpupperboundweak}, we will divide the interval $[0,1]^2$ in a dyadic partition. Let again $\Square$ denote the square $[(i-1)/n,i/n] \times [(j-1)/n,j/n]$ and
\[
\SquareDiag = \bigcup_{i=1}^{n} \SquareDiagI.
\]

Define $W_n$ to be a copy of $\widehat{W}_n$ but with the element $\delta_0$ on $\SquareDiag$ (see \cite[Remark 6.2]{abraham2023probabilitygraphons}). As already noticed, the Lebesgue measure of the set $\SquareDiag$ tends to zero as $n \to \infty$, therefore $d_{\Box}(W_n, \widehat{W}_n) \to 0$. We then need to prove that for any $\varepsilon \in (0,1)$,
\begin{equation}\label{eq:finalapprox}
\liminf_{n \to \infty} \frac{1}{n^2} \log \mu_{n,\nu}(B(W_n,\varepsilon)) \geq -I_\nu(W).     
\end{equation}

Being $W_n(i,j;dz)$ the value of  $W_n$ on $\Square$ we can sample from $W_n$ a realization of a weighted graph  that we will call $h_{n,M}$ where we have that $M=(m_{i,j}\sim W_n(i,j;dz))$ (this is equivalent to $\mathbb{H}(x,W_n)$ in \cite{abraham2023probabilitygraphons}). As we have already done before, we will interpret $h_{n,M}$ as a $\mathcal{P}(\Space)$-valued random graph where we associate $\delta_{m_{i,j}}$ to each edge. $h_{n,M}$ lives on some abstract probability space $(\Omega, \mathscr{F}, \mathbb{P})$ and is constructed by considering a weight between vertices $i$ and $j$ sampled from $W_n(i,j;dz)$, for every $1 \leq i < j \leq n$.
 $h_{n,M}$ can be naturally interpreted as a probability graphon, $W_{h_{n,M}}$, and will have a specific law that we denote by $\mathbb{P}_{n,W}$ given that it depends on the initial probability graphon $W$ and its $n-$th level approximation.

Is easy to realize that the support of $\mathbb{P}_{n,W}$ is contained in the support of $\mu_{n,\nu}$ (assuming the support of $\nu$ is all of $\Space$), i.e. all probability graphons arising from weighted (or better said
 $\Space-$edge decorated) graphs of size $n$; we call this space $\mathcal{W}_1\mid_{n,\nu}$. We will prove the following:

\begin{lemma}\label{lem:radonnykodym}
Let $\mathbb{P}_{n,W}$ and $\mathcal{W}_1\mid_{n,\nu}$ be defined as above. And let $V$ be any element of $\mathcal{W}_1\mid_{n,\nu}$, (notice that $\mathcal{W}_1\mid_{n,\nu}$ is the support of $\mathbb{P}_{n,W}$). Then
\[
\lim_{n \to \infty} \frac{2}{n^2} \int_{\mathcal{W}_1} \log \frac{d\mathbb{P}_{n,W}}{d\mu_{n,\nu}}(V) \, d\mathbb{P}_{n,W}(V) = I_\nu(W).
\]
\end{lemma}

\begin{proof}
Recall Definition \ref{DefRelatProbGraphon}. From the result of Section \ref{subsection:RelatEntroProbGra}, we know that we can write for the class $\widetilde{W}$ 
\[
\begin{split}
I_\nu(\widetilde{W})&=\widetilde{\mathcal{H}}(\widetilde{W}\mid \nu)=\int_{[0,1]^2} \mathcal{H}(\widetilde{W}(x,y)\mid \nu) dxdy\\
&=\int_{[0,1]^2} \lep \int_\Space\log\lep\frac{W(x,y;dz)}{ \nu(dz)}\rip W(x,y;dz) \rip dx dy,
\end{split}
\]
where $W$ is any representative of the class $\widetilde{W}$. We are going to show that 
\[
\frac{2}{n^2} \int_{\mathcal{W}_1} \log \frac{d\mathbb{P}_{n,W}}{d\mu_{n,\nu}}(V) \, d\mathbb{P}_{n,W}(V)=I_\nu(W_n)
\]
and that
\[\lim_{n\to\infty} I_\nu(W_n)=I_\nu(W)=\int_{[0,1]^2} \lep \int_\Space\log\lep\frac{W(x,y;dz)}{ \nu(dz)}\rip W(x,y;dz) \rip dx dy.
\]
First notice that we can rewrite 
\[
\frac{2}{n^2} \int_{\mathcal{W}_1} \log \frac{d\mathbb{P}_{n,W}}{d\mu_{n,\nu}}(V) \, d\mathbb{P}_{n,W}(V)=\frac{2}{n^2} \int_{\mathcal{W}_1\mid_{n,\nu}}\log \frac{d\mathbb{P}_{n,W}}{d\mu_{n,\nu}}(V) \, d\mathbb{P}_{n,W}(V).
\]
Now, for $V\in \mathcal{W}_1\mid_{n,\nu}$ we have $V(i,j;dz)=\delta_{f_{ij}}(dz)$ (for some $f_{ij}\in \Space)$ which denotes the value of $V$ in $\Square.$ Moreover, we observe that for $\Space=\{v_0=0\}\cup\{v_1,\ldots,v_{s}\}$ 

\[
\mathbb{P}_{n,W}(V)=\prod_{1 \leq i < j \leq n} \prod^s_{k=0}W_n(i,j;\{v_k\})^{\delta_{f_{ij}}(\{v_k\})}
\]
and
\[
\mu_{n,\nu}(V)=\prod_{1 \leq i < j \leq n} \prod^s_{k=0}\nu(\{v_k\})^{\delta_{f_{ij}}(\{v_k\})}.
\]

Furthermore, notice that 
\[
\begin{aligned}
&\frac{2}{n^2} \int_{\mathcal{W}_1\mid_{n,\nu}}\log \frac{d\mathbb{P}_{n,W}}{d\mu_{n,\nu}}(V) \, d\mathbb{P}_{n,W}(V)\\
&=\frac{2}{n^2} \int_{\mathcal{W}_1\mid_{n,\nu}} \log\left( \prod_{1\leq i<j\leq n}\prod^s_{k=0}\left(\frac{W_n(i,j;\{v_k\})}{\nu(\{v_k\})} \right)^{\delta_{f_{ij}}(\{v_k\})}\right)\, d\mathbb{P}_{n,W}(V)
\\
&=\frac{2}{n^2} \int_{\mathcal{W}_1\mid_{n,\nu}} \sum_{1\leq i<j\leq n}\sum^s_{k=0}{\delta_{f_{ij}}(\{v_k\})}\log\left( \frac{W_n(i,j;\{v_k\})}{\nu(\{v_k\})} \right)\, d\mathbb{P}_{n,W}(V)
\\
&=
\frac{2}{n^2}\int_{\mathcal{W}_1\mid_{n,\nu}}\sum_{1\leq i<j\leq n} \int_{\Space}V(i,j;dz)  \log \left(\frac{W_n(i,j;dz)}{\nu(dz)}\right) \, d\mathbb{P}_{n,W}(V)
\end{aligned}
\]

Interchanging the integrals, that we can exchange because they are integrals on finite spaces, and using that $\text{Leb}(\Square)= \lambda(S_i)\lambda(S_j)=\frac1{n^2}$, we have by definition of the probability $\mathbb{P}_{n,W}$
\[
\begin{split}
&\frac{2}{n^2} \sum_{1\leq i<j\leq n}\int_{\mathcal{W}_1\mid_{n,\nu}}\int_\Space V(i,j;dz)\log \frac{W_n(i,j;dz)}{\nu(dz)} \, d\mathbb{P}_{n,W}(V)\\
&=\frac{2}{n^2} \sum_{1\leq i<j\leq n}\int_{\Space}\log \frac{W_n(i,j;dz)}{\nu(dz)}\int_{\mathcal{W}_1\mid_{n,\nu}}V(i,j;dz) \, d\mathbb{P}_{n,W}(V)
\\&=\frac{2}{n^2} \sum_{1\leq i<j\leq n}\int_{\Space}\log \frac{W_n(i,j;dz)}{\nu(dz)}\,W_n(i,j;dz)
\\&=\frac{2}{n^2} \sum_{1\leq i<j\leq n}\mathcal{H}(W_n(i,j)\mid \nu)=\int_{[0,1]^2}\mathcal{H}(W_n(x,y)\mid \nu)\, dxdy=I_\nu(W_n)
\end{split}
\]
Now the proof is concluded by showing that 
\[
\lim_{n\to\infty}|I_\nu(W_n)-I_\nu(W)|=0.
\]
This follows directly from \ref{cor:ConvEntropyAproxim} and as $\widehat{W}_n$ and $W_n$ differ only on a set of Lebesgue measure $\tfrac{1}{n}$ and $I_\nu(W)$ is uniformly bounded on the probability simplex $\cP(\Space)$ when $\Space$ is finite and $\nu$ is fully supported on $\Space.$
\end{proof}

\begin{lemma}\label{lemm:auxBall1}
 For any $\varepsilon \in (0,1)$ and let $W_{h_{n,M}}$ be as above, the probability graphon obtained by sampling from $W_n.$ The following equality holds: 
 \begin{equation*}
     \lim_{n \to \infty} \mathbb{P}(d_{\Box}(W_{h_{n,M}}, W_n) > \varepsilon) =0.
 \end{equation*}
\end{lemma}
\proof
The proof follows directly combining \cite[Lemma 6.9]{abraham2023probabilitygraphons} and \cite[Proposition 1.2]{abraham2023probabilitygraphons}.
\endproof
We can prove the following
\begin{lemma}\label{lem:ballto1}
For any $\varepsilon \in (0,1)$, 
\[
\lim_{n \to \infty} \mathbb{P}_{n,W}(B(W_n, \varepsilon)) = 1.
\]
\end{lemma}

\begin{proof}
Fix $\varepsilon \in (0,1)$. Let $W_{h_{n,M}}$ as above, the probability graphon obtained by sampling from $W_n$. We have
\[
1 - \mathbb{P}_{n,W}(B(W_n,\varepsilon)) = \mathbb{P}(d_{\Box}(W_{h_{n,M}}, W_n) > \varepsilon) .
\]
Therefore, by Lemma \ref{lemm:auxBall1} we obtain the result.
\end{proof}

We are ready to prove our main result:

\begin{proof}[Proof of Theorem \ref{thm:lowerbound}]
The proof will take advantage of the reductions we showed we are allowed to do in this section: to complete the argument we just need to perform a standart tilting of the measure. We have that 
\[
\begin{split}
    \mu_{n,\nu}(B(W_n,\varepsilon)) &= \int_{B(W_n,\varepsilon)} d\mu_{n,\nu} 
= \int_{B(W_n,\varepsilon)} \exp\left( -\log \frac{d\mathbb{P}_{n,W}}{d\mu_{n,\nu}} \right) d\mathbb{P}_{n,W}\\
&= \mathbb{P}_{n,W}(B(W_n,\varepsilon)) \cdot \frac{1}{\mathbb{P}_{n,W}(B(W_n,\varepsilon))} \int_{B(W_n,\varepsilon)} \exp\left( -\log \frac{d\mathbb{P}_{n,W}}{d\mu_{n,\nu}} \right) d\mathbb{P}_{n,W}.
\end{split}
\]

By Jensen’s inequality we have:
\[
\log \mu_{n,\nu}(B(W_n,\varepsilon)) \geq \log \mathbb{P}_{n,W}(B(W_n,\varepsilon))
- \frac{1}{\mathbb{P}_{n,W}(B(W_n,\varepsilon))} \int_{B(W_n,\varepsilon)} \log \frac{d\mathbb{P}_{n,W}}{d\mu_{n,\nu}} d\mathbb{P}_{n,W}.
\]

Since $\mathbb{P}_{n,W}(B(W_n,\varepsilon)) \to 1$ by Lemma \ref{lem:ballto1}, this implies that
\[
\liminf_{n \to \infty} \frac{2}{n^2} \log \mu_{n,\nu}(B(W_n,\varepsilon)) 
\geq - \lim_{n \to \infty} \frac{2}{n^2} \int \log \frac{d\mathbb{P}_{n,W}}{d\mu_{n,\nu}} d\mathbb{P}_{n,W}.
\]

By Lemma \ref{lem:radonnykodym}, the expression on the right equals $-I_\nu(W)$. This proves by \eqref{eq:finalapprox} our theorem. 
\end{proof}

\section{Extension to compact spaces using Dawson--G\"artner theorem}\label{sec:dawsongartner}

We turn our attention now to the case where $\Space$ is not finite. Indeed in this section, we exploit the Dawson--G\"artner theorem (Theorem \ref{ThmDawsonGart}) to extend our LDP for probability graphons on a finite space $\Space$ (Theorem \ref{thm:LDP}) to the space of probability graphons on a general compact Polish space $\Space$, which constitutes our main Theorem \ref{thm_mainresult}, which we restate for convenience below.  

\MainTheorem*

To prove it we have to introduce a bit of notions in Section \ref{sec:disc_scheme} that we will need in Section \ref{sec:proj_system} to define a projective system and use the powerful result by Dawson and G\"artner.

\begin{remark}
\label{rmk:VertesWeights4}
A similar LDP can be established for the measure $\widetilde{\mu}_{n,\nu,\nu^{\prime}}$, with $\nu \in \mathcal{P}(\Space)$ and $\nu^{\prime} \in \mathcal{P}(\Space^{\prime})$, defined on the space of pairs $(W^{\mathrm{v}},W^{\mathrm{e}})$ as in Remarks \ref{Rmk:Vertex-Weights1}, \ref{rmk:VertexWeights2}, and \ref{rmk:VertexWeightedGraphs3}. This corresponds to sampling random vertex- and edge-weighted graphs, where vertices are drawn independently from $\nu^{\prime}$ and edges independently from $\nu$. In this setting, the rate function is given by the relative entropy
\[
\widetilde{\mathcal{H}}_{\mathrm{v},\mathrm{e}}\big((W^{\mathrm{v}},W^{\mathrm{e}})\,\big\|\,(\nu^{\prime},\nu)\big),
\]
as defined in \eqref{eq:RelEntropyVertexWeight}.
\end{remark}

\subsection{Discretization scheme}\label{sec:disc_scheme}
A discretization scheme partitions a continuous space into finite pieces, in such a way that one can take a suitable inverse (projective) limit of the discrete models to recover results for the continuous space.

Let $\Space$ be a Polish space, and let $(P_m)_{m \in \mathbb{N}}$ denote a sequence of finite partitions of $\Space$ into non-empty sets. For each $m \in \mathbb{N}$, we write
\[P_m = \{A_{m,i} : i = 1, \ldots, |P_m|\,\}.\]
The sequence $(P_m)_{m \in \mathbb{N}}$ is said to be \emph{nested} if, whenever $m \leq n$ and $A_{m,i} \in P_m$, there exists an index set $J \subset \{1, \ldots, |P_n|\}$ such that
\[A_{m,i} = \bigcup_{j \in J} A_{n,j}.\]

For any subset $A \subset \Space$, we define its diameter by
\[\operatorname{diam}(A) := \sup\{d(x,y) : x,y \in A\}\]
where $d$ is a chosen metric of $\Space$.
Finally, given a Borel measure $\mu$ on $\Space$, a set $A \subset \Space$ is called a continuity set of $\mu$ if $\mu(\partial A) = 0$. We have the following lemma, adapted to our case:

\begin{lemma}[Lemma 7.1 in \cite{bollobas2007phase}, Lemma 5.3 in \cite{andreis2023large}]\label{Lemm:AndreisApprox}
There exists a sequence of finite partitions $ (P_m)_{m \in \mathbb{N}}$ of $\Space$ with the following properties:
\begin{enumerate}
    \item For any $m \in \mathbb{N}$ and any $ i = 1, \ldots, |P_m| $, we have that the set $ A_{m,i} \in P_m $ is measurable and a continuity set of $ \mu.$ 
    
    \item The sequence $ (P_m)_{m \in \mathbb{N}} $ is nested; that is, for all \( m \), each set in $ P_{m+1} $ is contained in some set of $ P_m $.
    
    \item It holds that
    \begin{equation}
        \lim_{m \to \infty} \max_{1 \leq i \leq |P_m|} \operatorname{diam}(A_{m,i}) = 0.
    \end{equation}
\end{enumerate}
\end{lemma}

In the following, we always assume that $(P_m)_{m \in \mathbb{N}}$ has all the properties given in Lemma~\ref{Lemm:AndreisApprox}. For $m \in \mathbb{N}$ and any $A_{m,i} \in P_m$ we pick exactly one point $x_{m,i}$ from the set $A_{m,i}$, which we call the \emph{representative} of $A_{m,i}$. 
\begin{remark}
    The representative point may be chosen arbitrarily. For intuition, one may picture a subset of $\R$ partitioned into intervals, and within each interval select the representative according to a fixed rule—such as the rightmost point, the midpoint, or the median. 
\end{remark}
We define
\begin{equation}\label{eq:DiscTypeSpace}
\Space_m := \{x_{m,i} : i = 1, \ldots, |P_m|\}.
\end{equation}
and define the projection
\begin{equation}
\pi_m : \Space \to \Space_m, \quad x \mapsto x_{m,i}
\end{equation}
where $i$ is such that $A_{m,i}$ is the unique set in $P_m$ containing $x$.

The map $\pi_m$ naturally defines a projection in the space of probability measures on $\Space$, $\mathcal{P}(\Space)$. Indeed, for any $m \in \mathbb{N}$ we can write
\[
\pi_m : \mathcal{P}(\Space) \to \mathcal{P}(\Space_m), \quad \pi_m(\mu) := \mu \circ \pi_m^{-1}=\sum^{|P_m|}_{i=1} \mu(A_{i,m})\delta_{x_{m,i}}\quad \text{for } \mu \in \mathcal{P}(\Space),
\]
where, for the sake of notation, we used the same symbol for the projection map on $\Space$ and $\mathcal{P}(\Space)$.
Once $\pi_m$ has been defined on $\mathcal{P}(\Space)$, we can lift the projection further to the space of probability graphons over $\Space,$ $\Graphon(\Space)$. Using again the same symbol, we can write for any $m \in \mathbb{N}$ 
\[
\pi_m : \Graphon(\Space) \to \Graphon(\Space_m), \quad \pi_m(W) (x,y):= \pi_m(W (x,y))\quad \text{for } W \in \Graphon(\Space)).
\]
This automatically defines a map also on the quotient space of probability graphons $\UGraphon(\Space)$,
\[
\pi_m : \UGraphon(\Space) \to \UGraphon(\Space_m), \quad \pi_m(\widetilde{W}) (x,y):= \pi_m(W )(x,y)\quad \text{for } \widetilde{W} \in \Graphon(\Space)),
\]
where $W$ is any representative of the class $\widetilde{W}$. Indeed the above formula is well defined as we show in the next lemma.

\begin{lemma}\label{lemm:RelabelPim}
For $W\in \Graphon(\Space)$, we have the equality $(\pi_{m}(W))^{\varphi}=\pi_{m}(W^{\varphi})$ for every $\varphi\in \Relabel .$ In particular, the map $\pi_m$ is well defined as a map from $\UGraphon(\Space)$ to $\UGraphon(\Space_m).$
\end{lemma}
\proof
For $W\in \Graphon(\Space_m)$, we have the equality $(\pi_{m}(W))^{\varphi}=\pi_{m}(W^{\varphi})$ for every $\varphi\in \Relabel $ as $(\pi_{m}(W))^{\varphi}(x,y)=\pi_{m}(\widetilde{W})(\varphi(x),\varphi(y))=\pi_{m}(\widetilde{W}(\varphi(x),\varphi(y)))=\pi_{m}(\widetilde{W}^{\varphi}(x,y))=\pi_{m}(\widetilde{W}^{\varphi})(x,y)$ for a.e.\ $(x,y)\in [0,1]^2.$ Therefore, the map $\pi_m$ is well defined.
\endproof
All the previous maps are well defined and which $\pi_m$ we will be using is always clear by the context.

Analogously to the projection $\pi_m : \Space \to \Space_m$, we can also define projections between partitions of different resolutions. For any $m,n \in \mathbb{N}$ with $m \leq n$, set
\begin{equation}
\pi_{m,n} : \Space_n \to \Space_m, \quad x_{n,j} \mapsto x_{m,i}
\end{equation}
where $x_{m,i}$ denotes the representative of $A_{m,i}$, and $A_{m,i} \in P_m$ is the unique element containing $x_{n,j}$. This map is well defined for all suitable $m$ and $n$ and enjoys the same properties of $\pi_m$.

In a similar fashion, we can extend the projection $\pi_{m,n}$ to measures, graphons, and unlabeled graphons, and use the same symbol for all cases (the intended meaning being clear from the context):
\begin{align}
&\pi_{m,n} : \mathcal{P}(\Space_n) \to \mathcal{P}(\Space_m), \quad &&\pi_{m,n}(\mu) := \mu \circ \pi_{m,n}^{-1},\\
&\pi_{m,n} : \Graphon(\Space_n) \to \Graphon(\Space_m), \quad &&(\pi_{m,n}W)(x,y) := \pi_{m,n}(W(x,y)), \quad W \in \Graphon(\Space),\\
&\pi_{m,n} : \UGraphon(\Space_n) \to \UGraphon(\Space_m), \quad &&(\pi_{m,n}\widetilde{W})(x,y) := (\pi_{m,n}W)(x,y), \quad \widetilde{W} \in \UGraphon(\Space),
\end{align}
where in the last line $W$ is any representative of the equivalence class $\widetilde{W}$. This latter definition is well posed, as stated in the following lemma.

\begin{lemma}\label{lemm:RelabelPiMandn}
For $W\in \Graphon(\Space_n)$, we have the equality $(\pi_{m,n}(W))^{\varphi}=\pi_{m,n}(W^{\varphi})$ a.e. in $[0,1]^2$ for every $\varphi\in \Relabel .$ In particular, the map $\pi_{m,n}$ is well defined as a map from $\UGraphon(\Space_n)$ to $\UGraphon(\Space_m).$
\end{lemma}
\proof
For $W\in \Graphon(\Space_n)$, we have the equality $(\pi_{m,n}(W))^{\varphi}=\pi_{m,n}(W^{\varphi})$ for every $\varphi\in \Relabel $ as 
\[
\begin{split}
 (\pi_{m,n}(W))^{\varphi}(x,y)&=\pi_{m,n}(W)(\varphi(x),\varphi(y))=\pi_{m,n}(W(\varphi(x),\varphi(y)))\\&=\pi_{m,n}(W^{\varphi}(x,y))=\pi_{m,n}(W^{\varphi})(x,y)  
\end{split}
\] 
for a.e.\ $(x,y)\in [0,1]^2.$ Therefore, the map $\pi_{m,n}$ is well defined.
\endproof
\begin{remark}\label{rmk:proj_subprob}
    In general the maps $\pi_m$ and $\pi_{m,n}$ are well defined for any element of $\mathcal{M}_{\leq1}(\Space)$. This will be used in the next section.
\end{remark}
\subsection{Projective system and main result}\label{sec:proj_system}
Building on the discretization scheme introduced in the previous subsection, we will now construct a projective system that falls within the scope of the Dawson--Gärtner theorem (Theorem \ref{ThmDawsonGart}) in large deviations theory. This framework will enable us to lift the large deviation principle for probability graphons on a finite space $\Space$ (Theorem \ref{thm:LDP}) to the setting of probability graphons on a general compact Polish space $\Space$.

Recall the notions from the preceding subsection: in order to apply the projections introduced above, we require them to form a \emph{projective system}. The relevant notion is as follows.

\begin{definition}\label{def:projsyst}
    A family $(L_m, \pi_{m,n})_{m \leq n}$ is called a \emph{projective system} if
\begin{enumerate}
\renewcommand{\labelenumi}{\roman{enumi})}
    \item for any $m \in \mathbb{N}$, the space $L_m$ is a Hausdorff topological space,
    \item for any $m, n \in \mathbb{N}$ with $m \leq n$, the mapping $\pi_{m,n} : L_n \to L_m$ is continuous,
    \item for any $m, n, p \in \mathbb{N}$ with $m \leq n \leq p$, we have $\pi_{m,p} = \pi_{m,n} \circ \pi_{n,p}$. 
\end{enumerate}
\end{definition}
As we already said, the aim of a projective system is to approximate some larger space in a meaningful way. We expect that as the partition gets finer, the approximation gets better, till recovering the main features of the space in the limit. This is achieved through the \emph{projective limit} of the projective system $(L_m, \pi_{m,n})_{m \leq n}$. This is by definition the subset of the product space $\prod_{m \in \mathbb{N}} L_m$ that contains all elements $(\lambda_m)\in L_m$ with $m \in \mathbb{N}$ satisfying the following relation
\[
\pi_{m,n}(\lambda_n) = \lambda_m \quad \text{for all } m, n \in \mathbb{N} \text{ with } m \leq n,
\]
i.e.\ the space containing all the elements consistent with the maps $\pi_{m,n}$.
We denote the projective limit by
\[
\varprojlim L_m.
\]

The projective limit comes with a natural topology induced by the product topology on $\prod_{m \in \mathbb{N}} L_m$, called the \emph{projective limit topology}. 
In particular, a sequence $\lambda^{(n)} \in \varprojlim L_m$ converges to $\lambda \in \varprojlim L_m$ as $n \to \infty$ if and only if for every $m \in \mathbb{N}$ we have $\lambda^{(n)}_m \to \lambda_m$ as $n \to \infty$.

Projective limits are fundamental for applying \DG theorem. We will therefore first prove that
\[
(\UGraphon(\Space_m), \pi_{m,n})_{m \leq n}
\]
forms a projective system, then we will prove that its projective limits exists and we will identify it. To this end, we first introduce some auxiliary notation and establish a few preparatory lemmas.

\begin{lemma}
Given the family $P_m$, fix $s\in \N$. Let $\mu$ be any element of $\mathcal{M}_{\leq 1}(\Space_s).$ Let $(\mu_n)_n$ be a sequence of measures such that $\mu_n\in \mathcal{M}_{\leq 1}(\Space_s)$ for all $n\in \mathbb{N}.$ Let's assume that $\mu_n$ converges weakly to $\mu$. Then $\pi_{m,s}(\mu_n)$ converges weakly to $\pi_{m,s}(\mu)$ for all $m\in \mathbb{N}.$ 
\end{lemma}
\proof
Recall Remark \ref{rmk:proj_subprob}. Since $\Space_s$ carries the discrete topology, the projection map $\pi_{m,s} : \Space_s \to \Space_m$ is continuous for every $m \in \mathbb{N}$, $m\leq s$. 
Let's consider now the map $\pi_{m,s} : \mathcal{M}_{\leq 1}(\Space_s) \to \mathcal{M}_{\leq 1}(\Space_m)$ (see Remark \ref{rmk:proj_subprob}). We claim that $\pi_{m,s}$ is continuous with respect to weak convergence. Given some $\mu \in \mathcal{M}_{\leq 1}(\Space_s)$ and a sequence $(\mu_n)_{n \in \mathbb{N}} $ in $\mathcal{M}_{\leq 1}(\Space_s)$ such that $\mu_n \to \mu$ weakly as $n \to \infty$, and some (continuous bounded) function $f : \Space_m \to \mathbb{R}$, we clearly have that
\[
\int_{\Space_{s}} f \, d\pi_{m,s}(\mu_n) = \int_{\Space} f \circ \pi_{m,s} \, d\mu_n \to \int_{\Space} f \circ \pi_{m,s} \, d\mu = \int_{\Space_{s}} f \, d\pi_{m,s}(\mu),
\]
as $n \to \infty$, since $f \circ \pi_{m,s}$ is continuous and bounded. Hence, $\pi_{m,s}(\mu_n) \to \pi_{m,s}(\mu)$ weakly as $n \to \infty$.
\endproof

Therefore, we directly obtain the following corollary.
\begin{corollary}\label{cor:ineqLevyProkProjSys}
For every $\varepsilon>0$ there exists $\delta>0$ such that for all $\mu,\mu'\in \mathcal{M}_{\leq1}(\Space_s)$ if $d_{\mathcal{LP}}(\mu,\mu')<\delta$ then $d_{\mathcal{LP}}(\pi_{m,s}(\mu),\pi_{m,s}(\mu'))<\varepsilon.$ Moreover, $\pi_{m,s}$ is uniformly continuous on $\mathcal{M}_{\leq 1}(\Space_s)$.
\end{corollary}

\begin{lemma}\label{lemm:equal_pimProjSys}
Let $R$ be a measurable subset in $[0,1]^2$ and let $W\in \Graphon(\Space_n).$ The following equality holds
    \begin{equation*}
        \pi_{m,n}\left(\int_{R} W(x,y;\cdot) \mathrm{d}x \mathrm{d}y\right)=  \int_{R} \pi_{m,n}\left(W\right)(x,y;\cdot) \mathrm{d}x \mathrm{d}y.
    \end{equation*}
\end{lemma}
\proof
Let $A$ be any subset of $\Space_{n}.$ The lemma follows from the chain of equalities below
\begin{equation*}
\begin{aligned}
    \pi_{m,n}\left(\int_{R} W(x,y;\cdot) \mathrm{d}x \mathrm{d}y\right)(A)&=\left(\int_{R} W(x,y;\cdot) \mathrm{d}x \mathrm{d}y\right)(\pi^{-1}_{m,n}(A))\\&=\left(\int_{R} W(x,y;\pi^{-1}_{m,n}(A)) \mathrm{d}x \mathrm{d}y\right)\\
    &=\left(\int_{R} \pi_{m,n}\left(W(x,y;\cdot)\right)(A) \mathrm{d}x \mathrm{d}y\right)
    \\
    &=\left(\int_{R} \pi_{m,n}\left(W\right)(x,y;A) \mathrm{d}x \mathrm{d}y\right).
    \end{aligned}
\end{equation*}
\endproof

Now we can finally show that $(\UGraphon(\Space_m), \pi_{m,n})_{m \leq n}$ is a projective system.

\begin{lemma}\label{lemm:W1ProjSyst}
    The system $(\UGraphon(\Space_m), \pi_{m,n})_{m \leq n}$ is a projective system.
\end{lemma}
\proof
Property $i)$ of Definition \ref{def:projsyst} is immediate. Indeed the space $\UGraphon(\Space_m)$ is a metric space, therefore an Hausdorff topological space. 

We now turn to the proof of $ii)$, that the projections $\pi_{m,n}$ are continuous. 
Recall that for any $\varphi\in \Relabel$ we have$(\pi_{m,n}(W))^{\varphi}=\pi_{m,n}(W^{\varphi})$ for a.e.\ $(x,y)\in [0,1]^2$ from Lemma \ref{lemm:RelabelPiMandn}.
For the projections $\pi_{m,n},$ using  Lemma \ref{cor:ineqLevyProkProjSys} and Lemma \ref{lemm:equal_pimProjSys}, we obtain that for $W,U\in \Graphon(\Space)$ and any $\varphi\in \Relabel$ we have that for every $\varepsilon>0$ and all measurable sets $S,T\subset[0,1]$ there exists a $\delta>0$  such that 
\begin{equation*}
    d_\mathcal{LP}\left(\int_{S\times T}  (W(x,y;\cdot))^{\varphi} \mathrm{d}x \mathrm{d}y,\int_{S\times T} U(x,y;\cdot) \mathrm{d}x \mathrm{d}y\right)<\delta
\end{equation*}
implies \begin{equation*}
      d_\mathcal{LP}\left(\int_{S\times T}  (\pi_{m,n}\left(W\right)(x,y;\cdot))^{\varphi} \mathrm{d}x \mathrm{d}y,\int_{S\times T}  \pi_{m,n}\left(U\right)(x,y;\cdot) \mathrm{d}x \mathrm{d}y\right)<\varepsilon,
\end{equation*}
where we used the uniform continuity of $\pi_{m,n}$ from Corollary \ref{cor:ineqLevyProkProjSys}.

Therefore, taking the supremum over measurable subsets $S,T\subset[0,1]$ on both sides we get:
\begin{equation*}
    d_{\square}(W^{\varphi},U)< \delta,
\end{equation*}
which in turn implies
\begin{equation*}
    d_{\square}((\pi_{m,n}(W))^{\varphi},(\pi_{m,n}(U)))=d_{\square}(\pi_{m,n}(W^{\varphi}),\pi_{m,n}(U))<  \varepsilon
\end{equation*}

Now, using Lemma \ref{lemm:RelabelPiMandn} we have that
\begin{equation*}
    \delta_{\square}(\widetilde{W},\widetilde{U})<\delta. 
\end{equation*}
This implies that
\begin{equation*}
      \delta_{\square}(\pi_{m,n}(\widetilde{W}),\pi_{m,n}(\widetilde{U}))< \varepsilon.
\end{equation*}
We therefore proved that the functios $\pi_{m,n}$ are continuous. 

Finally, to prove $iii)$, we have
\begin{equation*}
\begin{aligned}
      (\pi_{m,n} \circ \pi_{n,p})(\widetilde{W})(x,y)&= (\pi_{m,n} \circ \pi_{n,p})(\widetilde{W}(x,y)))\\&=(\pi_{m,n} (\pi_{n,p}(\widetilde{W}(x,y))))\\&=\pi_{m,p}(\widetilde{W}(x,y))\\&=\pi_{m,p}(\widetilde{W})(x,y),
\end{aligned}
  \end{equation*}
i.e. $\pi_{m,n} \circ \pi_{n,p}=\pi_{m,p}$. This concludes the proof.
\endproof

For notation's compactness, we will denote with $\mathcal{L}_{\infty}$ the projective limit of the projective system $(\UGraphon(\Space_m), \pi_{m,n})_{m \leq n}$. We now want to identify this limit with a suitable subset of $\UGraphon(\Space)$. This is an essential step in order to apply the theory. 

In order to this, we will need the following notation and results.

Let's define the following set
\begin{equation}
\mathcal{N}_m := \left\{ \mu \in \mathcal{M}_{\leq 1}(\Space) : \mu(\partial A_{m,i}) = 0 \text{ for all } i = 1, \dots, |P_m| \right\}.
\end{equation}
\begin{lemma}
Let $\mu$ be a probability measure in $\mathcal{N}_m.$ Let $(\mu_n)_n$ be a sequence of probability measures such that $\mu_n\in \mathcal{N}_m$ for all $n\in \mathbb{N}.$ Let's assume that $\mu_n$ converges weakly to $\mu$ then $\pi_m(\mu_n)$ converges weakly to $\pi_m(\mu)$ for all $m\in \mathbb{N}.$
\end{lemma}
\proof
By definition, the projection maps $\pi_m : \Space \to \Space_m$ are continuous on 
\[
\Space \setminus \bigcup_{i=1}^{|P_m|} \partial A_{m,i}.
\]
We now consider the induced map $\pi_m : \mathcal{M}_{\leq 1}(\Space) \to \mathcal{M}_{\leq 1}(\Space_m)$. We claim that $\pi_m$ is continuous on the subset
\[
\mathcal{N}_m := \bigl\{ \nu \in \mathcal{M}_{\leq 1}(\Space) : \nu(\partial A_{m,i}) = 0 \ \text{for all } i=1,\dots, |P_m| \bigr\},
\]
with respect to weak convergence.  

Let $\mu \in \mathcal{N}_m$ and $(\mu_n)_{n\in\mathbb{N}}$ in $\mathcal{N}_m$ be such that $\mu_n \Rightarrow \mu$ weakly. For any bounded continuous $f : \Space_m \to \mathbb{R}$ we have
\[
\int_{\Space_m} f \, d\pi_m(\mu_n) 
= \int_{\Space} f \circ \pi_m \, d\mu_n.
\]
Since $f \circ \pi_m$ is continuous $\mu$-almost everywhere (by the definition of $\mathcal{N}_m$), the weak convergence $\mu_n \Rightarrow \mu$ implies
\[
\int_{\Space} f \circ \pi_m \, d\mu_n \longrightarrow \int_{\Space} f \circ \pi_m \, d\mu 
= \int_{\Space_m} f \, d\pi_m(\mu).
\]
Hence $\pi_m(\mu_n) \Rightarrow \pi_m(\mu)$ weakly, proving the claim.
\endproof

As a consequence, we obtain the following corollary.

\begin{corollary}\label{cor:ineqLevyProkPim}
 For every $\varepsilon>0$ there exists $\delta>0$ such that for all $\mu,\mu'\in \mathcal{N}_m$ if $d_{\mathcal{LP}}(\mu,\mu')<\delta$ then $d_{\mathcal{LP}}(\pi_m(\mu),\pi_m(\mu'))<\varepsilon.$ Moreover, given that $\mathcal{N}_m$ is compact (a consequence of $\Space$ being compact) then $\pi_m$ is uniformly continuous on $\mathcal{N}_m$.
\end{corollary}

We define the spaces of measures

\begin{equation}
    \mathcal{N}=\cap^{\infty}_{m=1}\mathcal{N}_m
\end{equation}

\begin{equation}
    \mathcal{M}=\{\mu\in \mathcal{P}(\Space): \ \mu \ll \nu\}\subset \mathcal{N},
\end{equation}
where $\nu$ is the measure in the statement of Theorem \ref{thm_mainresult}.
We denote the probability graphons taking values for almost every $(x,y)\in [0,1]^2$ in $\mathcal{N}\cap\mathcal{P}(\Space)$ (or $\mathcal{M}$) with $\widetilde{\mathcal{W}}_{\mathcal{N}}$ (or $\widetilde{\mathcal{W}}_{\mathcal{M}}$) respectively. The following proposition is a key step towards our main result.

\begin{proposition}\label{Prop:ProjectiveLimDawsonGart}
 The space of probability graphons $\widetilde{\mathcal{W}}_{\mathcal{N}}\subset\UGraphon(\Space)$ can be identified with the projective limit $\mathcal{L}_{\infty}.$ Moreover, the projective limit topology on $\mathcal{L}_{\infty}$ is equivalent to the topology induced by the unlabeled cut metric $\delta_{\square}$ on $\widetilde{\mathcal{W}}_{\mathcal{N}}.$
\end{proposition}

In order to prove this proposition we will to go through some lemmas.

\begin{lemma}\label{lemm:Pi1}
 Let $\widetilde{W}\in \UGraphon(\Space)$. Then for any $m \in \mathbb{N}$, we have that $\pi_m(\widetilde{W}) \in \UGraphon(\Space_m)$. Furthermore, the sequence $(\pi_m(\widetilde{W}))_{m \in \mathbb{N}}$ is an element of the projective limit $\mathcal{L}_\infty$. Consequently, the operator defined by
 \begin{equation*}
\Pi:\,\widetilde{\mathcal{W}}_{\mathcal{N}}\rightarrow \mathcal{L}_\infty, \quad \widetilde{W}\mapsto  (\pi_m(\widetilde{W}))_{m \in \mathbb{N}}
 \end{equation*}
 is well defined.
\end{lemma}
\proof
Recall that by the definition of probability graphon for a.e.\ $(x,y)\in [0,1]^2$ we have $\widetilde{W}(x,y)\in \mathcal{P}(\Space).$ Moreover, the map $\pi_m$ is a measureable map from $ \mathcal{P}(\Space)$ to $ \mathcal{P}(\Space_m).$ 
Therefore, for any $\widetilde{W}\in \UGraphon(\Space)$ we have that $\pi_m(\widetilde{W})\in \UGraphon(\Space_m)$ as $\pi_m(\widetilde{W})(x,y)=\pi_m(\widetilde{W}(x,y))\in \mathcal{P}(\Space_m)$ for a.e.\ $(x,y)\in [0,1]^2$ and  $\pi_m(\widetilde{W})$ is a measurable function as composition of measurable functions. This automatically implies that $(\pi_m(\widetilde{W}))_{m \in \mathbb{N}}$ is an element of the projective limit $\mathcal{L}_\infty$. The well definedness of $\Pi$ is straightforward.
\endproof

\begin{lemma}\label{lemm:equal_pimInt}
Let $R$ be a measurable subset of $[0,1]^2$ and $W\in \Graphon(\Space).$The following equality holds
    \begin{equation*}
        \pi_m\left(\int_{R} W(x,y;\cdot) \mathrm{d}x \mathrm{d}y\right)=  \int_{R}  \pi_m\left(W\right)(x,y;\cdot) \mathrm{d}x \mathrm{d}y.
    \end{equation*}
\end{lemma}
\proof
Let $A$ be any measurable subset of $\Space.$ The lemma follows from the chain of equalities below
\begin{equation*}
\begin{aligned}
    \pi_m\left(\int_{R} W(x,y;\cdot) \mathrm{d}x \mathrm{d}y\right)(A)&=\left(\int_{R} W(x,y;\cdot) \mathrm{d}x \mathrm{d}y\right)(\pi^{-1}_m(A))\\&=\left(\int_{R} W(x,y;\pi^{-1}_m(A)) \mathrm{d}x \mathrm{d}y\right)\\
    &=\left(\int_{R} \pi_m\left(W(x,y;\cdot)\right)(A) \mathrm{d}x \mathrm{d}y\right)
    \\
    &=\left(\int_{R} \pi_m\left(W\right)(x,y;A) \mathrm{d}x \mathrm{d}y\right).
    \end{aligned}
\end{equation*}
\endproof

\begin{lemma}\label{lemm:contPiProjective}
 For any $m \in \mathbb{N}$, the mapping $\pi_m:\widetilde{\mathcal{W}}_{\mathcal{N}}\rightarrow \UGraphon(\Space_m)$ are continuous with respect to the topologies induced by the respective unlabeled cut metrics $\delta_{\square}.$ Consequently, the mapping $\Pi$ from the previous lemma is continuous.
\end{lemma}
\proof
Recall that for any $\varphi\in \Relabel$ we have $(\pi_m(W))^{\varphi}=\pi_m(W^{\varphi})$ for a.e.\ $(x,y)\in [0,1]^2$ from Lemma \ref{lemm:RelabelPim}.
For the projections $\pi_{m},$ using Corollary \ref{cor:ineqLevyProkPim} and Lemma \ref{lemm:contPiProjective}, we obtain that for $W,U\in \Graphon(\Space)$ and any $\varphi\in \Relabel$ we have that for every $\varepsilon>0$ there exists a $\delta>0$  such that 
\begin{equation*}
    d_\mathcal{LP}\left(\int_{S\times T}  (W(x,y;\cdot))^{\varphi} \mathrm{d}x \mathrm{d}y,\int_{S\times T} U(x,y;\cdot) \mathrm{d}x \mathrm{d}y\right)<\delta
\end{equation*}
implies \begin{equation*}
      d_\mathcal{LP}\left(\int_{S\times T}  (\pi_m\left(W\right)(x,y;\cdot))^{\varphi} \mathrm{d}x \mathrm{d}y,\int_{S\times T}  \pi_m\left(U\right)(x,y;\cdot) \mathrm{d}x \mathrm{d}y\right)<\varepsilon.
\end{equation*}
Therefore, the following equality follows taking the supremum over measurable subsets $S,T\subset[0,1]$ on both sides:

\begin{equation*}
    d_{\square}(W^{\varphi},U)\leq \delta.
\end{equation*}
Similarly to what we said above, this implies that
\begin{equation*}
    d_{\square}((\pi_{m}(W))^{\varphi},(\pi_{m}(U)))=d_{\square}(\pi_{m}(W^{\varphi}),\pi_{m}(U))\leq  \varepsilon.
\end{equation*}

Therefore, taking the infimum over $\varphi\in \Relabel$ on both sides of the previous expression, we have that
\begin{equation*}
    \delta_{\square}(\widetilde{W},\widetilde{U})\leq \delta 
\end{equation*}
implies
\begin{equation*}
      \delta_{\square}(\pi_{m}(\widetilde{W}),\pi_{m}(\widetilde{U}))\leq  \varepsilon.
\end{equation*}
Therefore, the functios $\pi_{m}$ are continuous. 
\endproof

Here is a basic property of the mappings $\pi_m$ and $\pi_{m,n}$, which we will need in the following.

\begin{lemma}\label{lemm:compatProjPi}
    Let $m \leq n$. Then the equality $\pi_m = \pi_{m,n} \circ \pi_n$ holds.
\end{lemma}

\proof For what concerns the maps between partitions, i.e.  $\pi_m:\,\Space \to \Space_m$, the equality $\pi_m = \pi_{m,n} \circ \pi_n$ follows from the fact that $(P_m)_{m \in \mathbb{N}}$ is a sequence of nested partitions. Therefore going from $\Space \to \Space_m$ or from $\Space \to \Space_n\to \Space_m$ is equivalent. When considering $\pi_m$ as a map between probability measure, the identity follows from the fact that $\pi_m$ assigns weights to the representative points, and, again, the partition are nested in an ordered way. For $\pi_m$ acting on probability graphons, the identity follows from the previous two and the definition of probability graphon. 
\endproof

In order to identify $\mathcal{L}_\infty$ with $\widetilde{\mathcal{W}}_{\mathcal{N}}$, we need to find which elements of $\widetilde{\mathcal{W}}_{\mathcal{N}}$ are the inverse of elements of $\mathcal{L}_\infty$ through the inverse projection $\Pi^{-1}$. Furthermore we need to verify that $\Pi^{-1}$ is continuous. This requires, for any element of $\mathcal{L}_\infty$, we can find an element of $\widetilde{\mathcal{W}}_{\mathcal{N}}$. Concretely, given an element $(\widetilde{W}_m)_{m\in\N} \in \mathcal{L}_\infty$ we should be able to associate an element $\widetilde{W}$ such that $\pi_m(\widetilde{W})=\widetilde{W}_m$, for all $m\in\N$.  We say that any measure $\mu$ on some measurable space $X$ and any measurable set $U \subset X$  is \emph{concentrated on} $U$ if $\mu(X \setminus U) = 0$. We will need the following set, for $m \in \mathbb{N}$:
\begin{equation}
    \mathcal{P}_m(\Space) := \{ \mu \in \mathcal{P}(\Space) : \mu \text{ is concentrated on } \Space_m \}.
\end{equation}

By our definition of $\mathcal{P}(\Space_m)$, it is clear that  there exists a natural identification $\iota_m:\mathcal{P}(\Space_m)\to \mathcal{P}_m(\Space)$ of $\mathcal{P}(\Space_m)$ with $\mathcal{P}_m(\Space).$ Indeed given the definition of $\Space_m$ via the representative points of each partition, we have $\Space_m \subset \Space$.

For any $\mu_m \in \mathcal{P}(\Space_m)$, we will denote the corresponding element of $\mathcal{P}_m(\Space)$ by $\bar{\mu}_m$, with $\bar{\mu}_m=\iota_m(\mu_m)$. This defines an analogue of the map $\pi_m$ which we call $\bar{\pi}_m$. This is defined by composing the action of $\pi_m$ with $\iota_m$, and acts as 
\[
\bar{\pi}_m : \mathcal{P}(\Space) \to \mathcal{P}_m(\Space).
\]

The mappings $\bar{\pi}_m$ can be naturally extended to the space of graphons in the following way.

\[
\bar{\pi}_m : \Graphon(\Space) \to \Graphon(\Space) , \quad \bar{\pi}_m(\widetilde{W}) (x,y):= \bar{\pi}_m(\widetilde{W} (x,y))\quad \text{for } \widetilde{W} \in \Graphon(\Space).
\]

\begin{lemma}[Lemma 5.11 in \cite{andreis2023large}]
  \label{ApproxIdentProbMeas}  On $\mathcal{P}(\Space)$, the mapping
\[
\bar{\pi}_m:\,\mathcal{P}(\Space) \to \mathcal{P}(\Space), \quad \nu \mapsto \bar{\pi}_m(\nu),
\]
weakly converges uniformly to the identity as $m \to \infty$.
\end{lemma}

\begin{lemma}\label{ApproxIdentProbGraphons}
On $\Graphon(\Space),$ the mapping 
\[
\bar{\pi}_m:\,\Graphon(\Space) \to \Graphon(\Space), \quad W \mapsto \bar{\pi}_m(W),
\]
converges in labelled cut metric $d_{\square}$ to the identity as $m \to \infty$.
\end{lemma}
\proof
From Proposition 3.13 in \cite{abraham2023probabilitygraphons}, we have that for $U,W\in \Graphon(\Space)$  any measurable subsets $S,T\subset [0,1],$
\[
d_{\mathcal{LP}}\left(\int_{S\times T}U(x,y;\cdot), \int_{S \times T}W(x,y;\cdot)\right) \leq \operatorname*{ess\,sup}_{(x,y)\in S\times T} d_{\mathcal{LP}}(U(x,y;\cdot), W(x,y;\cdot)).
\]
Therefore, we obtain that 
\begin{equation*}
\begin{aligned}
     d_{\square}(U,W)&=\sup_{S,T\subset [0,1]^2}d_{\mathcal{LP}}(\int_{S\times T}U(x,y;\cdot), \int_{S \times T}W(x,y;\cdot))
    \\& 
    \leq \operatorname*{ess\,sup}_{(x,y)\in [0,1]^2} d_{\mathcal{LP}}(U(x,y;\cdot), W(x,y;\cdot)).
\end{aligned}
\end{equation*}
In particular, choosing $U=\bar{\pi}_m(W),$ we have 
\begin{equation*}
\begin{aligned}
     d_{\square}(\bar{\pi}_m(W),W)& 
    \leq \operatorname*{ess\,sup}_{(x,y)\in [0,1]^2} d_{\mathcal{LP}}(\bar{\pi}_m(W)(x,y;\cdot), W(x,y;\cdot))\\
    & \leq \sup_{\nu\in \mathcal{P}(\Space)} d_{\mathcal{LP}}(\bar{\pi}_m(\nu), \nu),
\end{aligned}
\end{equation*}
and from Lemma \ref{ApproxIdentProbMeas}, we have that for every $\varepsilon>0$ there exists an $M>0$ big enough such that for every $\nu \in \mathcal{P}(\Space)$ and $m>M$ we have $ d_{\mathcal{LP}}(\bar{\pi}_m(\nu), \nu)<\varepsilon.$ Therefore, we obtain that the map $ W \mapsto \bar{\pi}_m(W)$ is uniformly convergent to the identity on $\Graphon(\Space)$ in the topology of $d_{\square}.$ 
\endproof

\begin{lemma}\label{lemm:PiInv1}
Let $(\widetilde{W}_m)_m\in \mathcal{L}_{\infty}.$ There exists a unique $\widetilde{W}\in \mathcal{W}_{\mathcal{N}}\subset\UGraphon(\Space)$ such that $\pi_m(\widetilde{W})=\widetilde{W}_m.$
Therefore, the function  $\Pi$ is bijective and the inverse of $\Pi,$ that we denote with $\Pi^{-1},$ is well defined. 
\end{lemma}
\proof
Fix an element of  $\mathcal{L}_{\infty}$ and call it $(\widetilde{W}_m)_{m \in \mathbb{N}} $. We will proceed in the following way.
First we will identify, for any $m \in \mathbb{N}$, the probability graphon $\widetilde{W}_m\in \UGraphon(\Space_m)$ with the element $\overline{W}_m \in \widetilde{\mathcal{W}}_{\mathcal{P}_m(\Space)}$ in an unique way. Afterwards we will prove that the sequence  
$(\overline{W}_m)_{m \in \mathbb{N}}$ has a limit point in $\mathcal{\widetilde{W}}_{\mathcal{N}}$, which we will denote by $\overline{W}$.  
It then suffices to show that $\pi_m(\overline{W}) =\widetilde{W}_m$ holds for any $m \in \mathbb{N}$.
For the sequence ($\overline{W}_m )_m$ in $\widetilde{\mathcal{W}}_{\mathcal{N}}$ such that  $\overline{W}_m $ is in $\widetilde{\mathcal{W}}_{\mathcal{P}_m(\Space)}$ for every $m\in \mathbb{N}$ and from the compactness of $\widetilde{\mathcal{W}}_{\mathcal{N}}$ there exists a subsequence $(\overline{W}_{m_i} )_{m_i}$ converging to an element $\overline{W} \in \widetilde{\mathcal{W}}_{\mathcal{N}}.$ First notice that the space $\widetilde{\mathcal{W}}_{\mathcal{N}}$ is compact. This follows from $\widetilde{\mathcal{W}}_1(\Space)$ being compact and by the fact that $\widetilde{\mathcal{W}}_{\mathcal{N}}$ being closed. Indeed this follows from its definition by the restriction on $\mathcal{N}$. Now we fix $m\in \mathbb{N}$ and our goal is to show that $\pi_m(\overline{W}) =\widetilde{W}_m.$ Together with the consistency of the system and Lemma \ref{lemm:compatProjPi}, we get for $n \geq m$ that
\[
\widetilde{W}_m = \pi_{m,n}(\widetilde{W}_n) = \pi_{m,n}(\pi_n(\overline{W}_n)) = \pi_m(\overline{W}_n).
\]
Choosing a subsequence $(\overline{W}_{n_i})_{i \in \mathbb{N}}$ that converges to $\overline{W}$, we obtain
\[
\widetilde{W}_m = \lim_{i \to \infty} \pi_m(\overline{W}_{n_i}) = \pi_m(\overline{W}),
\]
where we used the continuity of the mapping $\pi_m$ proved in Lemma \ref{lemm:contPiProjective}. This concludes the proof.
\endproof

\begin{lemma}\label{lemm:PiInvCont}
 Let $\widetilde{W}_n$ be a sequence in $\UGraphon(\Space).$ Assume that $\pi_m(\widetilde{W}_n)$ converges to $\pi_m(\widetilde{W})$ in unlabeled cut metric $\delta_{\square}$ for all $m\in \mathbb{N},$ then $\widetilde{W}_n$ converges to $\widetilde{W}$ in unlabeled cut metric $\delta_{\square}.$ Consequentely, the function $\Pi^{-1}$ is continuous. 
\end{lemma}
\proof
In this proof we will use the overlay functionals of a probability graphon $\widetilde{W}\in \UGraphon(\Space)$ against a $\CbFunct-$valued kernel $A$, that we denote with $\mathcal{C}(\widetilde{W},A)$, see Definition \ref{ref:OverlayFun} and Remark \ref{Rmk:OverlayClasses}. In particular, the convergence of the overlay functionals is equivalent to the convergence in unlabeled cut metric $\delta_{\square},$ see Theorem \ref{ThmEquivalenceConvQuotientOverlayVersion2}.

 Let $\widetilde{W}_n$ be a sequence in $\UGraphon(\Space)$ and let's assume now that $\pi_m(\widetilde{W}_n)$ converges to $\pi_m(\widetilde{W})$ in unlabeled cut metric $\delta_{\square}$ for all $m\in \mathbb{N}.$

For any $\CbFunct-$valued kernel $A,$ i.e.\ a (strongly) measurable function $A:[0,1]^2\rightarrow \CbFunct$ (see Definition  \ref{DefCBValuedKernel} for more details) we have that
\begin{equation*}
\begin{aligned}
   & \mathcal{C}(\widetilde{W},A)-  \mathcal{C}(\widetilde{W}_n,A)\\&
    =\mathcal{C}(W,A)-  \mathcal{C}(W_n,A)
    \\&
    =\mathcal{C}(W,A)-\mathcal{C}(\pi_m(W),A)+ \mathcal{C}(\pi_m(W),A) - \mathcal{C}(\pi_m(W_n),A)+ \\
   & \qquad \qquad  \qquad \qquad \qquad  \qquad \mathcal{C}(\pi_m(W_n),A)-  \mathcal{C}(W_n,A)\\
   &
   \leq \left|\mathcal{C}(W,A)-\mathcal{C}(\pi_m(W),A)\right|+ \left|\mathcal{C}(\pi_m(\widetilde{W}),A) - \mathcal{C}(\pi_m(\widetilde{W}_n),A)\right|+ \\
   & \qquad \qquad  \qquad \qquad \qquad  \qquad \left|\mathcal{C}(\pi_m(W_n),A))-  \mathcal{C}(W_n,A)\right|,
\end{aligned}
\end{equation*}
where $W_n$ and $W$ are elements of the class of $\widetilde{W}_n$ and $\widetilde{W}$ respectively.
The first and third terms of the last inequality converge to zero as $m$ goes to infinity, using Lemma \ref{ApproxIdentProbGraphons} (the convergence is uniform as the overlay function is a continuous function on a compact space and therefore uniformly continuous) and as from convergence in cut metric $d_{\square}$ (that implies convergence in unlabeled cut metric $\delta_{\square}$) implies convergence of the overlay functionals, Lemma \ref{CorSequencesEqOv}. 

Moreover, by assumption, the second term of the last inequality converges to zero for every $m$ as $n$ goes to infinity, as convergence in unlabeled cut metric $\delta_{\square}$ of $\pi_m(\widetilde{W}_n)$ to $\pi_m(\widetilde{W})$ implies convergence of the overlay functionals, Theorem \ref{ThmEquivalenceConvQuotientOverlayVersion2}. 

As the argument is symmetric, we obtain that $\mathcal{C}(\widetilde{W}_n,A)$ converges to $\mathcal{C}(\widetilde{W},A)$ for every $\CbFunct-$valued kernel $A.$ However, this is equivalent, Theorem \ref{ThmEquivalenceConvQuotientOverlayVersion2}, to the convergence of $\widetilde{W}_n$ to $\widetilde{W}$ in unlabeled cut metric $\delta_{\square}$ and this concludes the proof of the lemma.
\endproof

We have therefore proved Proposition~\ref{Prop:ProjectiveLimDawsonGart}:

\begin{proof}[Proof of Proposition~\ref{Prop:ProjectiveLimDawsonGart}]
The proof follows directly from Lemma \ref{lemm:Pi1}, Lemma \ref{lemm:contPiProjective}, Lemma \ref{lemm:PiInv1} and Lemma \ref{lemm:PiInvCont} above.
\end{proof}

We can finally apply the Dawson--G\"artner theorem to obtain the large deviation principle on a general compact space $\Space.$ 

\begin{theorem}[Dawson--G\"artner, Theorem 4.6.1 in \cite{dembo2009large}]\label{ThmDawsonGart}
Let $\{\mu_n\}_{n \in\N }$ be a family of probability measures on a topological space $X:=\varprojlim L_m$ which is the \emph{projective limit} of the \emph{projective system}  $(L_m, \pi_{m,n})_{m \leq n}$. Suppose that for any $m \in \N$ the Borel probability measures $\mu_n \circ \pi_m^{-1}$ on $L_m$ satisfy the large deviation principle (LDP) with good rate function $I_m(\cdot)$. 

Then $\{\mu_n\}_{n\in\N }$ satisfies the LDP on $X$ with the good rate function
\[
I(x) = \sup_{m \in \N} \{ I_m(\pi_m(x)) \}, \quad x \in X.
\]
\end{theorem}

From Theorem \ref{ThmDawsonGart}, we directly obtain the LDP for probability graphons on general $\Space$. This follows because we can identify $\mathcal{L}_{\infty}=(\UGraphon(\Space_m), \pi_{m,n})_{m \leq n}$ with the space $\widetilde{\mathcal{W}}_{\mathcal{N}}$ and these two spaces have the same topology (Lemma \ref{lemm:W1ProjSyst}). The rate function is given by

\begin{equation}
    I_\nu(\widetilde{W})=\sup_m \widetilde{\mathcal{H}}(\pi_m(\widetilde{W})| \pi_m(\nu)).
\end{equation}
We are left with identifying the rate function $I_\nu$ above with the relative entropy.

The proof of our main theorem is now concluded by the following lemma.
\begin{lemma}\label{Lemm:EntropyDawsonGart}
 Let $\widetilde{W}$ be a probability graphon $\widetilde{W}\in \widetilde{\mathcal{W}}_{\mathcal{N}}.$ We have the following equality 
    \begin{equation*}
        \sup_m \widetilde{\mathcal{H}}(\pi_m(\widetilde{W})| \pi_m(\nu))=\widetilde{\mathcal{H}}(\widetilde{W}| \nu).
    \end{equation*}
\end{lemma}
\proof
For two probability measures $\omega,\nu\in \mathcal{P}(\Space),$ the sequence $\mathcal{H}(\pi_m(\omega)| \pi_m(\nu))$ is monotone increasing in $m,$ i.e.\ 
\begin{equation*}
  \mathcal{H}(\pi_m(\omega)| \pi_m(\nu))\leq \mathcal{H}(\pi_n(\omega)| \pi_n(\nu))
\end{equation*}
for $n\geq m,$ see proof of Lemma 4.1 in \cite{baldasso2023proof}. Therefore, \begin{equation}\label{eq:IneqRelEntPiM}
\mathcal{H}(\pi_m(\widetilde{W})(x,y)| \pi_m(\nu))\leq \mathcal{H}(\pi_n(\widetilde{W})(x,y)| \pi_n(\nu))
\end{equation} for all $n\geq m$ and all $(x,y)\in [0,1]^2.$ Moreover, from Lemma 4.1 in \cite{baldasso2023proof} we have the following equality
\begin{equation*}
    \sup_m \mathcal{H}(\pi_m(\omega)| \pi_m(\nu))=\lim_{m\rightarrow \infty} \mathcal{H}(\pi_m(\omega)| \pi_m(\nu))=\mathcal{H}(\omega| \nu).
\end{equation*}
Thus, we obtain 
\begin{equation*}
\begin{aligned}
\sup_{m} \int_{[0,1]^2}  \mathcal{H}(\pi_m(\widetilde{W}(x,y))| \pi_m(\nu))  &=\lim_{m\rightarrow\infty} \int_{[0,1]^2}  \mathcal{H}(\pi_m(\widetilde{W}(x,y))| \pi_m(\nu))\\&=
\int_{[0,1]^2} \mathcal{H}(\widetilde{W}(x,y)| \nu)\\& =\widetilde{\mathcal{H}}(\widetilde{W}| \nu),
\end{aligned}
\end{equation*}
where the first equality follows from \eqref{eq:IneqRelEntPiM} and the second equality follows from the Monotone Convergence theorem. This concludes the proof of the lemma.
\endproof

\begin{proof}[Proof of Theorem \ref{thm_mainresult}]
Let $\Space$ be a compact Polish space. From the Dawson--G\"artner theorem (Theorem \ref{ThmDawsonGart}), Proposition \ref{Prop:ProjectiveLimDawsonGart}, Lemma \ref{Lemm:EntropyDawsonGart} and the LDP for finite $\Space$ (Theorem \ref{thm:LDP}) we directly obtain the LDP for the space $(\widetilde{\mathcal{W}}_{\mathcal{N}},\delta_{\square})$ with rate function $\widetilde{\mathcal{H}}(\cdot|\nu).$ As $\widetilde{\mathcal{W}}_{\mathcal{N}}$ is a closed subset of $\UGraphon(\Space)$ using Lemma 4.1.5 in \cite{dembo2009large} we directly obtain the LDP on $(\UGraphon(\Space),\delta_{\square})$ with the rate function $\widetilde{\mathcal{H}}(\cdot|\nu),$ as if $\widetilde{W}$ does not take values in $ \mathcal{N}$ for a.e.\ $(x,y)\in [0,1]^2,$ then  $\widetilde{\mathcal{H}}(\widetilde{W}|\nu)=\infty,$ because probability measures outside $ \mathcal{N}$ are not absolutely continuous w.r.t.\ $\nu.$
\end{proof}

We directly obtain the following corollary.
\begin{corollary}\label{cor:SemicontinuityFinal}
Let $\nu$ be a probability measure on a compact Polish space $\Space$. The relative entropy $\widetilde{\mathcal{H}}(\cdot| \nu)$ is a good rate function on $\UGraphon(\Space)$. In particular, the realtive entropy is lower semicontinuous.
\end{corollary}

\section{Conditional distributions and continous graph parameters.}\label{section:CondDist}
One might think that the principal feature of having a large deviation principle on the space of probability graphons is to have estimates on events the type given in Theorem \ref{thm_mainresult}. Perhaps, as important as estimating these tail probabilities for weighted graphs, the LDP allows us to answer another very important useful question. Once a rare event happened, how does a typical weighted graph does look like? Using the notation of the introduction, let's say we have the random weighted graph $g_{n,\nu}$, the graph  on $n$ vertices sampled from $\mu_{n,\nu}$. Let $W_{g_{n,\nu}}$ be the probability graphon of $g_{n,\nu}$ and let $\widetilde{W}_{g_{n,\nu}}$ its equivalence class. Suppose that $\widetilde{W}_{g_{n,\nu}}\in \widetilde{\mathcal{B}}\subseteq\widetilde{\mathcal{W}}_1$, where $\widetilde{\mathcal{B}}$ is a closed set such that the following holds for its interior $\widetilde{\mathcal{B}}^\circ$:
\begin{equation}\label{eq:closedsubset}
\inf_{\widetilde{W} \in \widetilde{\mathcal{B}}^\circ} I_\nu(\widetilde{\mathcal{B}}\,) = \inf_{\widetilde{W} \in \widetilde{\mathcal{B}}} I_\nu(\widetilde{W}) > 0.
\end{equation}
Let be $\mathscr{M}_{\widetilde{\mathcal{B}}}$ the set of minimizers of $I_\nu$ on $\widetilde{\mathcal{B}}$. Given that the space $\widetilde{\mathcal{W}}_1$ is compact, and given that the functional $I_\nu$ is lower semicontinous we have that $\mathscr{M}_{\widetilde{\mathcal{B}}}$ is non empty. Define
\begin{equation}\label{eq:cutmetricsets}
\delta_{\square}(\widetilde{W}, \mathscr{M}_{\widetilde{\mathcal{B}}}) := \inf_{\widetilde{V} \in \mathscr{M}_{\widetilde{\mathcal{B}}}} \delta_{\square}(\widetilde{W},\widetilde{V}). 
\end{equation}
We will prove our second main Theorem \ref{thm:CondDist}, which we restate for convenience: 

\SecondTheorem*

\begin{proof}
Since $\widetilde{\mathcal{W}}_1(\Space)$ is compact and $\widetilde{\mathcal{B}}$ is a closed subset, $\widetilde{\mathcal{B}}$ is also compact. Since $I_\nu$ is a lower semicontinuous function on $\widetilde{\mathcal{B}}$  and $\widetilde{\mathcal{B}}$ is compact, it must attain its minimum on $\widetilde{\mathcal{B}}$. Thus, $\mathscr{M}_{\widetilde{\mathcal{B}}}$ is non-empty. By the lower semicontinuity of $I_\nu$, $\mathscr{M}_{\widetilde{\mathcal{B}}}$ is closed (and hence compact). Fix $\varepsilon > 0$ and let
\[
\widetilde{\mathcal{B}}_\varepsilon := \{\widetilde{V} \in \widetilde{\mathcal{B}} : \delta_{\square}(\widetilde{V}, \mathscr{M}_{\widetilde{\mathcal{B}}}\,) \geq \varepsilon \}.
\]
Then $\widetilde{\mathcal{B}}_\varepsilon$ is again a closed subset. Observe that
\[
\mathbb{P}(\delta_{\square}(\widetilde{W}_{g_{n,\nu}}, \mathscr{M}_{\widetilde{\mathcal{B}}}\,) \geq \varepsilon \mid \widetilde{W}_{g_{n,\nu}} \in \widetilde{\mathcal{B}}) = \frac{\mathbb{P}(\widetilde{W}_{g_{n,\nu}} \in \widetilde{\mathcal{B}}_\varepsilon)}{\mathbb{P}(\widetilde{W}_{g_{n,\nu}} \in \widetilde{\mathcal{B}})}.
\]

Thus, with
\[
I_1 := \inf_{\widetilde{V} \in \widetilde{\mathcal{B}}} I_\nu(\widetilde{V}), \quad I_2 := \inf_{\widetilde{V} \in \widetilde{\mathcal{B}}_\varepsilon} I_\nu(\widetilde{V}),
\]
Theorem \ref{thm_mainresult} and condition \eqref{eq:closedsubset} give
\[
\limsup_{n \to \infty} \frac{1}{n^2} \log \mathbb{P}(\delta_{\square}(\widetilde{W}_{g_{n,\nu}}, \mathscr{M}_{\widetilde{\mathcal{B}}}\,) \geq \varepsilon \mid \widetilde{W}_{g_{n,\nu}} \in \widetilde{\mathcal{B}}) \leq I_1 - I_2.
\]

The proof will be complete if it is shown that $I_1 < I_2$. Now clearly, $I_1 \leq I_2$. If $I_1 = I_2$, the compactness of $\widetilde{\mathcal{B}}_\varepsilon$ implies that there exists $\widetilde{V} \in \widetilde{\mathcal{B}}_\varepsilon$ satisfying $I_{\nu}(\widetilde{V}\,) = I_2$. However, this means that $\widetilde{V} \in \mathscr{M}_{\widetilde{\mathcal{B}}}$ and hence $\widetilde{\mathcal{B}}_\varepsilon \cap \mathscr{M}_{\widetilde{\mathcal{B}}}\, \neq \emptyset$, which is impossible. 
\end{proof} 

\appendix

\section{General results from large deviations theory}

We collect here some general results about large deviations theory.
\begin{theorem}[Theorem 4.1 in \cite{chatterjee2017large}]\label{ThmGeneralLargeDevChatterjee4.1}
Let $\mathscr{X}$ be a real topological vector space whose topology has the Hausdorff property. Let $\mathscr{X}^*$ be the dual space of $\mathscr{X}$. Let $\mathscr{B}$ be the Borel sigma-algebra of $\mathscr{X}$ and let $\{\mu_n\}_{n \geq 1}$ be a sequence of probability measures on $(\mathscr{X}, \mathscr{B})$. Define the logarithmic moment generating function $\Lambda_n : \mathscr{X}^* \to (-\infty, \infty]$ of $\mu_n$ as
\[
\Lambda_n(\lambda) := \log \int_{\mathscr{X}} e^{\lambda(x)} \, d\mu_n(x).
\]

Let $\{\varepsilon_n\}_{n \geq 1}$ be a sequence of positive real numbers tending to zero. Define a function $\bar{\Lambda} : \mathscr{X}^* \to [-\infty, \infty]$ as
\begin{equation*}
\bar{\Lambda}(\lambda) := \limsup_{n \to \infty} \varepsilon_n \Lambda_n(\lambda/\varepsilon_n).
\end{equation*}
The Fenchel--Legendre transform of $\bar{\Lambda}$ is the function $\bar{\Lambda}^* : \mathscr{X} \to [-\infty, \infty]$ defined as
\begin{equation}
\bar{\Lambda}^*(x) := \sup_{\lambda \in \mathscr{X}^*} \left( \lambda(x) - \bar{\Lambda}(\lambda) \right).
\end{equation}
For any compact set $\Gamma \subseteq \mathscr{X}$,
\[
\limsup_{n \to \infty} \varepsilon_n \log \mu_n(\Gamma) \leq -\inf_{x \in \Gamma} \bar{\Lambda}^*(x).
\]
\end{theorem}

\begin{lemma}[Lemma 4.1 in \cite{chatterjee2017large}]\label{lemma:upperbound_balls}
Let $\{\mu_n\}_{n \geq 1}$ be a sequence of probability measures on a metric space $\mathscr{X}$ (equipped with its Borel sigma-algebra), and $\{\varepsilon_n\}_{n \geq 1}$ be a sequence of positive real numbers tending to zero. Let $B(x, \eta)$ denote the closed ball of radius $\eta$ around a point $x$, and let $\Gamma \subseteq \mathscr{X}$ be a compact set. Suppose that $I : \Gamma \to [-\infty, \infty]$ is a function such that for every $x \in \Gamma$,
\[
\lim_{\eta \to 0} \limsup_{n \to \infty} \varepsilon_n \log \mu_n(B(x, \eta)) \leq -I(x).
\]
Then for any closed set $F \subseteq \Gamma$,
\[
\limsup_{n \to \infty} \varepsilon_n \log \mu_n(F) \leq -\inf_{x \in F} I(x).
\]
\end{lemma}

\section{Some useful lemmas}

To proceed we need the following lemma:
\begin{lemma}\label{lem:BoundLogMomen}
Let $\Space$ be a Polish space. Let $f,g$ two functions from $\Space$ to $\R$ and $X$ a $\Space-$valued random variable. 

We have that
\[|\log\EE e^{f(X)}-\log\EE e^{g(X)}|\leq|f(X)-g(X)|\]
\end{lemma}
\begin{proof}
Let's consider two functions $f,g:\Space \to \mathbb{R}$ and a $\Space-$ valued random variable with distribution $\mu,$ $\mu\in \mathcal{P}(\Space).$ Let's denote:
\[
\bar{r} := \sup_{u\in \mathbb{R}} |f(u) - g(u)|,
\]
where $\sup$ here is the essential supremum with respect to $\mu.$
Let's define
\[
r(u) = f(u) - g(u).
\]
First of all, observe that the exponential function is monotone increasing. Therefore:
\[
e^{r(u)} \leq e^{\bar{r}}.
\]

We now, using the fact that the integral of the exponential is non negative, we have:
\begin{align*}
    &\log \int_{\mathbb{R}} e^{f(u)} \, d\mu(u) - \log \int_{\mathbb{R}} e^{g(u)} \, d\mu(u)=\log \int_{\mathbb{R}} e^{g(u) + r(u)} \, d\mu(u) - \log \int_{\mathbb{R}} e^{g(u)} \, d\mu(u) \leq\\
    &\leq\log \int_{\mathbb{R}} e^{g(u) + \bar{r}} \, d\mu(u) - \log \int_{\mathbb{R}} e^{g(u)} \, d\mu(u)=\log (e^{\bar{r}})= \bar{r} = \sup_{u} |f(u) - g(u)|.
\end{align*}
On the other hand we get the other sign of the inequality by:
\begin{align*}
    &\log \int_{\mathbb{R}} e^{g(u)} \, d\mu(u) - \log \int_{\mathbb{R}} e^{f(u)} \, d\mu(u) =\log \int_{\mathbb{R}} e^{f(u) - r(u)} \, d\mu(u) - \log \int_{\mathbb{R}} e^{f(u)} \, d\mu(u) \leq\\
    &\leq\log \int_{\mathbb{R}} e^{f(u) + \bar{r}} \, d\mu(u) - \log \int_{\mathbb{R}} e^{f(u)} \, d\mu(u) = \bar{r} = \sup_{u} |f(u) - g(u)|
\end{align*}

which proves the statement.
\end{proof}

Another useful lemma that we need in our work is the following.

\begin{lemma}\label{lemm:AppendCutStrongerWeak}
Let $\Space$ finite and $|\Space|=s$. The topology induced by the labelled cut metric $d_{\square}$ is stronger than the weak topology of $L^2([0,1]^2;\R^s)\cong  L^2([0,1]^2;\SignedMeas)$ on the space of probability graphons $\Graphon(\Space).$
\end{lemma}
\proof
Recalling that $L^2([0,1]^2;\R^s)$ is an Hilbert space, the proof is a simple adaptation of \cite[Proposition 3.3]{chatterjee2017large}. 
\endproof

\section{Overlay functionals for probability graphons}

In this section, we introduce the notion of overlay functional for probability graphons that we use in the proof of Lemma \ref{lemm:PiInvCont}. The interested reader can find more details about the notions mentioned in this section in \cite{zucal2024probabilitygraphonsrightconvergence}.

In order to introduce the notion of overlay functional we will need the definition of a $\CbFunct$-valued kernel. 

\begin{definition}[$\CbFunct$-valued kernels]\label{DefCBValuedKernel}
Let $\Space$ be a Polish space. A \emph{$\CbFunct$-valued kernel} is a map $A$ from $[0,1]^2$  to $\CbFunct$, 
such that:
\begin{enumerate}
\item   $A$  is a function belonging to $\CbFunct$ in the $z\in \Space$ coordinate, i.e.\ for (almost) every $(x,y) \in [0,1]^2$, $W(x,y;\cdot)$ belongs to $ \CbFunct$.
\item $A$ is \emph{Bochner measurable} in $(x,y).$
\item $A$ is \emph{bounded}, i.e.\
\begin{equation}
  \label{eq:def:TMCb}
\sup\|A\|_{\infty}=    \sup_{x,y\in [0, 1] }\, \TotalMass{A(x, y; \cdot)} <+\infty 
  \end{equation}where the supremum is meant to be an essential supremum (almost everywhere with respect to the Lebesgue measure) and we recall that $\TotalMass{\cdot}$ denotes the infinity norm for a function in $\CbFunct$.

\end{enumerate}
\end{definition}

\begin{definition}\label{ref:OverlayFun}
For a (finite linear combination of) probability graphons $W$ and a $\CbFunct-$valued kernel $A$, the overlay functional of $W$ with respect to $A,$ is
\begin{equation}
\mathcal{C}( W,A)=\sup_{\varphi\in \Relabel} \int_{[0,1]^2}   \left(\int_{\Space}A(x,y;z)W(\varphi(x), \varphi(y);dz)\right) \mathrm{d}x \mathrm{d}y. 
\end{equation}
\end{definition}
These objects were introduced in \cite[Section 4.1]{zucal2024probabilitygraphonsrightconvergence} as generalizations of overlay functionals for real-graphons \cite{borgs2011convergentAnnals}, see also \cite[Section 12]{LovaszGraphLimits}. Overlay fuctionals generalize important notions from optimization, computer science and statistical physics, see \cite{borgs2011convergentAnnals,LovaszGraphLimits}. 

\begin{remark}\label{Rmk:OverlayClasses}
It is easy to see that the notion of overlay functional it is consistent with the classes of weakly isomorphic probability graphons. Therefore, this definition can be adapted trivially to $\widetilde{W}\in \UGraphon(\Space)$ defining  $\mathcal{C}(\widetilde{W}, A)=\mathcal{C}(W, A),$ where $W$ is a representative of the class $\widetilde{W}.$
\end{remark}

The functional $\mathcal{C}(W, A)$ is sub-additive in the first variable (but not linear), i.e.\ for $ U,W \in \Graphon(\Space)$ and $A$ a function  $\CbFunct-$valued kernels, we have
\begin{equation}\label{subaddOverlay}
    \mathcal{C}( U+ W, A) \leq \mathcal{C}( U, A)+\mathcal{C}( W, A).
\end{equation}

Convergence in the labelled cut metric $d_{\square}$ implies the convergence of the overlay functionals as stated in the following lemma.

\begin{lemma}[Corollary 4.10.1 in \cite{zucal2024probabilitygraphonsrightconvergence}]\label{CorSequencesEqOv}
 Let $W\in\Graphon(\Space)$ and $(W_n)_n$ be a sequence such that $W_n\in \Graphon(\Space).$ If $\delta_{\square}( W_n, W) \rightarrow 0$ as $n \rightarrow \infty$, then for every $\CbFunct$-valued kernel $A$ we have $\mathcal{C}(W_n, A) \rightarrow \mathcal{C}(W, A)$.
\end{lemma}

A very important fact about overlay functionals that we exploit in this work is that they characterize completely convergence in unlabeled cut metric.

\begin{theorem}[Theorem 4.30 in \cite{zucal2024probabilitygraphonsrightconvergence}]\label{ThmEquivalenceConvQuotientOverlayVersion2}
For any sequence $( W_n)$ of probability graphons and $ W$ a probability graphon, the following are equivalent:
\begin{enumerate}

\item  The sequence $( W_n)$ is convergent to $ W$ in the unlabeled cut distance $\delta_{\square};$
\item The overlay functional values $\mathcal{C}( W_n, U)$ converge to $\mathcal{C}( W, U)$ for every $\CbFunct-$valued kernel $U.$
\end{enumerate}
\end{theorem}
 Again, see Theorem 4.27 and Theorem 4.30 in \cite{zucal2024probabilitygraphonsrightconvergence} for more details.  

\medskip
\newpage
\section*{Bibliography}
\bibliographystyle{plain}
\bibliography{graphons}

\end{document}